\documentclass{amsart}

\usepackage{amsfonts}
\usepackage{graphicx}
\usepackage{algorithmic}
\usepackage{cleveref}

\usepackage{amsopn}

\usepackage{scalefnt}
\usepackage{pgfplots}
\usepackage{stmaryrd}
\usepackage{longtable}
\usepackage{txfonts}
\usepackage{geometry}
\usepackage{supertabular}

\usepackage{amssymb}
\usepackage{tikz}
\usepackage{mathrsfs}
\usepackage{booktabs}
\usepackage{multicol}
\usepackage{placeins}
\usepackage{caption} 
\usepackage{subcaption} 
\usepackage{lscape}

\usetikzlibrary{arrows, automata, positioning, calc, decorations.pathreplacing}
\newtheorem{thm}{Theorem}[section]

\newtheorem{lemma}[thm]{Lemma}
\newtheorem{conj}[thm]{Conjecture}

\theoremstyle{definition}
\newtheorem{defn}[thm]{Definition}

\theoremstyle{remark}

\newcommand\figwidth{0.31}

\begin{document}

	\title[Second Minimal Orbits]{Second Minimal Orbits, Sharkovski Ordering and Universality in Chaos}\thanks{Department of Mathematics, Florida Institute of Technology, Melbourne, FL 32901 (abdulla@fit.edu, http://my.fit.edu/$\sim$abdulla/).}\thanks{This work was funded by the NSF grant \#1359074: REU Site on PDEs \& Dynamical Systems}
	\author[U.G.Abdulla]{Ugur G. Abdulla}
	
	\author[R.U.Abdulla]{Rashad U. Abdulla}	
		
	\author[M.U.Abdulla]{Muhammad U. Abdulla}
				
	\author[N.H.Iqbal]{Naveed H. Iqbal}
\begin{abstract}
\label{sec:abstract}
      This paper introduces the notion of second minimal $n$-periodic orbit of the continuous map on the interval according as if $n$ is a successor of the minimal period of the map in Sharkovski ordering. We pursue classification of second minimal $7$-orbits in terms of cyclic permutations and digraphs. It is proved that there are 9 types of second minimal orbits with accuracy up to inverses. The result is applied to the problem on the distribution of periodic windows within the chaotic regime of the bifurcation diagram of the one-parameter family of unimodal maps. It is revealed that by fixing the maximum number of appearances of the periodic windows there is a universal pattern of distribution. In particular, the first appearance of all the orbits is always a minimal orbit, while the second appearance is a second minimal orbit. It is observed that the second appearance of 7-orbit is a second minimal 7-orbit with Type 1 digraph. The reason for the relevance of the Type 1 second minimal orbit is the fact that the topological structure of the unimodal map with single maximum is equivalent to the structure of the Type 1 piecewise monotonic endomorphism associated with the second minimal 7-orbit. Yet another important report of this paper is the revelation of the universal pattern dynamics with respect to increased number of appearances. 
\end{abstract}
\maketitle
	


\section{Introduction and Main Result}
\label{sec:introduction}
Let $f: I\rightarrow I$ be a continuous endomorphism, and $I$ be a non-degenerate interval on the real line. Let $f^n: I \rightarrow I$ be an $n$th iteration of $f$. A point $c\in I$ is called a periodic point of $f$ with period $m$ if $f^m(c)=c$, $f^k(c)\neq c$ for $1\leq k <m$. The set of $m$ distinct points
\[ c, f(c), \cdots , f^{m-1}(c) \]
is called the orbit of $c$, or briefly $m$-orbit or periodic $m$-cycle. In his celebrated paper \cite{sharkovsky64}, Sharkovski discovered a law on the coexistence of periodic orbits of continuous endomorphisms on the real line.

\begin{thm}[Sharkovskii]
	\label{thm:sharkovskii}
	\cite{sharkovsky64} Let the positive integers be totally ordered in the following way:
	\begin{equation}\label{sharkovskiordering}
	1\triangleleft 2\triangleleft 2^{2} \triangleleft 2^{3}\triangleleft\dots\triangleleft 2^{2}\cdot 5\triangleleft2^{2}\cdot 3\triangleleft\dots\triangleleft 2\cdot 5\triangleleft 2\cdot 3 \triangleleft\dots\triangleleft 9\triangleleft 7\triangleleft 5\triangleleft 3.
	\end{equation}

	\noindent If a continuous endomorphism, $f:I\rightarrow I$, has a cycle of period $n$ and $m\triangleleft n$, then $f$ also has a periodic orbit of period $m$.
\end{thm}
This result played a fundamental role in the development of the theory of discrete dynamical systems. A conceptually novel proof was given in \cite{block1979}. Following the standard approach, we characterize each periodic orbit with cyclic permutations and directed graphs of transitions or {\it digraphs}. Consider $m$-orbit:
\[ {\bf B}=\{\beta_{1}<\beta_{2}<\cdots <\beta_{m} \} \]
\begin{defn}
	If $f(\beta_i)=\beta_{s_i}$ for $1\le s_i \le m$, with $i=1,2,...,m$, then ${\bf B}$ is associated with cyclic permutation
	\[\pi=
	\begin{bmatrix}
	1&2&\dots&m\\
	s_1&s_2&\dots&s_m
	\end{bmatrix}
	\]
\end{defn}
\begin{defn}
	Let $\omega$ be the order reversing permutation
	
	\[
	\omega =\begin{bmatrix}
	1 & 2 & \dots & m-1 & m\\
	m & m-1 & \dots & 2 & 1
	\end{bmatrix}
	\]	
	
	Then, given a cyclic permutation $\pi$, it's inverse is obtained as $\pi^{-1} = \omega\circ \pi\circ \omega$.
\end{defn}
In the sequel $<a,b>$ means either $[a,b]$ or $[b,a]$.
\begin{defn}
	Let $J_i=[\beta_i, \beta_{i+1}]$. The digraph of $m$-orbit is a directed graph of transitions with vertices $J_1,J_2,\cdots,J_{m-1}$
	and oriented edges $J_i \rightarrow J_s$ if $J_s \subset \  <f(\beta_i), f(\beta_{i+1})>$.
\end{defn}
\begin{defn}
         The inverse digraph of $m$-orbit is a digraph associated with inverse cyclic permutation $\pi^{-1}$. Equivalently, inverse digraph is obtained from the digraph of $m$-orbit by replacing each $J_i$ with $J_{m-i}$. 
\end{defn}
Proof of the Sharkovskii's theorem significantly uses the concept of {\it minimal orbit}.
\begin{defn}
	$n$-orbit of $f$ is called minimal if $n$ is the minimal period of $f$ in Sharkovski's ordering.
\end{defn}

\definecolor{Red}{rgb}{1.0, 0, 0}
\begin{defn}
	Digraph of the $m$-orbit contains the red edge 
	$J_i{\color{Red} \rightarrow}J_s$ if $J_s=<f(\beta_i), f(\beta_{i+1})>$. 
\end{defn}
The structure of the minimal orbits is well understood \cite{Stefan1977, alseda1984, Block1986, block92}. Minimal odd orbits are called Stefan orbits, due to the following characterization:
\begin{thm}[Stefan]
	\label{thm:stefan}
		\cite{Stefan1977, block92} The digraph of a $m=2k+1$ minimal odd orbit has the unique structure given in Figure \ref{fig:minOddDigraph} up to an inverse. 
		
		\begin{figure}[htpb]
			\centering						
			\resizebox{8cm}{!}{\begin{tikzpicture}[
            > = stealth, 
            shorten > = 1pt, 
            auto,
            node distance = 1.5cm, 
            semithick 
        ]
        
				\node[draw=none,fill=none] (jk2) {$J_{k+2}$};
        \node[draw=none,fill=none] (jk) [below of=jk2] {$J_{k}$};
				\node[draw=none,fill=none] (jk1) [left=1cm of {$(jk2)!0.5!(jk)$}] {$J_{k+1}$};				
        \node[draw=none,fill=none] (jk3) [right of=jk2] {$J_{k+3}$};
        \node[draw=none,fill=none] (jkm1) [below of=jk3] {$J_{k-1}$};
        \node[draw=none,fill=none] (dots1) [right of=jk3] {$\dots$};
        \node[draw=none,fill=none] (dots2) [right of=jkm1] {$\dots$};
        \node[draw=none,fill=none] (j2km1) [right of=dots1] {$J_{4}$};
        \node[draw=none,fill=none] (j3) [right of=dots2] {$J_{3}$};
        \node[draw=none,fill=none] (j2k) [right of=j2km1] {$J_{2k}$};
				\node[draw=none,fill=none] (j2) [below of=j2k] {$J_{2}$};
        \node[draw=none,fill=none] (j1) [right=1cm of {$(j2k)!0.5!(j2)$}] {$J_{1}$};        

				\path[->] (jk1) edge[loop left] node{} (jk1);
        \path[->] (jk1) edge node{} (jk);
        \path[->] (jk) edge[red] node{} (jk2);
        \path[->] (jk2) edge[red] node{} (jkm1);
        \path[->] (jkm1) edge[red] node{} (jk3);
				\path[dashed,->] (jk3) edge[red] node{} (dots2);
        \path[dashed,->] (dots2) edge[red] node{} (dots1);
        \path[dashed,->] (dots1) edge[red] node{}  (j3);
        \path[->] (j3) edge[red] node{}  (j2km1);
        \path[->] (j2km1) edge[red] node{}  (j2);
        \path[->] (j2) edge[red] node{}  (j2k);
				\path[->] (j2k) edge[red] node{}  (j1);
				\draw [<-] (-1.6,-0.2) -- (-1.6, 0.2);
				\draw (-1.6,0.2) arc (180:90:8mm) (-0.8,1) -- (6.5,1) (6.5,1) arc (90:0:8mm) (7.3,0.2) -- (7.3, -0.2) -- cycle;
				\draw [<-] (0,0.4) -- (0,1);
				\draw [<-] (1.5,0.4) -- (1.5,1);
				\draw [dashed, <-] (3,0.4) -- (3,1);
				\draw [<-] (4.5,0.4) -- (4.5,1);
				\draw [<-] (6,0.4) -- (6,1);
			\end{tikzpicture}} 		
		\caption{Digraph of Minimal Odd Orbit}
		\label{fig:minOddDigraph}
	\end{figure}	
\end{thm}
Similar characterization of $2(2k+1)$-orbits ($k>1$) is given in \cite{abdulla2013}.
\begin{thm}
\label{thm:abdulla13}	
\cite{abdulla2013} The digraph of a minimal $2(2k+1)$-orbit ($k>1$) has one of four types up to their inverses (Type I is shown in Figure \ref{fig:min2TimesOdd}).
	\begin{figure}[htpb]
		\centering
		\resizebox{10cm}{!}{%
		\begin{tikzpicture}[
								> = stealth, 
								shorten > = 0.1pt, 
								auto,
								node distance = 1.0cm, 
								semithick 
						]
			\node[draw=none,fill=none] (j2k+1) {$J_{2k+1}$};
			\node[draw=none,fill=none] (dots1) [below of=j2k+1] {$\dots$};
			\node[draw=none,fill=none] (dots2) [below of=dots1] {$\dots$};
			\node[draw=none,fill=none] (dots3) [below of=dots2] {$\dots$};
			\node[draw=none,fill=none] (j2k) [below of=dots3] {$J_{2k}$};
			\node[draw=none,fill=none] (jk+2) [right of=dots2] {$J_{k+2}$};
			\node[draw=none,fill=none] (jk-1) [right of=jk+2] {$J_{k-1}$};
			\node[draw=none,fill=none] (jk+1) [right of=jk-1] {$J_{k+1}$};
			\node[draw=none,fill=none] (jk) [right of=jk+1] {$J_{k}$};
			\node[draw=none,fill=none] (j2k-2) [left of=dots2] {$J_{2k-2}$};
			\node[draw=none,fill=none] (j2) [left of=j2k-2] {$J_{2}$};
			\node[draw=none,fill=none] (j2k-1) [left of=j2] {$J_{2k-1}$};
			\node[draw=none,fill=none] (j1) [left of=j2k-1] {$J_{1}$};
			\node[draw=none,fill=none] (j2k+2) [above of=j1] {$J_{2k+2}$};
			\node[draw=none,fill=none] (j2k+3) [above of=j2] {$J_{2k+3}$};
			\node[draw=none,fill=none] (j3k) [above of=jk-1] {$J_{3k}$};
			\node[draw=none,fill=none] (j3k+1) [above of=jk] {$J_{3k+1}$};
			\node[draw=none,fill=none] (j22k+1-1) [below of=j1] {$J_{2(2k+1)-1}$};
			\node[draw=none,fill=none] (j22k+1-2) [below of=j2] {$J_{2(2k+1)-2}$};
			\node[draw=none,fill=none] (j3k+3) [below of=jk-1] {$J_{3k+3}$};
			\node[draw=none,fill=none] (j3k+2) [below of=jk] {$J_{3k+2}$};

			\path[->] (j2k+1) edge[loop above] node{} (j2k+1);
			\path[->] (j2k+1) edge node{} (j2k+2);
			\path[->] (j2k+1) edge node{} (j2k+3);
			\path[dashed,->] (j2k+1) edge node{} (dots1);
			\path[->] (j2k+1) edge node{} (j3k);
			\path[->] (j2k+1) edge node{} (j3k+1);
			\path[->] (j3k+1) edge[red] node{} (jk+1);
			\path[->] (j3k) edge[red] node{} (jk+2);
			\path[dashed,->] (dots1) edge[red] node{} (j2k-2);
			\path[->] (j2k+3) edge[red] node{} (j2k-1);
			\path[->] (j2k+2) edge[red] node{} (j2k);
			\path[->] (j1) edge[red] node{} (j2k+2);
			\path[->] (j2) edge[red] node{} (j2k+3);
			\path[dashed,->] (dots2) edge[red] node{} (dots1);
			\path[->] (jk-1) edge[red] node{} (j3k);
			\path[->] (jk) edge node{} (j3k+1);
			\path[->] (j2k-1) edge[red] node{} (j22k+1-1);
			\path[->] (j2k-2) edge[red] node{} (j22k+1-2);
			\path[dashed,->] (jk+2) edge[red] node{} (dots3);
			\path[dashed,->] (jk+1) edge[red] node{} (j3k+3);
			\path[->] (jk) edge[bend left] node{} (j3k+2);
			\path[->] (j22k+1-1) edge[red] node{} (j1);
			\path[->] (j22k+1-2) edge[red] node{} (j2);
			\path[dashed,->] (dots3) edge[red] node{} (dots2);
			\path[->] (j3k+3) edge[red] node{} (jk-1);
			\path[->] (j3k+2) edge[red] node{} (jk);
			\path[->] (j2k) edge node{} (j22k+1-1);
			\path[->] (j2k) edge node{} (j22k+1-2);
			\path[->] (j2k) edge node{} (j3k+3);
			\path[->] (j2k) edge node{} (j3k+2);
			\path[dashed,->] (j2k) edge node{} (dots3);
		\end{tikzpicture}
		}
		\caption{Type I digraph of minimal $2(2k+1)$-orbit}
		\label{fig:min2TimesOdd}
	\end{figure}
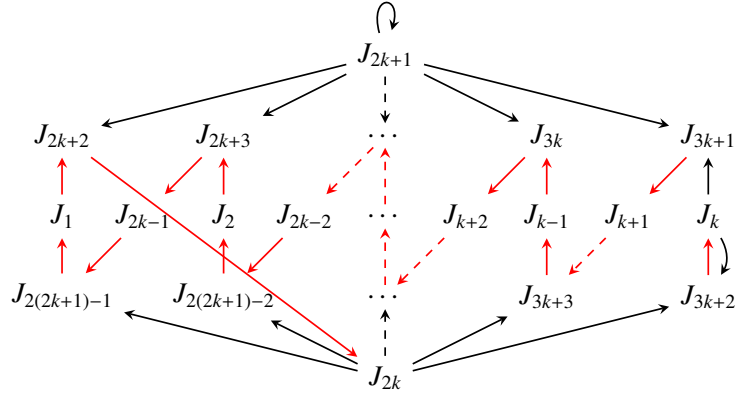
\end{thm}

The main idea of the constructive proof of \cite{abdulla2013} is based on the fact that each half of the minimal $2(2k+1)$-orbit is minimal $2k+1$ orbit of $f^2$. Therefore, the digraph of the minimal $2(2k+1)$-orbit is designed as one of the possible four "unions" of two Stefan digraphs of $f^2$. 
The result of Theorem \ref{thm:abdulla13} can be generalized as follows:

\begin{thm}
	\label{thm:odd2n}
	The digraph of any minimal $2^n(2k+1)$-orbit, $k>1$,  has one of $2^{2^{n+1}-2}$ types up to their inverses. Furthermore, each digraph is strongly simple and can be constructed from the digraphs of two minimal and strongly simple $2^{n-1}(2k+1)$-orbits in $f^2$.
\end{thm}

The main goal of this paper is the characterization of second minimal odd orbits.
\begin{defn} 
	\label{thm:secondmin}
	An $n$-orbit, $n\geq 7$, of $f$ is called second minimal if $n$ is the successor of the minimal orbit of $f$ in the Sharkovskii ordering. 
\end{defn}

For example, if map has a second minimal $7$-orbit, then it has a minimal $5$-orbit, but no $3$-orbit. Our main result reads: 

\begin{thm}
	\label{thm:secondmin7}
	The second minimal $7$-orbit has one of $9$ possible types up to their inverses. The associated cyclic permutations are listed in Table \ref{tab:allsecmin7}; digraphs and piecewise linear representatives are demonstrated in Appendix 1.
	
		\begin{table}%
			\begin{tabular}{ccc}
					\refstepcounter{equation}\label{eq:secmin7_1}
					$\left(\begin{array}{ccccccccccccc}
						1 & 2 & 3 & 4 & 5 & 6 & 7 \\ 
						4 & 5 & 7 & 6 & 3 & 2 & 1 \\
					\end{array} \right)^{1}$ & 
					\refstepcounter{equation}\label{eq:secmin7_2}
					$\left(\begin{array}{ccccccccccccc}
						1 & 2 & 3 & 4 & 5 & 6 & 7 \\ 
						3 & 7 & 5 & 6 & 4 & 2 & 1
					\end{array} \right)^{2}$ &
					\refstepcounter{equation}\label{eq:secmin7_3}
					$\left(\begin{array}{ccccccccccccc}
						1 & 2 & 3 & 4 & 5 & 6 & 7 \\ 
						6 & 4 & 7 & 5 & 3 & 2 & 1 \\
					\end{array} \right)^{3}$ \\
					\refstepcounter{equation}\label{eq:secmin7_4}
					$\left(\begin{array}{ccccccccccccc}
						1 & 2 & 3 & 4 & 5 & 6 & 7 \\ 
						7 & 4 & 6 & 5 & 3 & 1 & 2 \\
					\end{array} \right)^{4}$ &		
					\refstepcounter{equation}\label{eq:secmin7_5}
					$\left(\begin{array}{ccccccccccccc}
						1 & 2 & 3 & 4 & 5 & 6 & 7 \\ 
						4 & 6 & 7 & 5 & 2 & 3 & 1 \\
					\end{array} \right)^{5}$ & 
					\refstepcounter{equation}\label{eq:secmin7_6}
					$\left(\begin{array}{ccccccccccccc}
						1 & 2 & 3 & 4 & 5 & 6 & 7 \\ 
						4 & 6 & 7 & 5 & 3 & 1 & 2 \\
					\end{array} \right)^{6}$ \\
					\refstepcounter{equation}\label{eq:secmin7_7}
					$\left(\begin{array}{ccccccccccccc}
						1 & 2 & 3 & 4 & 5 & 6 & 7 \\ 
						4 & 7 & 6 & 5 & 2 & 1 & 3 \\
					\end{array} \right)^{7}$ &			
					\refstepcounter{equation}\label{eq:secmin7_8}
					$\left(\begin{array}{ccccccccccccc}
						1 & 2 & 3 & 4 & 5 & 6 & 7 \\ 
						3 & 7 & 6 & 5 & 2 & 4 & 1 \\
					\end{array} \right)^{8}$ &
					\refstepcounter{equation}\label{eq:secmin7_9}
					$\left(\begin{array}{ccccccccccccc}
						1 & 2 & 3 & 4 & 5 & 6 & 7 \\ 
						4 & 7 & 5 & 6 & 2 & 3 & 1 \\
					\end{array} \right)^{9}$	\\			
			\end{tabular}
		\caption{All Second Minimal $7$ cycles}
		\label{tab:allsecmin7}
		\end{table}			
\end{thm}

The method of the proof of Theorem \ref{thm:secondmin7} is extended to prove that the second minimal $9$, $11$, and $13$ orbits have respectively $13$, $17$, and $21$ possible types up to their inverse. We conjecture the following result:

\begin{conj}{}
	The digraph of any second minimal $(2k+1)$-orbit, $k\geq 3$, has one of $4k-3$ possible types up to their inverses.
	\label{thm:secondminodd}	
\end{conj}

We adress the proof of the Conjecture \ref{thm:secondminodd} in a forthcoming paper.

The structure of the remainder of the paper is as follows: In Section 2, we recall some preliminary facts. Theorem \ref{thm:secondmin7} is proved in Section 3. In Section 4, we describe a new universal law of the distribution of periodic windows within the chaotic regime of the bifurcation diagram of the one-parameter family of unimodal maps. First we recall the celebrated Feigenbaum scenario of the transition from periodic to chaotic behaviour through successful period doublings and outline the rigorous universality theory in the class of $\mathscr{C}^1$-unimodal maps \cite{Collet1980a, Campanino1981, Collet1980}. In subsection 4.1, we report the numerical reslts which reveal fascinating pattern of distribution of all the superstable periodic orbits when parameter changes in the range between the Feigenbaum transition point to chaos and the value when superstable 3-orbit appears for the first time. In fact, this parameter range is divided into infinitely many Sharkovski $s$-blocks where all the $2^s(2k+1)$-orbits are distributed and the pattern is independent of $s$. Subsection 4.2 demonstartes that the convergence of the successive parameter values for superstable $2^s(2k+1)$-orbits within each $s$-block is exponential with the rate which is independent of the appearance index. Finally, in subsection 4.3, we report the numerical results which demonstrate that any superstable odd orbit in the indicated parameter range is going through successful period doublings according to the Feigenbaum scenario when the parameter decreases to a critical transition point. This indicates that Feigenbaum Universality is true in more general classes of maps, which are the $(2k+1)$st iteration of the class of $\mathscr{C}^1$-unimodal maps. We end Section 4 with the brief outline of the anticipated rigorous universality theory in general classes of maps.       
	
\section{Preliminary Results}
\label{sec:prelim}

\begin{lemma}\label{preliminarylemma}
	The digraph of an $m$-orbit, ${\bf B}=\left\{\beta_{1}<\beta_{2}<\cdots <\beta_{m} \right\}$, $m>2$, possesses the following properties \cite{block92}:
	\begin{enumerate}
		\item The digraph contains a loop: $\exists r_{\ast}$ such that $J_{r_{\ast}}\rightarrow J_{r_{\ast}}$.
		\item $\forall r$, $\exists {r}'$ and ${r}''$ such that $J_{{r}'}\rightarrow J_{r} \rightarrow J_{{r}''}$; moreover, it is always possible to choose ${r}'\neq r$ unless $m$ is even and $r=m/2$, and it is always possible to choose ${r}''\neq r$ unless $m=2$.
		\item If $\left [ {\beta}', {\beta}''\right ]\neq\left [ \beta_{1}, \beta_{m}\right ]$, ${\beta}', {\beta}''\in {\bf B}$, then $\exists J_{{r}'}\subset\left [ {\beta}', {\beta}''\right ]$ and $\exists J_{{r}'}\nsubseteq\left [ {\beta}', {\beta}''\right ]$ such that $J_{{r}'}\rightarrow J_{{r}''}$.
		\item The digraph of a cycle with period $m>2$ contains a subgraph $J_{r_{\ast}}\rightarrow\cdots J_{r}$ for any $1\leq r\leq m-1$.
	\end{enumerate}
\end{lemma}

\begin{defn} 
	A cycle in a digraph is said to be primitive if it does not consist entirely of a cycle of smaller length described several times.
\end{defn}

\begin{lemma}[Straffin]
	\label{thm:straffin}
	\cite{Straffin1978, block92} If $f$ has a periodic point of period $n>1$ and it's associated digraph contains a primitive cycle $J_{0}\rightarrow J_{1}\rightarrow\dots\rightarrow J_{m-1}\rightarrow J_{0}$ of length $m$, then f has a periodic point $y$ of period $m$ such that $f^{k}(y) \in J_{k}, (0\leq k < m)$.
\end{lemma}

\section{Proof of Theorem \ref{thm:secondmin7}}
\label{subsec:secondmindigraphs}

Let $f:I\rightarrow I$ be a continuous endomorphism that has a $7$-orbit which is second minimal. Let $B = \left \{ \beta_{1}< \beta_{2} < \cdots < \beta_{7}  \right \}$ be the ordered elements of this orbit; Let $r_{\ast} = \max\left \{ i\mid f(\beta_{i}) > \beta_{i} \right \}$. Such an $r_{\ast}$ exists since $f(\beta_{1}) > \beta_{1}$ and $f(\beta_{2k+1}) < \beta_{2k+1}$. So, we have a loop: $J_{r_{\ast}}\rightarrow J_{r_{\ast}}$; Let

\[ B^{-} = \left \{ \beta\in B\mid \beta \leq \beta_{r_{\ast}} \right \}, \ \ B^{+} = \left \{ \beta\in B\mid \beta > \beta_{r_{\ast}} \right \}. \] 
Then, $\left | B^{-} \right | + \left | B^{+} \right | = 7$, where $|X|$ denotes the number of elements of the set $X$.  Hence, $\left | B^{-} \right |\neq\left | B^{+} \right |$. Assume that $\left | B^{-} \right | > \left | B^{+} \right |$. Then let $r=\max\left \{ i < r_{\ast}\mid f(\beta_{i})\leq\beta_{r_{\ast}} \right \}$ so $f(\beta_{r})\leq \beta_{r_{\ast}}$; $f(\beta_{r+1})> \beta_{r_{\ast}}$ $\Rightarrow J_{r}\rightarrow J_{r_{\ast}}$. According to Lemma \ref{preliminarylemma} there is a subgraph

\begin{equation}
	\circlearrowright J_{r_{\ast}} \rightarrow \cdot \cdot \cdot \rightarrow J_r \rightarrow J_{r_{\ast}}
\label{eq:secondminfc_orig}
\end{equation}

Assume that this is the shortest path. Then its length is at most $7$, since there are $6$ different intervals, and if any interval is repeated twice, one can get shorter path by removing all the intervals between the repetitions (including one of the repetitions). From another side the length is at least $5$, since if it is $4$ we will deduce by Lemma \ref{thm:straffin} the existence of $3$-orbit. The same conclusion is true if the length is shorter than $4$, since if necessary we can always add $J_{r_{\ast}}$ to the right end of the subgraph \eqref{eq:secondminfc_orig}. Hence, the length can be $5$, $6$, or $7$; In the sequel $<a,b>$ indicates either $[a,b]$ or $[b,a]$; $\begin{matrix}a\\ b \end{matrix}$ or $a\wedge b$ imply either of the entries $a$ or $b$ are valid choices for mappings of a given node;
$J_{r_{i}}\rightarrow [a,b]$ means $f(\beta_{r_{i}}) = a$ and $f(\beta_{r_{i+1}}) = b$.

\begin{enumerate}
	\item [Case 1] length is $7$ $\Rightarrow$ all $6$ intervals are represented in the cycle \eqref{eq:secondminfc_orig}. Choose $r_{1}=r_{\ast}$, $r_{6}=r$ and write
		\begin{equation}
			\circlearrowright J_{r_{1}} \rightarrow J_{r_{2}} \rightarrow J_{r_{3}} \rightarrow J_{r_{4}} \rightarrow J_{r_{5}} \rightarrow J_{r_{6}} \rightarrow J_{r_{1}}
		\label{eq:secondminfc_orig7}		
		\end{equation}		
		
		\noindent Since $\circlearrowright J_{r_{1}} \rightarrow J_{r_{2}}$, but $J_{r_{1}} \not\rightarrow J_{r_{j}}$, $j=3,\cdots, 6\Rightarrow$ $J_{r_{2}}$ must be adjacent to $J_{r_1}$, so either Figure \ref{fig:case11} or Figure \ref{fig:case12} is relevant.
		
		\begin{figure}[ht]%
		\begin{minipage}[t]{0.45\linewidth}
			\centering
				\begin{tikzpicture}
					\draw (1,0)--(3.0,0);
					\foreach \num/\idx in {1,2,3}
						{
						\draw (\num,0.2)--(\num,-0.2); 
						\node (\num) at (\num,0) {};
						}
					\foreach \label/\loc in {J_{r_{2}}/1.5, J_{r_{1}}/2.5}
						{
						\node[above] at (\loc, 0) {$\label$};
						}
					\path[->] (2) edge[bend left=90] node{} (3);
					\path[->] (3) edge[bend left=50] node{} (1);
				\end{tikzpicture}	
				\caption{Case $1.1$}%
				\label{fig:case11}%
			\end{minipage}		
			\hspace{0.1cm}
			\begin{minipage}[t]{0.45\linewidth}
			\centering
			\begin{tikzpicture}
					\draw (1,0)--(3.0,0);
					\foreach \num/\idx in {1,2,3}
						{
						\draw (\num,0.2)--(\num,-0.2); 
						\node (\num) at (\num,0) {};
						}
					\foreach \label/\loc in {J_{r_{1}}/1.5, J_{r_{2}}/2.5}
						{
						\node[above] at (\loc, 0) {$\label$};
						}
					\path[->] (2) edge[bend left=90] node{} (1);
					\path[->] (1) edge[bend left=60] node{} (3);		
			\end{tikzpicture}				
			\caption{Case $1.2$}%
			\label{fig:case12}%
		\end{minipage}
		\end{figure}				
		\noindent Continuing in this manner we get either Figure \ref{fig:case111} or Figure \ref{fig:case122}. Both are Stefan orbits, and the first one is the right one satisfying $\left | B^{-} \right | = 4 > 3 = \left | B^{+} \right |$, while the second one is its inverse satisfying $\left | B^{+} \right | = 4 > 3 = \left | B^{-} \right |$. But Stefan orbit excludes $5$-orbit and so we dismiss this case as irrelevant.
		
		\begin{figure}[ht]%
		\begin{minipage}[t]{0.44\linewidth}
			\centering
				\begin{tikzpicture}
					\draw (1,0)--(7.0,0);
					\foreach \num/\idx in {1,2,...,7}
						{
						\draw (\num,0.2)--(\num,-0.2); 
						\node (\num) at (\num,0) {};
						\node[below] at (\num,-0.21) {\num};
						}
					\foreach \label/\loc in {J_{r_{6}}/1.5, J_{r_{4}}/2.5, J_{r_{2}}/3.5, J_{r_{1}}/4.5, J_{r_{3}}/5.5, J_{r_{5}}/6.5}
						{
						\node[above] at (\loc, 0) {$\label$};
						}

					\path[->] (1) edge[bend left=60] node{} (4);
					\path[->] (2) edge[bend left=60] node{} (7);
					\path[->] (3) edge[bend left=60] node{} (6);
					\path[->] (4) edge[bend left=90] node{} (5);
					\path[->] (5) edge[bend left=65] node{} (3);
					\path[->] (6) edge[bend left=60] node{} (2);		
					\path[->] (7) edge[bend left=60] node{} (1);	
				\end{tikzpicture}	
				\caption{Case $1.1$}%
				\label{fig:case111}%
			\end{minipage}		
			\hspace{0.8cm}
			\begin{minipage}[t]{0.44\linewidth}
			\centering
				\begin{tikzpicture}
					\draw (1,0)--(7.0,0);
					\foreach \num/\idx in {1,2,...,7}
						{
						\draw (\num,0.2)--(\num,-0.2); 
						\node (\num) at (\num,0) {};
						\node[below] at (\num,-0.21) {\num};
						}
					\foreach \label/\loc in {J_{r_{5}}/1.5, J_{r_{3}}/2.5, J_{r_{1}}/3.5, J_{r_{2}}/4.5, J_{r_{4}}/5.5, J_{r_{6}}/6.5}
						{
						\node[above] at (\loc, 0) {$\label$};
						}

					\path[->] (1) edge[bend right=60] node{} (7);
					\path[->] (2) edge[bend right=60] node{} (6);
					\path[->] (3) edge[bend right=65] node{} (5);
					\path[->] (4) edge[bend right=90] node{} (3);
					\path[->] (5) edge[bend right=60] node{} (2);
					\path[->] (6) edge[bend right=60] node{} (1);		
					\path[->] (7) edge[bend right=60] node{} (4);
				\end{tikzpicture}		
			\caption{Case $1.2$}%
			\label{fig:case122}%
		\end{minipage}
		\end{figure}				
		
	\item [Case 2] length is $6$; Choose $r_{1}=r_{\ast}$, $r_{5}=r$ and write
		\begin{equation}
			\circlearrowright J_{r_{1}} \rightarrow J_{r_{2}} \rightarrow J_{r_{3}} \rightarrow J_{r_{4}} \rightarrow J_{r_{5}} \rightarrow J_{r_{1}}
		\label{eq:secondminfc_orig6}		
		\end{equation}	
		
		\noindent We have
		\begin{subequations}
			\begin{align}
				J_{r_{1}}&\rightarrow J_{r_{1}}, J_{r_{1}}\rightarrow J_{r_{2}}, J_{r_{1}}\not\rightarrow J_{r_{3}}, J_{r_{4}}, J_{r_{5}} \\
				J_{r_{2}}&\rightarrow J_{r_{3}}, J_{r_{2}}\not\rightarrow J_{r_{1}}, J_{r_{4}}, J_{r_{5}} \\
				J_{r_{3}}&\rightarrow J_{r_{4}}, J_{r_{3}}\not\rightarrow J_{r_{1}}, J_{r_{5}}\,\,(J_{r_{3}}\rightarrow J_{r_{4}}\,\,\mathrm{optional}) \\
				J_{r_{4}}&\rightarrow J_{r_{5}}, J_{r_{4}}\not\rightarrow J_{r_{1}}, J_{r_{2}}\,\,(J_{r_{4}}\rightarrow J_{r_{5}}\,\,\mathrm{optional}) \\ 
				J_{r_{5}}&\rightarrow J_{r_{1}}, J_{r_{5}}\not\rightarrow J_{r_{3}}\,\,(J_{r_{5}}\rightarrow J_{r_{2}},J_{r_{4}}\,\,\mathrm{optional})
			\end{align}	
			\label{eq:length6rules}
		\end{subequations}
		
		\noindent Hence, we have two possible orders among five intervals $J_{r_{i}}$, $i=1\cdots 5$. Either
		
		\begin{figure}[ht]%
		\begin{minipage}[t]{0.45\linewidth}
			\centering
				\begin{tikzpicture}
					\draw (1,0)--(6.0,0);
					\foreach \num/\idx in {1,2,...,6}
						{
						\draw (\num,0.2)--(\num,-0.2); 
						\node (\num) at (\num,0) {};
						}
					\foreach \label/\loc in {J_{r_{4}}/1.5, J_{r_{2}}/2.5, J_{r_{1}}/3.5, J_{r_{3}}/4.5, J_{r_{5}}/5.5}
						{
						\node[above] at (\loc, 0) {$\label$};
						}
					\path[->] (1) edge[bend left=60] node{} (6);
					\path[->] (2) edge[bend left=60] node{} (5);
					\path[->] (3) edge[bend left=90] node{} (4);
					\path[->] (4) edge[bend left=65] node{} (2);
					\path[->] (5) edge[bend left=60] node{} (1);
					\path[dashed,->] (6) edge[bend left=60] node{} (4);
				\end{tikzpicture}	
				\caption{Case $2.1$, The dashed path demonstrates $J_{r_{5}}\to J_{r_{1}}$}%
				\label{fig:case21}%
			\end{minipage}		
			\hspace{0.1cm}
			\begin{minipage}[t]{0.45\linewidth}
			\centering
			\begin{tikzpicture}
					\draw (1,0)--(6.0,0);
					\foreach \num/\idx in {1,2,...,6}
						{
						\draw (\num,0.2)--(\num,-0.2); 
						\node (\num) at (\num,0) {};
						}
					\foreach \label/\loc in {J_{r_{5}}/1.5, J_{r_{3}}/2.5, J_{r_{1}}/3.5, J_{r_{2}}/4.5, J_{r_{4}}/5.5}
						{
						\node[above] at (\loc, 0) {$\label$};
						}
					\path[dashed,->] (1) edge[bend right=60] node{} (3);					
					\path[->] (2) edge[bend right=60] node{} (6);
					\path[->] (3) edge[bend right=65] node{} (5);
					\path[->] (4) edge[bend right=90] node{} (3);
					\path[->] (5) edge[bend right=60] node{} (2);
					\path[->] (6) edge[bend right=60] node{} (1);		
			\end{tikzpicture}				
			\caption{Case $2.2$, The dashed path demonstrates $J_{r_{5}}\to J_{r_{1}}$}%
			\label{fig:case22}%
		\end{minipage}
		\end{figure}		
		
		\noindent Since there are $6$ different intervals, only one interval is missing. Let us denote this interval $\tilde{J}$, and try to find its place. We have $J_{r_{5}}\rightarrow J_{r_{1}}$ but $J_{r_{5}}\not\rightarrow J_{r_{3}}$. This implies that the missing interval $\tilde{J}$ must be between $J_{r_{1}}$ and $J_{r_{3}}$. Case $2.2$ corresponds to $\left | B^{-} \right | > \left | B^{+} \right |$ so we restrict our discussion to this case.
		
		\begin{equation}
			\left(\begin{array}{ccccccccccccc}
				1 & 2 & 3 & 4_{\ast} & 5 & 6 & 7 \\ 
				\begin{matrix}3\\ 4 \end{matrix} & <5, & 7> & 6 & \begin{matrix}3\\ 4 \end{matrix} & 2 & 1 \\
			\end{array} \right)			
			\label{eq:cycpermlen6}
		\end{equation}			
		
		Hence, constructing a general cyclic permutation from the rules we have the the $2\times 7$ matrix \ref{eq:cycpermlen6}. It follows that either $f(5)=3$ or $f(5)=4$. Now, if $f(5)=3$, according to the rules we must have $f(3)=7$ and $f(1)=4$ and this leads to a valid second minimal $7$ orbit. Alternatively, if $f(5)=4$ then we cannot have $f(2)=5$ else we have a closed $4$ cycle. Thus, $f(2)=7$ and we have a another valid second minimal $7$ orbit. Both of these are displayed in Table \ref{tab:allsecmin7} indexed as $1$ and $2$ respectively and the digraph for the case $f(2)=7$ is presented in Figure \ref{fig:case221}.	
				
		\begin{figure}[ht]%
			\centering						
				\begin{tikzpicture}[
            > = stealth, 
            shorten > = 1pt, 
            auto,
            node distance = 1.5cm, 
            semithick 
        ]
				
					\node[draw=none,fill=none] (j4) {$J_{4}$};
					\node[draw=none,fill=none] (j5) [above right of=j4] {$J_{5}$};
					\node[draw=none,fill=none] (j3) [below right of=j4] {$J_{3}$};				
					\node[draw=none,fill=none] (j6) [right of=j5] {$J_{6}$};
					\node[draw=none,fill=none] (j2) [right of=j3] {$J_{2}$};
					\node[draw=none,fill=none] (j1) [below right of=j6] {$J_{1}$};			

					\path[->] (j4) edge[loop left] node{} (j4);
					\path[->] (j4) edge node{} (j5);
					\path[->] (j5) edge[red] node{} (j3);
					\path[->] (j3) edge[bend left] node{} (j5);
					\path[->] (j5) edge node{} (j2);
					\path[->] (j2) edge[bend left] node{} (j5);
					\path[->] (j2) edge node{} (j6);
					\path[->] (j6) edge[red] node{} (j1);
									
					\draw (0,1.2) arc (180:90:8mm) (0.8,2) -- (2.7,2) (2.7,2) arc (90:0:8mm) (3.5,1.2) -- (3.5, 0.6) -- cycle;
					\draw [<-] (0,0.6) -- (0, 1.2);
					\draw [<-] (1,1.25) -- (1,2);
					\draw [<-] (2.5,1.25) -- (2.5,2);
				\end{tikzpicture}
				\caption{Case $2.2$ Digraph when $f(2)=7$}%
				\label{fig:case221}%
		\end{figure}
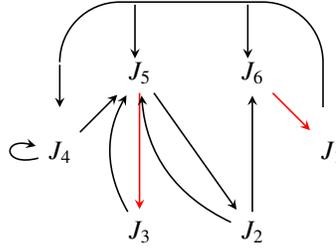			
		
		Finally, in order to show that the orbit above is indeed a valid second minimal odd orbit, it must be proven that there are no odd orbits of length less than $2k-1$. In the case of a second minimal $7$ orbit, it must be proven that no $3$ orbits are present. Assuming that there is an orbit of length $3$, we immediately have two options. Either $J_{1}$ is included in the $3$ orbit, or it isn't included. If $J_{1}$ isn't included, then $J_{6}$ can't be included either, because it only maps to $J_{1}$. $J_{4}$ also can't be included because only $J_{1}$ and $J_{4}$ map to $J_{4}$. If $J_{1}$ isn't in the orbit, and $J_{4}$ is, then $J_{4}$ will only map to itself in the form: $J_{4} \rightarrow J_{4} \rightarrow J_{4}$, as no other orbit will map back to it. The only remaining orbits are $\left \{ J_{2}, J_{3}, J_{5}  \right \}$. Note that these intervals can only form orbits of even length, by splicing together various combinations of the two orbits: $J_{3} \rightarrow J_{5} \rightarrow J_{3}$, and $J_{2} \rightarrow J_{5} \rightarrow J_{2}$. Thus, it is impossible for a $3$ orbit to be present in the above digraph, which does not contain $J_{1}$.
		
		Suppose now, that the assumed $3$ orbit does contain $J_{1}$. $J_{1}$ can map to either: $J_{6}$, $J_{5}$, or $J_{4}$. If $J_{1}$ maps to $J_{6}$, the shortest path back to $J_{1}$ is: $J_{1} \rightarrow J_{6} \rightarrow J_{1}$, which has a length of $2$. If $J_{1}$ maps to $J_{5}$, the shortest path back to $J_{1} $ is $J_{1} \rightarrow J_{5} \rightarrow J_{2} \rightarrow J_{6} \rightarrow J_{1}$, which has a length of $4$. If $J_{1}$ maps to $J_{4}$, the shortest path back to $J_{1}$ is: $J_{1} \rightarrow J_{4} \rightarrow J_{5} \rightarrow J_{2} \rightarrow J_{6} \rightarrow J_{1}$, which has a length of $5$. Thus, it is impossible to form a $3$ orbit, regardless if $J_{1}$ is or isn't contained. If a $3$ orbit is proven impossible, and a $5$ orbit and a $7$ orbit where observed during construction of the orbit, then the above cyclic permutation represents a valid second minimal $7$ orbit.
		
		It is in this way that the validity of the second minimal orbits are proven. Note that, for all cyclic permutations depicted in Table \ref{tab:allsecmin7}, this same method can be used to effectively prove the fact that no $3$ orbits exist. This can also be done by simple observation, as digraphs for a $7$ orbit can only have a finite number of interactions.

	\item [Case 3] length is $5$; (four different intervals are included and two are missing.) Let $r_{1} = r_{\ast}$, $r_{4} = r$, so we have
		\begin{equation}
			\circlearrowright J_{r_{1}} \rightarrow J_{r_{2}} \rightarrow J_{r_{3}} \rightarrow J_{r_{4}} \rightarrow J_{r_{1}}
		\label{eq:chain}
		\end{equation}
		
		\noindent then we have
		\begin{subequations}
			\begin{align}
				J_{r_{1}}&\rightarrow J_{r_{1}}, J_{r_{1}}\rightarrow J_{r_{2}}, J_{r_{1}}\not\rightarrow J_{r_{3}}, J_{r_{4}} \label{eq:l5rulesa} \\
				J_{r_{2}}&\rightarrow J_{r_{3}}, J_{r_{2}}\not\rightarrow J_{r_{1}}, J_{r_{2}}, J_{r_{4}} \label{eq:l5rulesb} \\
				J_{r_{3}}&\rightarrow J_{r_{4}}, J_{r_{3}}\not\rightarrow J_{r_{1}}, J_{r_{3}}\,\,(J_{r_{3}}\rightarrow J_{r_{2}}\,\,\mathrm{optional}) \label{eq:l5rulesc} \\ 
				J_{r_{4}}&\rightarrow J_{r_{1}}, J_{r_{4}}\not\rightarrow J_{r_{2}}, J_{r_{4}}\,\,(J_{r_{4}}\rightarrow J_{r_{3}}\,\,\mathrm{optional}) \label{eq:l5rulesd}
			\end{align}	
			\label{eq:length5rules}
		\end{subequations}
		
		\noindent So we have either Case $3.1$: $J_{r_{4}}J_{r_{2}}J_{r_{1}}J_{r_{3}}$ or Case $3.2$: $J_{r_{3}}J_{r_{1}}J_{r_{2}}J_{r_{4}}$ 

		\begin{figure}[ht]%
		\begin{minipage}[t]{0.45\linewidth}
			\centering
				\begin{tikzpicture}
					\draw (1,0)--(5.0,0);
					\foreach \num/\idx in {1,2,...,5}
						{
						\draw (\num,0.2)--(\num,-0.2); 
						\node (\num) at (\num,0) {};
						\node[below] at (\num,-0.21) {\num};
						}
					\foreach \label/\loc in {J_{r_{4}}/1.5, J_{r_{2}}/2.5, J_{r_{1}}/3.5, J_{r_{3}}/4.5}
						{
						\node[above] at (\loc, 0) {$\label$};
						}
					\path[->] (1) edge[bend left=64] node{} (3);
					\path[->] (2) edge[bend left=60] node{} (5);
					\path[->] (3) edge[bend left=90] node{} (4);
					\path[->] (4) edge[bend left=60] node{} (2);
					\path[->] (5) edge[bend left=60] node{} (1);
				\end{tikzpicture}	
				\caption{Case $3.1$, Where we have $3$ points in $B^{-}$, $2$ in $B^{+}$; $2$ remaining points can't go to $B^{+}$	}%
				\label{fig:case31}%
			\end{minipage}		
			\hspace{0.1cm}
			\begin{minipage}[t]{0.45\linewidth}
			\centering
			\begin{tikzpicture}
				\draw (1,0)--(5.0,0);
				\foreach \num/\idx in {1,2,...,5}
					{
					\draw (\num,0.2)--(\num,-0.2); 
					\node (\num) at (\num,0) {};
					\node[below] at (\num,-0.21) {\num};
					}
				\foreach \label/\loc in {J_{r_{3}}/1.5, J_{r_{1}}/2.5, J_{r_{2}}/3.5, J_{r_{4}}/4.5}
					{
					\node[above] at (\loc, 0) {$\label$};
					}
					\path[->] (1) edge[bend right=60] node{} (5);
					\path[->] (2) edge[bend right=60] node{} (4);
					\path[->] (5) edge[bend right=64] node{} (3);
					\path[->] (4) edge[bend right=60] node{} (1);
					\path[->] (3) edge[bend right=90] node{} (2);				
			\end{tikzpicture}				
			\caption{Case $3.2$, Where we have $2$ points in $B^{-}$, $3$ in $B^{+}$; so we need to add both points to $B^{-}$ and at least one of them should be mapped to $B^{-}$.}%
			\label{fig:case32}%
		\end{minipage}
		\end{figure}		
		
		\noindent Consider Case $3.1$, where should the remaining intervals go (call them $\tilde{J}$, $\hat{J}$)? Case $3.1.1$, assume one is between $J_{r_{1}}$, $J_{r_{3}}$ and the other is between $J_{r_{4}}$, $J_{r_{2}}$. Now, we adjust Figure \ref{fig:case31} in one of the two ways illustrated in Figures \ref{fig:case311} and \ref{fig:case312}.
		
		\begin{figure}[ht]%
		\begin{minipage}[t]{0.44\linewidth}
			\centering
				\begin{tikzpicture}
					\draw (1,0)--(7.0,0);
					\foreach \num/\idx in {1,2,...,7}
						{
						\draw (\num,0.2)--(\num,-0.2); 
						\node (\num) at (\num,0) {};
						\node[below] at (\num,-0.21) {\num};
						}
					\foreach \label/\loc in {J_{r_{4}}/1.5, \hat{J}/2.5, J_{r_{2}}/3.5, J_{r_{1}}/4.5, \tilde{J}/5.5, J_{r_{3}}/6.5}
						{
						\node[above] at (\loc, 0) {$\label$};
						}
					\path[->] (4) edge[bend left=70] node{} (6);
					\path[->] (5) edge[bend left=65] node{} (3);
				\end{tikzpicture}	
				\caption{Case $3.1.1$}%
				\label{fig:case311}%
			\end{minipage}		
			\hspace{0.8cm}
			\begin{minipage}[t]{0.44\linewidth}
			\centering
				\begin{tikzpicture}
					\draw (1,0)--(7.0,0);
					\foreach \num/\idx in {1,2,...,7}
						{
						\draw (\num,0.2)--(\num,-0.2); 
						\node (\num) at (\num,0) {};
						\node[below] at (\num,-0.21) {\num};
						}
					\foreach \label/\loc in {J_{r_{4}}/1.5, \hat{J}/2.5, J_{r_{2}}/3.5, J_{r_{1}}/4.5, \tilde{J}/5.5, J_{r_{3}}/6.5}
						{
						\node[above] at (\loc, 0) {$\label$};
						}
					\path[->] (4) edge[bend left=70] node{} (6);
					\path[->] (5) edge[bend left=65] node{} (2);
				\end{tikzpicture}		
			\caption{Case $3.1.2$}%
			\label{fig:case312}%
		\end{minipage}
		\end{figure}	
				
\noindent Now, we can continue to construct possible orbits graphically in this way and demonstrate which settings for $\hat{J}$ and $\tilde{J}$ result in valid second minimal $7$ cycles however, to better communicate the possible settings, we adopt a slightly different tactic - we will study the cyclic permutations associated with each possible setting in order to extract the relevant cycles. To construct these cyclic permutations first observe in Figures \ref{fig:case31} and \ref{fig:case32} that there are $5$ possible locations in which to insert the extra two intervals $\hat{J}$ and $\tilde{J}$ and we can insert these, assuming we assign $\hat{J}$ first and $\tilde{J}$ second, as demonstrated below in \ref{eq:cases}:

\begin{equation}
	\begin{matrix}
	(1,1) & (1,2) & (1,3) & (1,4) & (1,5) \\ 
	 & (2,2) & (2,3) & (2,4) & (2,5) \\ 
	 &  & (3,3) & (3,4) & (3,5) \\ 
	 &  &  & (4,4) & (4,5) \\ 
	 &  &  &  & (5,5)
	\end{matrix}
	\label{eq:cases}
\end{equation}

\noindent To construct the cyclic permutation determine where each interval $J_{r_{i}}$, $i=1\cdots 4$, is mapped to then combine all the mappings. We begin with Case $3.1$. We furnish an example of how to construct the cyclic permutation for the setting $(2,4)$ which corresponds to the setting for Figures \ref{fig:case311} and \ref{fig:case312}. First, determine where each interval is mapped according to the rules

\begin{align*}
	\begin{matrix}
	J_{r_{1}}=[4, 5]\rightarrow \left [ \begin{matrix}
	5\\ 
	6
	\end{matrix}, \begin{matrix}
	2\\ 
	3
	\end{matrix} \right ] & J_{r_{3}}=[6,7]\rightarrow \left \langle 1 , \begin{matrix}
	2\\
	3\\ 
	4
	\end{matrix}  \right \rangle \\ 
	J_{r_{2}}=[3,4]\rightarrow \left [ 7, \begin{matrix}
	5\\ 
	6
	\end{matrix} \right ] & J_{r_{4}}=[1,2]\rightarrow \left \langle 4, \begin{matrix}
	5\\ 
	6\\
	7
	\end{matrix} \right \rangle
	\end{matrix}
	\label{eq:case24eg}
\end{align*}

\noindent Then, construct the associated cyclic permutation 

\begin{equation}
			\left(\begin{array}{ccccccccccccc}
					1 & 2 & 3 & 4_{\ast} & 5 & 6 & 7 \\ 
					<4, & \begin{matrix}5\\ 6 \end{matrix}> & 7 & \begin{matrix}6\\ 5 \end{matrix} & \begin{matrix}2\\ 3 \end{matrix} & <1, &\begin{matrix}3\\ 2 \end{matrix}> \\
				\end{array} \right)
\label{eq:cycperm24eg}
\end{equation}

		%

Since $f(3)=7\Rightarrow f(7)\neq 3\Rightarrow f(7)=1$ or $2$.
\begin{enumerate}
	\item Case $(2,4)_{1}$: $f(7)=1\Rightarrow$ either $f(6)=2$, $f(5)=3$ or $f(6)=3$, $f(5)=2$
	\item Case $(2,4)_{1,1}$: $f(7)=1$, $f(6)=2$, $f(5)=3$ or
		\begin{equation}
			\left(\begin{array}{ccccccccccccc}
					1 & 2 & 3 & 4 & 5_{\ast} & 6 & 7 \\ 
					<4, & \begin{matrix}5\\ 6 \end{matrix}> & 7 & \begin{matrix}6\\ 5 \end{matrix} & 3 & 2 & 1 \\
				\end{array} \right)
		\label{eq:case2411}	
		\end{equation}	
	\item Case $(2,4)_{1,1,1}$: $f(4)=6\Rightarrow J_{1}{\color{Red} \rightarrow} J_{4}$. Now, $f(1)=5$ implies period $4$-suborbit $\left \{ 1,3,5,7 \right \}$ which is a contradiction. If $f(1)=4$ we get the second minimal $7$ orbit with index $1$ in Table \ref{tab:allsecmin7}.
	\item Case $(2,4)_{1,1,2}$: $f(4)=5\Rightarrow$ either $f(1)=4$, $f(2)=6$ or $f(1)=6$, $f(2)=4$ however since $f(6)=2$ the former implies a period $2$-suborbit $\left \{ 2,6 \right \}$ which is a contradiction. The latter case implies the second minimal $7$ orbit with index $3$ in Table \ref{tab:allsecmin7}.
	\item Case $(2,4)_{1,2}$: $f(7)=1$, $f(6)=3$, $f(5)=2$ or
		\begin{equation}
			\left(\begin{array}{ccccccccccccc}
					1 & 2 & 3 & 4 & 5_{\ast} & 6 & 7 \\ 
					<4, & \begin{matrix}5\\ 6 \end{matrix}> & 7 & \begin{matrix}6\\ 5 \end{matrix} & 2 & 3 & 1 \\
				\end{array} \right)
		\label{eq:case2412}	
		\end{equation}		
		\item Case $(2,4)_{1,2,1}$: $f(4)=6\Rightarrow f(1)=5$, $f(2)=4$. The digraph of the associated cyclic permutation contains the subgraph $J_{2}\rightarrow J_{4}\rightarrow J_{5}\rightarrow J_{2}$ and by Straffin's lemma this implies the existence of a $3$-orbit, a contradiction.
		\item Case $(2,4)_{1,2,2}$] $f(4)=5\Rightarrow$ either $f(1)=4$ or $f(1)=6$. The latter implies a period $4$-suborbit $\left \{ 1,3,6,7 \right \}$, a contradiction. The former implies the second minimal $7$ orbit with index $5$ in Table \ref{tab:allsecmin7}.
		\item Case $(2,4)_{2}$: $f(7)=2\Rightarrow f(6)=1$, $f(5)=3$ or
			\begin{equation}
				\left(\begin{array}{ccccccccccccc}
						1 & 2 & 3 & 4 & 5_{\ast} & 6 & 7 \\ 
						<4, & \begin{matrix}5\\ 6 \end{matrix}> & 7 & \begin{matrix}6\\ 5 \end{matrix} & 3 & 1 & 2 \\
					\end{array} \right)
			\label{eq:case2421}	
			\end{equation}	
\end{enumerate} 	

Considering the alternative we have $f(6)=1\Rightarrow f(1) = 5$ or $f(1)=4$.
\begin{enumerate}		 
	  \item Case $(2,4)_{2,1}$: If $f(1)=5$ the digraph of the associated cyclic permutation contains the subgraph $J_{2}\rightarrow J_{4}\rightarrow J_{5}\rightarrow J_{2}$ which implies the existence of a $3$-orbit, a contradiction.
		\item Case $(2,4)_{2,2}$: $f(1)=4\Rightarrow$ either $f(2)=5$, $f(4)=6$ or $f(2)=6$, $f(4)=5$. In the former case we have a period $4$-suborbit $\left \{ 2,3,5,7 \right \}$, a contradiction. In the latter case we get the second minimal $7$ orbit with index $6$ in Table \ref{tab:allsecmin7}.
\end{enumerate}

\noindent Now, proceeding in this fashion we will analyze each of the $15$ settings to extract valid second minimal $7$ orbits.

\begin{enumerate}
	\item [Setting (1,1)] We have the cyclic permutation
		\begin{equation}
			\left(\begin{array}{ccccccccccccc}
					1 & 2 & 3 & 4 & 5_{\ast} & 6 & 7 \\ 
					\begin{matrix}2\\ 3 \end{matrix} & \begin{matrix}1\\ 3 \end{matrix} & 5 & 7 & 6 & 4 &\begin{matrix}1\\ 2\\ 3 \end{matrix} \\
				\end{array} \right)
		\label{eq:cycperm11}
		\end{equation}	
		
\noindent Observe, letting $f(7)=3$ would force a period $2$-suborbit $\left \{ 1,2 \right \}$ and period $5$-suborbit $\left \{ 3,4,5,6,7 \right \}$ so $f(7) = 1$ or $f(7)=2$ which implies $J_{2}\rightarrow [1,5]$ or $J_{2}\rightarrow [3,5]$ both of which lead to the subgraph $J_{4}\rightarrow J_{6}\rightarrow J_{2}\rightarrow J_{4}$. By Straffin's lemma this implies the existence of a $3$-orbit, a contradiction.

	\item [Setting (1,2)] From \ref{eq:l5rulesa} and \ref{eq:l5rulesb} if follows $f(5)=6$, $f(4)=7$, from \ref{eq:l5rulesa} it follows $J_{r_{4}} = [2,3]\rightarrow <5, 6\wedge 7>$ and hence either $f(2)=6$ or $f(3)=7$ which is a contradiction since three nodes are mapped to $6$ and $7$.
	
	\item [Setting (1,3)] From \ref{eq:l5rulesa} and \ref{eq:l5rulesb} $\Rightarrow f(5)=6$ , $f(6)=3$, and $J_{r_{2}}=[3,4]\rightarrow <6,7>$ which is a contradiction since three nodes $\left \{ 3,4,5\right\}$ are mapped to $\left \{ 6,7\right\}$.
	
	\item [Setting (1,4)] $J_{r_{1}}=[4,5]\rightarrow [5\wedge 6, 3]$, $J_{r_{4}}=[2,3]\rightarrow[4,7]$, and $J_{r_{3}}=[6,7]\rightarrow <3\wedge 4, 1\wedge 2>$; Since $f(5)=3\Rightarrow$ either $f(6)=4$ or $f(7)=4$; but we also have $f(2)=4$, a contradiction.
	
	\item [Setting (1,5)] We have the cyclic permutation
		\begin{equation}
			\left(\begin{array}{ccccccccccccc}
					1 & 2 & 3 & 4_{\ast} & 5 & 6 & 7 \\ 
					\begin{matrix}2\\ 6\\ 7 \end{matrix} & 4 & \begin{matrix}7\\ 6 \end{matrix} & 5 & 3 & \begin{matrix}1\\ 2 \end{matrix} &\begin{matrix}2\\ 1\\ 6 \end{matrix} \\
				\end{array} \right)
		\label{eq:cycperm15}
		\end{equation}		
		
		\begin{enumerate}
			\item Case $(1,5)_{1}$: $f(7)=6\Rightarrow f(3)=7 \Rightarrow f(1)=2$, and $f(6)=1$. The digraph of this cyclic permutation contains the subgraph $J_{6}\rightarrow J_{1}\rightarrow J_{3}\rightarrow J_{6}$ and by Straffin's lemma this implies the existence of a $3$-orbit, a contradiction.
			\item Case $(1,5)_{2}$: $f(7)=1\Rightarrow f(3)=6$ or $f(3)=7$
			\item Case $(1,5)_{2,1}$: $f(7)=1$, $f(3)=6\Rightarrow f(6)=2$, $f(1)=7$ and this implies the $2$-suborbit $\left \{ 1,7\right \}$ and the $5$-suborbit $\left \{ 2,3,4,5,6 \right \}$.
			\item Case $(1,5)_{2,2}$: $f(7)=1$, $f(3)=7\Rightarrow f(6)=2$, $f(1)=6$ and we get the valid second minimal $7$-cycle indexed as $3$ in Table \ref{tab:allsecmin7}.
			\item Case $(1,5)_{3}$: $f(7)=2\Rightarrow f(6)=1$; since $f(1)=6\Rightarrow 2$-suborbit $\left\{ 1,6\right\}$ so we must have $f(1)=7$ and $f(3)=6$ and this implies the valid cyclic permutation indexed as $4$ in Table \ref{tab:allsecmin7}.
		\end{enumerate}
		
	\item [Setting (2,2)] From \ref{eq:l5rulesb} $J_{r_{2}}=[4,5]\rightarrow [7, 6]$ and from \ref{eq:l5rulesd} $J_{r_{4}}=[1,2]\rightarrow <5, 6\wedge 7>$ which is a contradiction since we have $2$ nodes being mapped to $6$ and $7$.
	
	\item [Setting (2,3)] From \ref{eq:l5rulesa} and From \ref{eq:l5rulesb} $J_{r_{1}}=[5,6]\rightarrow [6, 2\wedge 3]$ and $J_{r_{2}}=[3,4]\rightarrow <6,7>$ so $f(5)=6$ and either $f(3)=6$ or $f(4)=6$, a contradiction.
	
	\item [Setting (2,4)] See above.	
		
	\item [Setting (2,5)] We have the cyclic permutation
		\begin{equation}
			\left(\begin{array}{ccccccccccccc}
					1 & 2 & 3 & 4_{\ast} & 5 & 6 & 7 \\ 
					<4, & \begin{matrix}6\\ 7 \end{matrix}> & \begin{matrix}7\\ 6 \end{matrix} & 5 & \begin{matrix}2\\ 3 \end{matrix} & 1 &\begin{matrix}3\\ 2 \end{matrix} \\
				\end{array} \right)
		\label{eq:cycperm25}
		\end{equation}			
		
		\begin{enumerate}
			\item Case $(2,5)_{1}$: If $f(5)=3\Rightarrow f(7)=2\Rightarrow f(2)\neq 7$ or we have period $2$-suborbit $\left\{ 2,7\right\}$. Now, either $f(3)=6$ or $f(3)=7$
			\item Case $(2,5)_{1,1}$: $f(5)=3$, $f(3)=6\Rightarrow f(2)=4,f(1)=7\Rightarrow$ valid second minimal $7$-orbit indexed as $4$ in Table \ref{tab:allsecmin7}.	
			\item Case $(2,5)_{1,2}$: $f(5)=3$, $f(3)=7\Rightarrow J_{r_{4}}=[1,2]\rightarrow <4,6>$.
			\item Case $(2,5)_{1,1,1}$: $f(5)=3$, $f(3)=7$, $f(2)=6\Rightarrow f(1)=4\Rightarrow$ valid second minimal $7$-orbit indexed as $6$ in Table \ref{tab:allsecmin7}.	
			\item Case $(2,5)_{1,1,2}$: $f(5)=3$, $f(3)=7$, $f(2)=4\Rightarrow f(1)=6\Rightarrow$	a period $2$-suborbit $\left\{1,6\right\}$ and a period $5$-suborbit $\left\{2,3,4,5,7\right\}$, a contradiction.		
			\item Case $(2,5)_{2}$: If $f(5)=2\Rightarrow f(7)=3\Rightarrow f(3)\neq 7$ or we get period $2$-suborbit $\left\{3,7\right\}$, so $f(3)=6$ and either $f(1)=4$ or $7$.
			\item Case $(2,5)_{2,1}$: $f(5)=2$, $f(7)=3$, $f(3)=6$, $f(1)=4\Rightarrow$ a valid second minimal cycle indexed as $7$ in Table \ref{tab:allsecmin7}.
			\item Case $(2,5)_{2,2}$: $f(5)=2$, $f(7)=3$, $f(3)=6$, $f(1)=7\Rightarrow$ period $3$-suborbit $\left\{ 2,4,5 \right\}$, a contradiction.
		\end{enumerate}
		
	
	\item [Setting (3,3)] From \ref{eq:l5rulesa}, \ref{eq:l5rulesb} $J_{r_{1}}=[5,6]\rightarrow [6,2]$ and $J_{r_{2}}=[2,3]\rightarrow <6,7>$ so $f(5)=6$ and either $f(2)=6$ or $f(3)=6$, a contradiction.	
	
	\item [Setting (3,4)] We have the cyclic permutation
\begin{equation}
	\left(\begin{array}{ccccccccccccc}
			1 & 2 & 3 & 4_{\ast} & 5 & 6 & 7 \\ 
			\begin{matrix}3\\ 4 \end{matrix} & <\begin{matrix}6\\ 5 \end{matrix}, & 7> & \begin{matrix}5\\ 6 \end{matrix} & 2 & <1, & \begin{matrix}4\\ 3 \end{matrix}> \\
		\end{array} \right)
\label{eq:cycperm34}
\end{equation}	

%

	\begin{enumerate}
		\item Case $(3,4)_{1}$: $f(4)=5$ and either $f(1)=3$ or $f(1)=4$.
		\item Case $(3,4)_{1,1}$: $f(4)=5$, $f(1)=4$ and $f(6)=1$ or $3$.
		\item Case $(3,4)_{1,1,1}$: $f(4)=5$, $f(1)=4$, $f(2)=6$, $f(6)=3\Rightarrow f(3)\neq 6$ or we get period $2$-suborbit $\left\{ 3,6\right\}$, thus $f(3)=7$ giving valid second minimal $7$ cycle indexed by $5$ in Table \ref{tab:allsecmin7}.
		\item Case $(3,4)_{1,1,2}$: $f(4)=5$, $f(1)=4$, $f(6)=1\Rightarrow f(7)=3\Rightarrow f(3)\neq 7$, or we get period $2$-suborbit $\left\{ 3,7\right\}$, thus $f(3)=6$ giving valid second minimal $7$ cycle indexed by $7$ in Table \ref{tab:allsecmin7}.			
		\item Case $(3,4)_{1,2}$: $f(4)=5$, $f(1)=3$ and either $f(6)=1$ or $f(6)=4$
		\item Case $(3,4)_{1,2,1}$: $f(4)=5$, $f(1)=3$, and $f(6)=1$ then $f(3)\neq 6$ or we have period $3$-suborbit $\left\{ 1,3,6\right\}$ so $f(3)=7$. Then the digraph of the cyclic permutations has the subgraph $J_{1}\rightarrow J_{3}\rightarrow J_{5}\rightarrow J_{1}$ and by Straffin's lemma this implies the existence of a period $3$-suborbit, a contradiction.
		\item Case $(3,4)_{1,2,2}$: $f(4)=5$, $f(1)=3$, and $f(6)=4$ then $f(2)\neq 6$ or we have a period $3$-suborbit $\left\{ 1,3,7\right\}$ so $f(2)=7\Rightarrow$ valid second minimal $7$ cycle indexed by $8$ in Table \ref{tab:allsecmin7}.
		\item Case $(3,4)_{2}$: $f(4)=6$ and $f(2)=5\Rightarrow$ period $2$-suborbit $\left\{ 2,5\right\}$ so $f(2)=7$ and $f(6)=1$ or $f(6)=3$ since $f(6)=4\Rightarrow$ period $2$-suborbit $\left\{ 4,6\right\}$.
		\item Case $(3,4)_{2,1}$: $f(4)=6$, $f(2)=7$, $f(6)=1\Rightarrow f(1)=3$ or $f(1)=4$.
		\item Case $(3,4)_{2,1,1}$: $f(4)=6$, $f(2)=7$, $f(6)=1$, $f(1)=3\Rightarrow$ digraph contains the subgraph $J_{1}\rightarrow J_{3}\rightarrow J_{5}\rightarrow J_{1}$, a contradiction.
		\item Case $(3,4)_{2,1,2}$: $f(4)=6$, $f(2)=7$, $f(6)=1$, $f(1)=4\Rightarrow$ a period $3$-suborbit $\left\{ 1,4,6\right\}$, a contradiction.
		\item Case $(3,4)_{2,2}$: $f(4)=6$, $f(2)=7$, $f(6)=3\Rightarrow f(7)=1\Rightarrow$ valid second minimal $7$ cycle indexed by $9$ in Table \ref{tab:allsecmin7}.
	\end{enumerate}


	\item [Setting (3,5)] From \ref{eq:l5rulesa}, \ref{eq:l5rulesb} $J_{r_{1}}=[4,5]\rightarrow [5,2]$ and $J_{r_{2}}=[2,3]\rightarrow <5, 6\wedge 7>$ so $f(4)=5$ and either $f(2)=5$ or $f(3)=5$, a contradiction.
	
	\item [Setting (4,4)] From \ref{eq:l5rulesa}, \ref{eq:l5rulesc}, \ref{eq:l5rulesd} we have $J_{r_{1}}=[3,4]\rightarrow [4\wedge 5\wedge 6, 2]$, $J_{r_{3}}=[6,7]\rightarrow <1, 2\wedge 3>$, and $J_{r_{4}}=[1,2]\rightarrow [3,7]$. So $f(4)=2$ and $f(1)=3$ but either $f(6)$ or $f(7)$ is $2$ or $3$, a contradiction.
	
	\item [Setting (4,5)] From \ref{eq:l5rulesa}, \ref{eq:l5rulesc}, \ref{eq:l5rulesd} we have $J_{r_{1}}=[3,4]\rightarrow [4\wedge 5, 2]$, $J_{r_{3}}=[5,6]\rightarrow <1, 2\wedge 3>$,  and $J_{r_{4}}=[1,2]\rightarrow [3,6\wedge 7]$. So $f(4)=2$ and $f(1)=3$ but either $f(5)$ or $f(6)$ is $2$ or $3$, a contradiction.
	
	\item [Setting (5,5)] We have the cyclic permutation
\begin{equation}
	\left(\begin{array}{ccccccc}
			1 & 2 & 3_{\ast} & 4 & 5 & 6 & 7 \\ 
			3 & 7 & 4 & 2 & 1 & 5 & 6 
		\end{array} \right)
\label{eq:cycperm55}
\end{equation}	

\noindent The digraph of this cycle admits several subgraphs of length $3$; one of which is $J_{1}\rightarrow J_{6}\rightarrow J_{5}\rightarrow J_{1}$ and by Straffin's lemma this implies a $3$-suborbit, a contradiction.
	
	\end{enumerate}				
\end{enumerate}

\noindent Proceeding in the same fashion for Case $3.2$ generates the inverses of the valid cycles already found. Counting all distinct valid second minimal $7$ orbits we see there are exactly $9$, unique up to an inverse. The topological structure and digraph associated with each of these cyclic permutations are listed in Appendix \ref{app:sevenDigraph}.
	
\section{Universality in Chaos}
\label{sec:UniversalityandChaos}

In this section we present some fascinating results pertaining to universal behavior in the route to chaos for a family of unimodal maps. Specifically, we study continuous endomorphisms, dependent on a parameter,  from an interval to itself: $f_{\lambda}:[0,1]\to[0,1]$ satisfying $f(0)=f(1) = 0$ with a single maximum at some point, $x_{max}$, interior to the interval $[0,1]$ under the iterative relation $x_{n+1} = f_{\lambda}(x_{n})$. We are interested in the asymptotic behavior of $x_n$ for $n\rightarrow \infty$ and how this behavior depends on the parameter $\lambda$. A prototypical example is the logistic map

\begin{equation}
	x_{n+1} = 4\lambda x_{n} \left ( 1-x_{n} \right )
	\label{eq:logmap}
\end{equation}

In 1978, Fiegenbaum \cite{Feigenbaum1978, Feigenbaum1979, Feigenbaum1983} discovered a universal transition mechanism to Chaos through successful period doubling bifurcations. As $\lambda$ increases, the behaviour of $x_n$ for large $n$ changes from periodic to chaotic via bifurcations from the $2^n$ periodic cycle to the $2^{n+1}$ periodic cycle. Two universal constants $\delta = 4.6692016...$ and $\alpha=-2.502907875...$ qualitatively characterize the universal transition route. Let  $\lambda^1_n$ be the value of the parameter when $2^n$-orbit is superstable, i.e. critical point $x_{max}$ is one of the elements of the orbit, and let $d^1_n$ be directed distance from $x_{max}$ to the closest element of the orbit: 
\[ d^1_n=x_{max}-f_{\lambda^1_n}^{2^{n-1}}(x_{max}). \]
Then $\lambda^1_n \uparrow \lambda_\infty$, and for a class of unimodal maps with a quadratic maximum of $1$ has 

\begin{equation}
	 \ \lim_{n\to\infty} \frac{\lambda^1_{n-1} - \lambda^1_{n-2}}{\lambda^1_{n}-\lambda^1_{n-1}}=\delta, \ \ \lim_{n\to\infty}\frac{d^1_n}{d^1_{n+1}}=\alpha.
	\label{eq:deltafeig}
\end{equation}
Having discovered the universality of $\delta$ and $\alpha$ numerically, Feigenbaum proposed the mechanism of it based on the renormalization group approach to critical phenomena in statistical mechanics. He revealed that both of these constants are related to a universal function that governs the period doubling route to chaos and expresses this function as the fixed point of some functional operator. The rigorous proof of Feigenbaum's suggested theory was completed for a class of unimodal maps with quadratic maximum in \cite{Collet1980a, Campanino1981, Lanford1980}. The following is the brief summary of the rigorous universality theory (\cite{Collet1980}).

Map $\psi:[-1,1] \to [-1,1]$ is called $\mathscr{C}^1$-unimodal, if $\psi \in C[-1,1], \psi(0)=1$; $\psi$ is strictly increasing on $[-1,0]$ and strictly decreasing on $[0,1]$; $\psi'(x)\neq 0$ if $x\neq 0$. Let $\varmathbb{P}$ be the space of symmetric $\mathscr{C}^1$-unimodal maps. Choose $\psi \in \varmathbb{P}$ and define
\[ a=a(\psi)=-\psi(1), \ b=b(\psi)=\psi(a). \]
Assume that 
\begin{equation}\label{condition for rescaling}
0<\psi(b)=\psi^2(a)<a<b<1.
\end{equation}
This condition guarantees that the second iteration $\psi^2$ maps $[-a,a]$ to itself. Therefore, {\it the doubling transformation}
\begin{equation}\label{doubling}
\mathscr{F}\psi(x)= -\frac{1}{a}\psi^2(-ax).
\end{equation}
maps $[-1,1]$ to itself. The following properties of $\mathscr{F}$ are key features of the universality theory:
\begin{itemize}
\item $\mathscr{F}$ has a fixed point $g$ with $a=-\alpha^{-1}$. Namely, $g$ solves the functional equation
\begin{equation}\label{g}
g(x)=\alpha g^2\Big (\frac{x}{\alpha}\Big ), \ \  g(0)=1.
\end{equation}
\item The Frechet derivative of $\mathscr{F}$ at the fixed point $g$ has a simple eigenvalue equal to $\delta$; the remainder of the spectrum is contained in the open unit disk. Therefore, $\mathscr{F}$ has a one-dimensional unstable manifold $W_u$ and a co-dimension one stable manifold $W_s$ at $g$.
\item $W_u$ intersects transversally the co-dimension one surface $\Sigma_1$ of maps with superstable 2-orbits:
\[ \Sigma_1=\{\psi: \psi(1)=0 \} \]
\item Consider a set $\Sigma_k$ of maps with superstable $2^k$-orbits (inverse images of $\Sigma_1$), i.e.
\[ \Sigma_k=\mathscr{F}^{-(k-1)}\Sigma_1=\{\psi: \psi=\mathscr{F}^{k-1}\phi, \phi \in \Sigma_1\}, \ k=2,3,... \]
Then the distance between $\Sigma_k$ and $W_s$ decreases like $\delta^{-k}$ for large $k$.
\item Consider an arbitrary one-parameter family $\mu \to \psi_\mu$ of maps and treat it as a curve in $\varmathbb{P}$. Assume that this curve crosses the stable manifold $W_s$ at $\mu_{\infty}$ with a non-zero transverse velocity. This implies that for all large $k$, there will be a unique $\mu_k$ near $\mu_{\infty}$, such that $\psi_{\mu_k} \in \Sigma_k$ is a map with superstable $2^k$-orbit. Then 
\[ \lim_{j\to \infty} \mathscr{F}^{j-k}\psi_{\mu_j}=g_k, k=1,2,3,... \  \  \lim_{j\to \infty} \mathscr{F}^{j}\psi_{\mu_{\infty}}=g  \]
where $g_k$ is an intersection of $\Sigma_k$ with $W_u$; $g$ is a fixed point of $\mathscr{F}$ which solves \eqref{g}. All the functions $g_k$ and $g$ are universal functions.
\end{itemize}
The rigorous theory was only developed for a special class of $\mathscr{C}^1$-unimodal maps of the form
\[ \psi(x)=f(|x|^{1+\epsilon}) \]
where the function $f$ is analytic in a complex neighborhood of $[0,1]$, $\epsilon >0$. The typical example would be
\[ \psi(x)=1-\mu |x|^{1+\epsilon}. \]
The perturbative analysis of \cite{Collet1980a} requires $\epsilon$ to be sufficiently small. The case $\epsilon =1$ was completed in \cite{Campanino1981}.

In \cite{Metropolis1973} periodic orbits are characterised through patterns, which is the sequence of R's and L's, the $k$th letter expressing the fact that the $k$th element of the cycle is on the right or left side of the critical point of the map. In particular, paper \cite{Metropolis1973} presents a table of relative position of periodic orbits of period $p\leq 11$ for the logistic map. Much of the work in this direction was inspired by the paper \cite{Milnor1988}, where the calculus for describing the qualitative behaviour of successive iterates of piecewise monotone maps of the interval was invented. We refer to \cite{Collet1980} which presents an extensive description of this approach.  

In this paper, in addition to the logistic map we will present numerical results for the sine map,

\begin{equation}
	f_{\lambda}(x) = \lambda\sin(\pi x)
	\label{eq:sinmap}
\end{equation}

\noindent the cubic map,

\begin{equation}
	f_{\lambda}(x) = \frac{3\sqrt{3}}{2}\lambda x (1-x^{2})
	\label{eq:cubmap}
\end{equation}

\noindent and the quartic map

\begin{equation}
	f_{\lambda}(x) = \lambda - \lambda(2x-1)^{4}
	\label{eq:quarticmap}f
\end{equation}

Note that $x_{max}=0.5$ in \eqref{eq:logmap}, \eqref{eq:sinmap}, \eqref{eq:quarticmap} and $x_{max}=1/ \sqrt{3}$ in \eqref{eq:cubmap}. Moreover, only logistic and sine maps are symmetric around $x_{max}$. All three maps demonstrate Feigenbaum transition route to chaos through successful period doubling from $2^n$ to $2^{n+1}$-orbits. Feigenbaum constants $\delta$ and $\alpha$ are the same for logistic, sine and cubic maps. For the quartic map we have $\delta=7.31...$, and $\alpha=-1.69...$. 

It is well-known that for $\lambda> \lambda_\infty$, one can observe all possible periodic orbits within the chaotic regime. Figure 15 demonstrates the {\it bifurcation diagram} - asymptotic behaviour of the sequence $x_n$ as $n\to +\infty$ (periodic orbits or {\it chaotic attractors}) as a function of the parameter of the map. One can clearly see periodic windows in the chaotic regime, the period 3 window being the largest. Let $\lambda=\lambda^3_0$ be the value of the parameter when superstable 3-orbit appears first time when $\lambda > \lambda_\infty$ increases. In fact, periodic orbits of all possible periods appear when $\lambda \in [\lambda_\infty, \lambda^3_0]$. Our goal in this section is to continue the results reported in a recent paper \cite{abdulla2013} and to reveal and analyze a fascinating pattern of distribution of all the periodic windows in this range of the parameter.

\begin{figure}[htb]
		\centering
		\subcaptionbox{Logistic map: $x_{n+1} = 4\lambda x_{n}(1-x_{n})$}{\begin{tikzpicture}
			\node[anchor=south west,inner sep=0] (image) at (0,0) {\includegraphics[width=0.8\textwidth, keepaspectratio=true, trim=310 85 250 55, clip,]{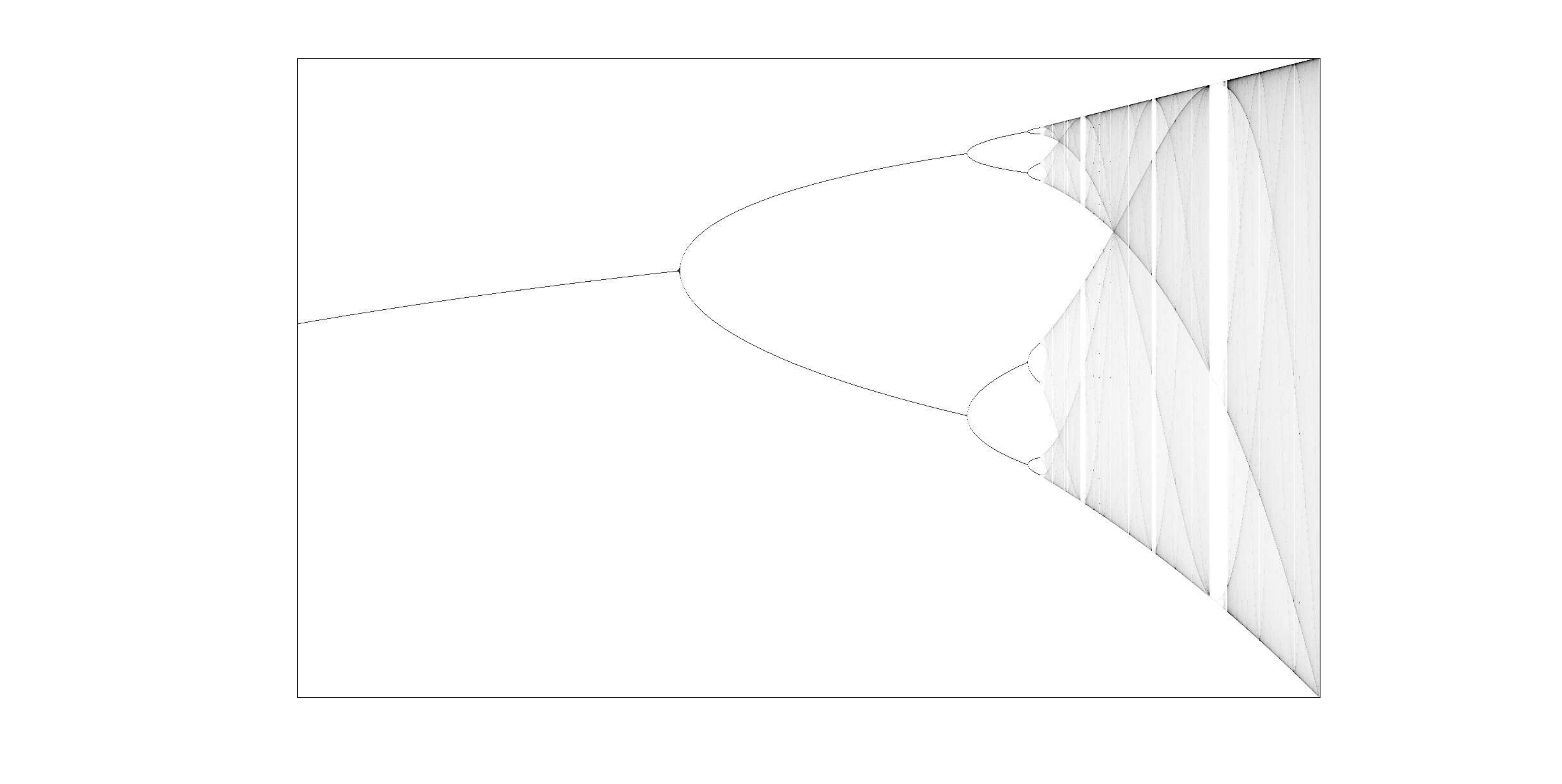}};
				\begin{scope}[x={(image.south east)},y={(image.north west)}]
				\draw (0,0)--(0,1) node[above left] {$x$};
				\draw (0,0)--(1,0) node[below] {$\lambda$};
					\foreach \y in {0, 0.1,0.2,0.3,0.4,0.5,0.6,0.7,0.8,0.9}
						{
							\draw (0.01,\y)--(-0.01,\y) node[left] {$\y$};
						};
				\end{scope}			
			\end{tikzpicture}} 		
			
		\subcaptionbox{Cubic Map}{\includegraphics[width=\figwidth\textwidth, keepaspectratio=true, trim=310 85 250 55, clip,]{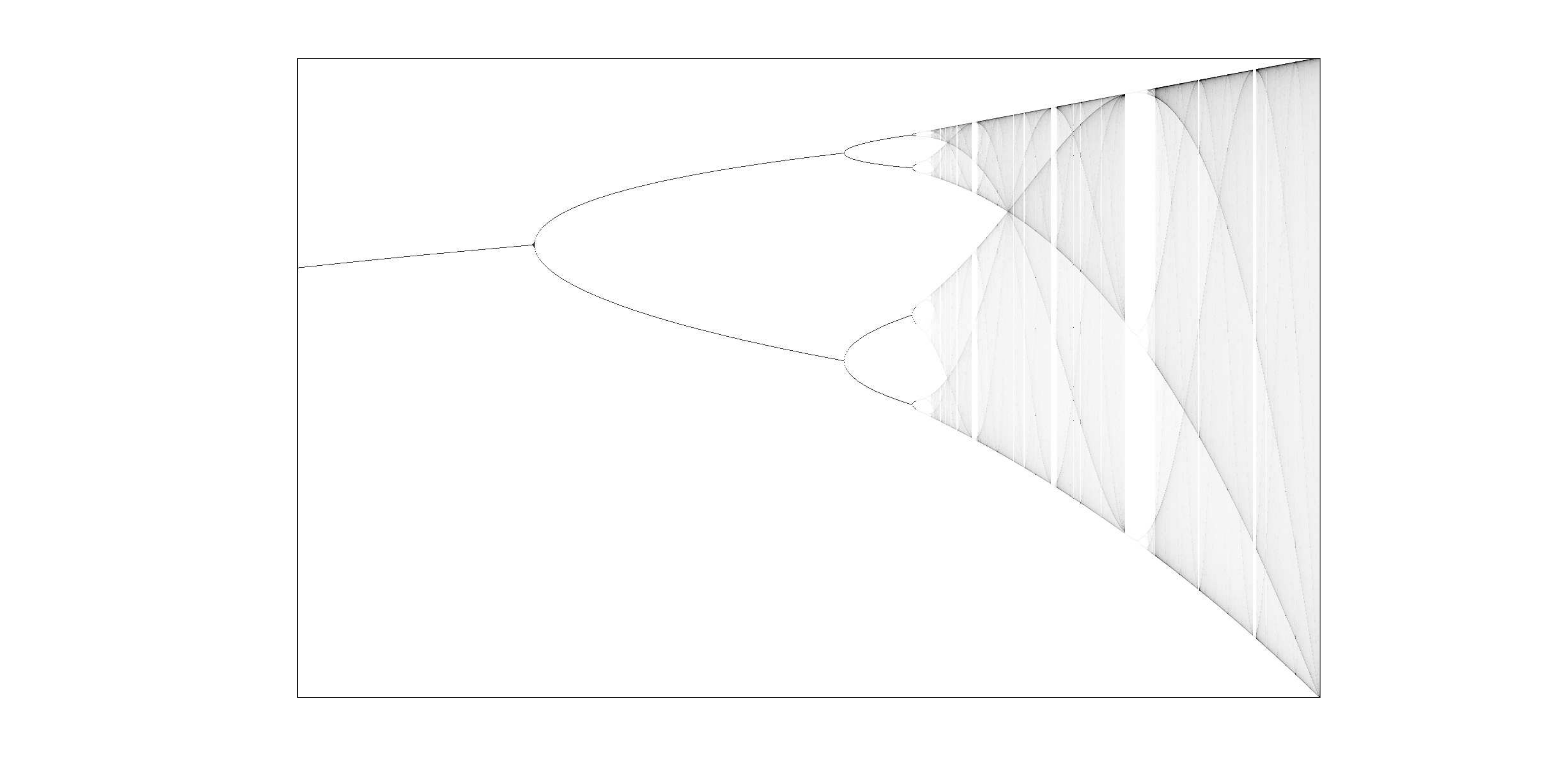}} 				
		\subcaptionbox{Sine Map}{\includegraphics[width=\figwidth\textwidth, keepaspectratio=true, trim=310 85 250 55, clip,]{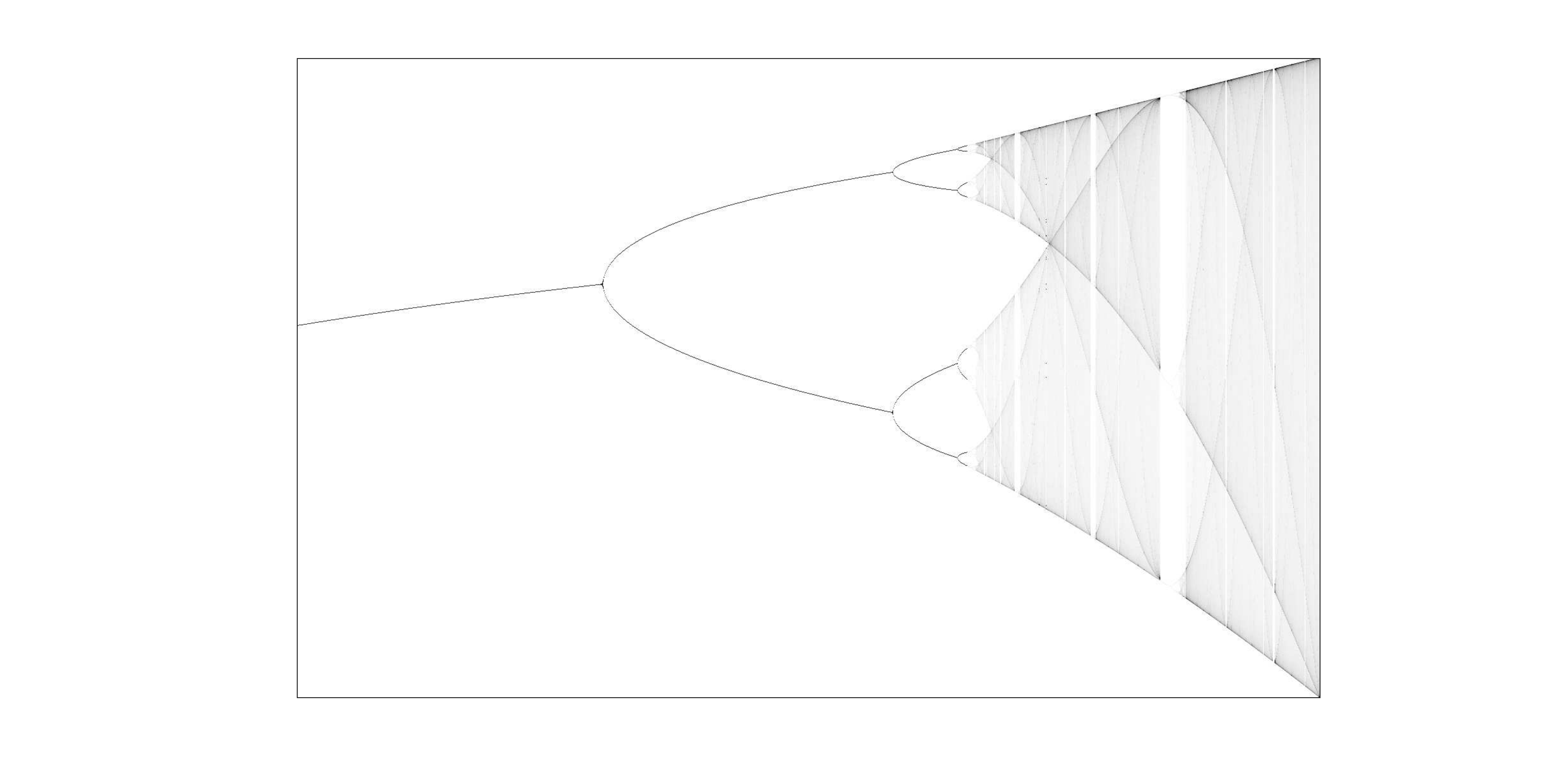}}
		\subcaptionbox{Quartic Map}{\includegraphics[width=\figwidth\textwidth, keepaspectratio=true, trim=310 85 250 55, clip,]{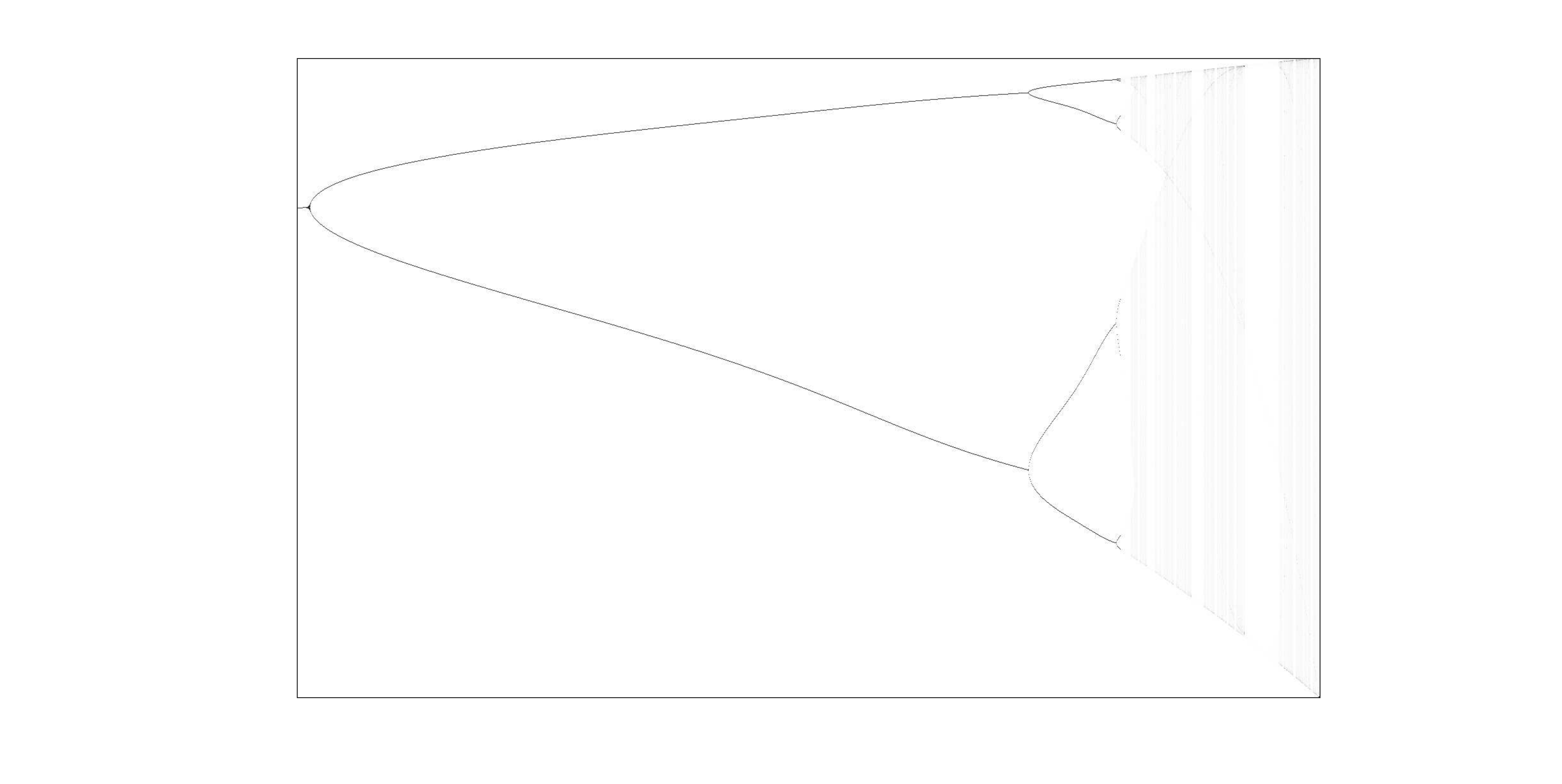}}
		\caption{Bifurcation Diagrams}
		\label{fig:bifdiags}
	\end{figure}
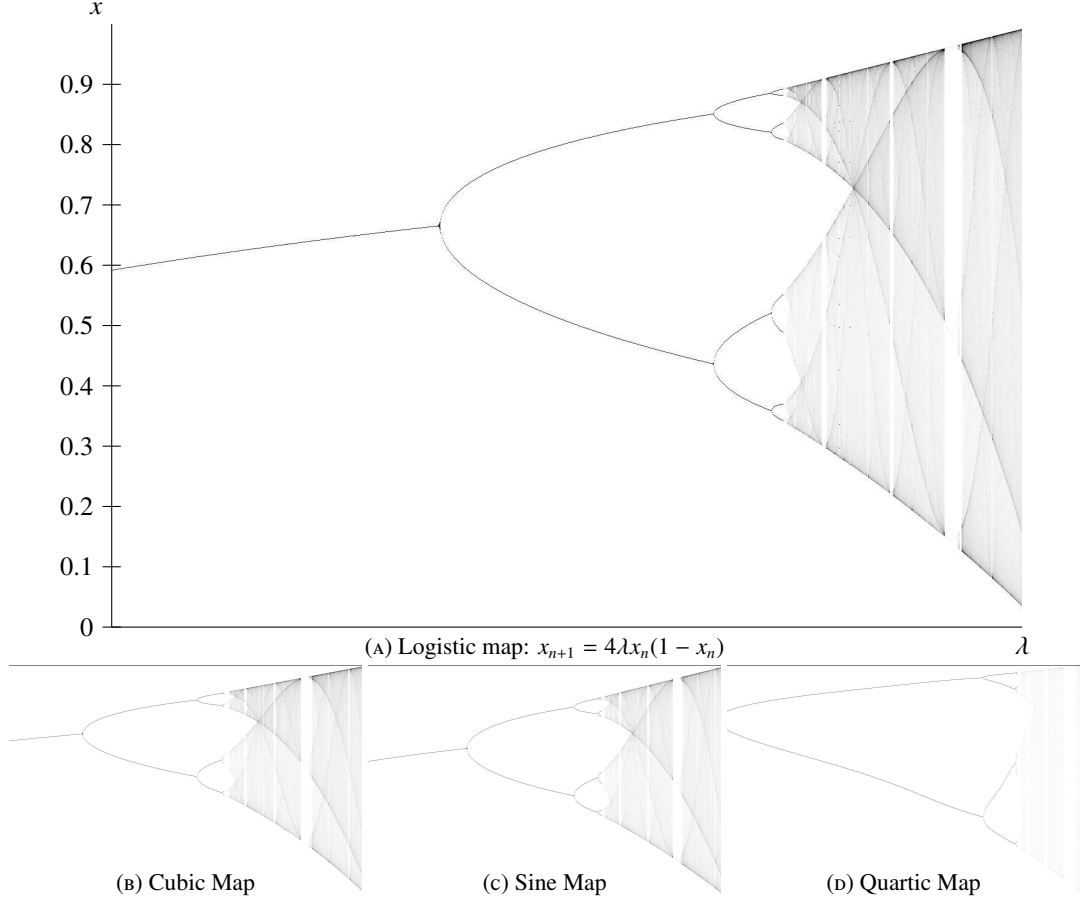		
	
\subsection{Ordering in Terms of the Number of Appearences of Orbits}
\label{subsec:patofpat}
For any odd number $q>1$, and nonnegative integer $s$, let $\Lambda^q_s$ denotes the set of values $\lambda \in [\lambda_\infty, \lambda^3_0]$ such that $f_\lambda$ has a superstable $2^sq$-orbit. In fact, the cardinality $|\Lambda^q_s|$ of $\Lambda^q_s$ is non-zero and finite for all $q$ and $s$. In particular,
\[ \Lambda^3_0=\{\lambda^3_0\}, \  |\Lambda^3_0|=|\Lambda^5_0|=1, \ |\Lambda^7_0|=2, \ |\Lambda^9_0|=4, \ |\Lambda^{11}_0|=9, ... \]
It is well known that the number of appearances of orbits exponentially increases by increasing the period. Let
\[ \Lambda^q_s=\Big \{ \lambda^q_{s,1}< \lambda^q_{s,2}<\cdots < \lambda^q_{s,|\Lambda^q_s|} \Big \}, \]
where $\lambda^q_{s,i}$ denotes the value of the parameter which marks the $i$th appearance of the superstable $2^sq$-orbit when the parameter $\lambda$ increases in the range $[\lambda_\infty, \lambda^3_0]$. Furthermore this orbit will be called $(2^sq)_i$-orbit. Note that $\lambda^3_{0,1}=\lambda^3_0$.
	
Assume that we are looking only first appearance of all the orbits in the indicated parameter range. Numerical results of \cite{abdulla2013} demonstrate that the first appearances of all the orbits are distributed according to Sharkovskii ordering \eqref{sharkovskiordering} when parameter $\lambda$ decreases from $\lambda^3_0$ to $\lambda_\infty$. This is reflected in the first row of the Table \ref{tbl:patofpat}. Moreover, the first appearance of all the orbits is always a minimal orbit. For example, the first appearance of all the odd orbits is always a Stefan orbit and its digraph is as in Figure 1 of Theorem \ref{thm:stefan}. First appearance of all the $2(2k+1)$-orbits always has Type I digraph as in Figure 2 of Theorem \ref{thm:abdulla13}. The reason of relevance of exactly Type I minimal $2(2k+1)$-orbit is hidden in the fact the topological structure of the unimodal map with single maximum is equivalent to the topological structure of the piecewise monotonic map associated with the Type I digraph of Figure 2. In fact, if we iterate the unimodal map with single minimum then inverse Type I digraph will be relevant. 

Assume now that we are looking to first and second appearances of all the $2^sq$-orbits with odd $q\geq 7$, and the first appearance of $2^sq$-orbits with $q=3,5$, while the parameter increases from $\lambda_\infty$ to $\lambda^3_0$. Numerical results of \cite{abdulla2013} and this paper demonstrate the distribution of periodic windows as in \eqref{2nd appear2^n}-\eqref{2nd appear2^0}. Note that we use a notation $n\leftarrow m$ meaning that parameter decreases from the value giving superstable $m$-orbit down to the value giving superstable $n$-orbit.
	
	\begin{align}
			\cdots\leftarrow (2^{n}11)_2 \leftarrow (2^{n}7)_1 \leftarrow (2^{n}9)_2 \leftarrow (2^{n}5)_1 \leftarrow (2^{n}7)_2 \leftarrow (2^{n}3)_1 \leftarrow\cdots\label{2nd appear2^n} \\
			\vdots\nonumber \\
			\cdots\leftarrow 36_1\leftarrow 44_2 \leftarrow 28_1 \leftarrow 36_2 \leftarrow 20_1 \leftarrow 28_2 \leftarrow 12_1 \leftarrow\cdots\label{2nd appear2^2} \\
		\cdots\leftarrow 18_1\leftarrow 22_2 \leftarrow 14_1 \leftarrow 18_2 \leftarrow 10_1 \leftarrow 14_2 \leftarrow 6_1 \leftarrow\cdots\label{2nd appear2^1} \ \\
\cdots\leftarrow 13_2\leftarrow 9_1\leftarrow 11_2 \leftarrow 7_1 \leftarrow 9_2 \leftarrow 5_1 \leftarrow 7_2 \leftarrow 3_1 \leftarrow\cdots\label{2nd appear2^0} \
	\end{align}	
	
	\noindent and we have the pattern $\lambda^{2k-1}_{n,1} < \lambda^{2k+1}_{n,2} < \lambda^{2k-3}_{n,1}$ for $k=3,4,\ldots; n = 0,1,\ldots$ and in particular notice that while decreasing the parameter $\lambda$, $(2^sq)$-orbits are changed with respect to $q$ according to pattern $+4-2$; while the index of appearance is changed according to the simple pattern $1,2,1,2,\dots$. This pattern is expressed in the second row of the Table \ref{tbl:patofpat}. Interestingly, the second appearance of all the odd orbits is second minimal odd orbit. In fact, this numerically observed fact was a motivation to introduce the notion of second minimal orbit as in Definition \ref{thm:secondmin}. In fact, in all four maps the second appearance of the 7-orbit is exactly Type I second minimal orbit with cyclic permutation and digraph demonstrated in Table 1 of Theorem \ref{thm:secondmin7} and Fig. 29 in Appendix. The reason of the relevance of exactly Type 1 second minimal 7-orbit is hidden in the fact that the topological structure of the single maximum unimodal map is equivalent to the topological structure of the piecewise monotonic map associated with Type 1 second minimal 7-orbit of Fig. 29. In fact, according to Theorem \ref{thm:secondmin7} among all possible 9 types of second minimal 7-orbits
(Figures 29-37), Type 1 7-orbit is the only one with a unimodal structure with a single maximum point. In fact, if we iterate the unimodal endomorphism with a single minimum point, then the inverse Type I digraph would be relevant.

Assume now that we are identifying up to third appearances of all the $2^sq$-orbits when the parameter increases from $\lambda_\infty$ to $\lambda^3_0$. Numerical results for all four maps demonstrate the distribution of periodic windows as in \eqref{3rd appear2^n}-\eqref{3rd appear2^0}.
	
	\begin{align}
			\cdots\leftarrow (2^{n}11)_3 \leftarrow (2^{n}9)_2 \leftarrow (2^{n}5)_1 \leftarrow (2^{n}9)_3  \leftarrow (2^{n}7)_2 \leftarrow (2^{n}3)_1 \leftarrow\cdots\label{3rd appear2^n} \\
			\vdots\nonumber \\
		\cdots\leftarrow 28_1\leftarrow 44_3 \leftarrow 36_2 \leftarrow 20_1\leftarrow 36_3 \leftarrow 28_2 \leftarrow 12_1 \leftarrow\cdots\label{3rd appear2^2} \\
		\cdots\leftarrow 14_1\leftarrow 22_3 \leftarrow 18_2 \leftarrow 10_1\leftarrow 18_3 \leftarrow 14_2 \leftarrow 6_1 \leftarrow\cdots\label{3rd appear2^1} \\
\cdots\leftarrow 13_2\leftarrow 9_1\leftarrow 13_3\leftarrow 11_2 \leftarrow 7_1 \leftarrow 11_3 \leftarrow 9_2 \leftarrow 5_1\leftarrow 9_3  \leftarrow 7_2 \leftarrow 3_1 \leftarrow\cdots\label{3rd appear2^0}
	\end{align}	
	
	\noindent Note that we have the pattern $\lambda^{2k-1}_{n,1} < \lambda^{2k+3}_{n,3} < \lambda^{2k+1}_{n,2} < \lambda^{2k-3}_{n,1}$ for $k=3,4,\ldots; n = 0,1,\ldots$ and in particular notice that while decreasing the parameter $\lambda$, $(2^sq)$-orbits are changed with respect to $q$ according to pattern +4+2-4; while index of appearance is changed according to pattern 1,2,3,... This pattern is expressed in the third row of the Table \ref{tbl:patofpat}. 

Continuing this process reveals the structure presented in Table \ref{tbl:patofpat}. As an example assume that we are identifying up to 9th appearances of all the $2^sq$-orbits when the parameter increases from $\lambda_\infty$ to $\lambda^3_0$. Numerical results for all four maps demonstrate the distribution of periodic windows according to the pattern expressed in the ninth row of the Table \ref{tbl:patofpat}. It is satisfactory to explain the pattern only for $q$-orbits, $q$ is odd number, since the pattern is preserved for $2^nq$-orbits. As it is demonstrated in \eqref{eq:ninthrowex},  when the parameter $\lambda$ decreases from $\lambda^3_0$ to $\lambda_\infty$, superstable $q$-orbits appear according to pattern +8-2+2-4+4-2+2+2-8 starting with superstable 3-orbit (written in red in \eqref{eq:ninthrowex}), while index of appearance changes according to pattern 1,8,4,7,2,6,3,5,9,...
		
\begin{equation}
				5_{1}\overset{\mathrm{\textcolor{red}{-8}}}{\leftarrow}13_{9}\overset{\mathrm{\textcolor{red}{+2}}}{\leftarrow}11_{5}\overset{\mathrm{\textcolor{red}{+2}}}{\leftarrow}9_{3}\overset{\mathrm{\textcolor{red}{-2}}}{\leftarrow}11_{6}\overset{\mathrm{\textcolor{red}{+4}}}{\leftarrow}7_{2}\overset{\mathrm{\textcolor{red}{-4}}}{\leftarrow}11_{7}\overset{\mathrm{\textcolor{red}{+2}}}{\leftarrow}9_{4}\overset{\mathrm{\textcolor{red}{-2}}}{\leftarrow}11_{8}\overset{\mathrm{\textcolor{red}{+8}}}{\leftarrow} 3_{1}
\label{eq:ninthrowex}
\end{equation}			

To construct the table in general, first consider only appearances that are powers of $2$. Now, say we wanted to construct the $2^{n}$ row of the table, then the two outermost entries, that is, the first and $2^{n}$-th entries are set to $+2(n+1)$ and $-2n$ respectively. Then, the two entries exactly in the middle of the $1$st and $2^{n}$th, namely the $2^{n-1}$st and $2^{n-1}+1$st entries are set to $-2(n-1)$ and $+2(n-1)$ respectively. Now, find the median entries between the two halves, $1$ to $2^{n-1}$ and $2^{n-1}+1$ to $2^n$ and set them to $-2(n-2)$ and $+2(n-2)$, and continue in this fashion setting each new set of median entries to $-2(n-i)$ and $+2(n-i)$ for $i=3,\cdots,n-1$ as illustrated in Figure \ref{fig:pat2n}.

To generate the $N$-th row that is not a power of $2$ say $2^{n} < N < 2^{n+1}$

	\begin{enumerate}
		\item Find the pattern for $2^n$ row
		\item Let $J = N - 2^{n}$
		\item Replace the last $J$ values, $\left \{ p_1, p_2, \cdots, p_j, \cdots, p_J \right \}$, of the $2^{n}$ pattern according to the following rule:
		\begin{enumerate}
			\item If $p_j > 0$, $p_j \rightarrow \left \{ p_{j}+2, -2 \right \}$
			\item If $p_j < 0$, $p_j \rightarrow \left \{ +2, p_{j}-2 \right \}$
		\end{enumerate}
	\end{enumerate}

The procedure to generate the indices is recursive. Given a pattern corresponding to row $i$ of the table to generate the row $i+1$, first counting from $1$, left to right, identify the position of $i$, say it's in position $m$ and insert the new one between positions $m-2$ and $m-1$, unless position $m-1$ is $1$, in which case insert the new (highest) index at the end of the list or in the $(i+1)^{th}$ position. For example, to go from row $7$ to row $8$, we start with row $7$ and observe that the highest index, $7$, is in position $3$ so we insert the new index, $8$, in between the index $1$ in position $1$ and the index $4$ in position $2$. However, in going from row $8$ to row $9$ observe that $8$ is in position $2$ so position $m-1$ is $1$. So, we insert $9$ at the end of the list in position $9$.

\subsection{Constant Shift in Appearences}
\label{subsec:univApp}

Numerical results demonstrate that for all four maps, parameter range $[\lambda_\infty, \lambda^3_0]$ is divided into infinitely many blocks. For arbitrary fixed appearance index $j=1,2,...$ we have
\begin{align}
\lambda^{2k+1}_{s,j} \downarrow \lambda^\infty_s, {\quad as} \ k \uparrow \infty; \ s=0,1,2,\dots, \label{s-blocks}\\
\lambda^3_0>\lambda^\infty_0>\cdots >\lambda^\infty_s > \lambda^\infty_{s+1} > \cdots >\lambda_\infty; \ \ \lambda^\infty_s \downarrow \lambda_\infty, {\quad as} \ s \uparrow \infty.\label{sharkovskiblocks}
\end{align}
Note that the limit values $\lambda^\infty_s$ in \eqref{s-blocks} are independent of $j$. Moreover, the results presented in a Table \ref{tbl:universalShift} demonstrate exponential convergence in \eqref{s-blocks}:
\begin{equation}\label{universalshift}
\lambda^{2k+1}_{s,j}-\lambda^\infty_s \sim C \delta_s^{-k}, {\quad as} \ k\uparrow \infty,
\end{equation}
where $C$ is some positive constant, and $\delta_s$ is a convergence rate. With the notation $\overset{m}{\delta}$ in Table \ref{tbl:universalShift}, we expressed the fact that $m$ is the highest period of orbit used for the approximation of the convergence rate $\delta$. For example, $\delta_0=2.817...$ is calculated for up to a 31-orbit and it is approximately the same up to the 5th appearance of all the odd orbits. This results demonstrates that for any fixed two appearance indices, the ratio of distances of parameter values for respective appearances of superstable $2^s(2k+1)$-orbits is an asymptotically positive constant for large $k$, i.e. for any fixed positive integers $i$ and $j$ we have
\[ \lim_{k\to\infty}\frac{\lambda^{2k+1}_{s,j}-\lambda^\infty_s}{\lambda^{2k+1}_{s,i}-\lambda^\infty_s}=C>0. \]

	\begin{table}[!htbp]
		\centering
		\caption{Convergence Rates $\left ( \frac{\lambda_{s-1} - \lambda_{s-2}}{\lambda_{s}-\lambda_{s-1}} \right )$ for $2^{s}(2k+1)$ orbits}	
			{\begin{tabular}{p{2.7cm} | l | l | l | l | l}\toprule
			$s$ & Appearance & \multicolumn{4}{c}{$\overset{\mathrm{highest\,\,orbit\,\,used}}{\mathrm{Convergence\,\,Rate}}$} \\ \midrule
			\multicolumn{2}{c}{}   & Logistic    & Sine        & Cubic       & Quartic \\ \midrule
				0 & 1 & $\overset{31}{2.81758}$ & $\overset{31}{2.93749}$ & $\overset{31}{2.96453}$ & $\overset{31}{3.95368}$ \\
				0 & 2 & $\overset{31}{2.81747}$ & $\overset{31}{2.93741}$ & $\overset{31}{2.96448}$ & $\overset{31}{3.95363}$ \\
				0 & 3 & $\overset{31}{2.81734}$ & $\overset{31}{2.93731}$ & $\overset{31}{2.96437}$ & $\overset{31}{3.95362}$ \\
				0 & 4 & $\overset{31}{2.81712}$ & $\overset{31}{2.93713}$ & $\overset{31}{2.96421}$ & $\overset{31}{3.95358}$ \\
				0 & 5 & $\overset{31}{2.81707}$ & $\overset{31}{2.93710}$ & $\overset{31}{2.96402}$ & $\overset{31}{3.95351}$ \\
				1 & 1 & $\overset{30}{2.92338}$ & $\overset{38}{2.94158}$ & $\overset{38}{2.94044}$ & $\overset{38}{4.54383}$ \\
				1 & 2 & $\overset{42}{2.95071}$ & $\overset{42}{2.93561}$ & $\overset{42}{2.93446}$ & $\overset{42}{4.53395}$ \\
				1 & 3 & $\overset{38}{2.91317}$ & $\overset{38}{2.89814}$ & $\overset{38}{2.89683}$ & $\overset{38}{4.52329}$ \\
				1 & 4 & $\overset{34}{2.73591}$ & $\overset{34}{2.72223}$ & $\overset{34}{2.60199}$ & $\overset{34}{4.32894}$ \\
				2 & 1 & $\overset{84}{2.94355}$ & $\overset{68}{2.93108}$ & $\overset{68}{2.93112}$ & $\overset{68}{4.40456}$ \\
				2 & 2 & $\overset{92}{2.94121}$ & $\overset{76}{2.91648}$ & $\overset{76}{2.91649}$ & $\overset{76}{4.40001}$ \\
				2 & 3 & $\overset{100}{2.94257}$ & $\overset{84}{2.92619}$ & $\overset{84}{2.92622}$ & $\overset{84}{4.40308}$ \\
				2 & 4 & $\overset{84}{2.94243}$ & $\overset{84}{2.92603}$ & $\overset{84}{2.92618}$ & $\overset{84}{4.40154}$ \\\bottomrule
			\end{tabular}}
		\label{tbl:universalShift}
	\end{table}	

\begin{landscape}
\begin{figure}[!htbp]
	\centering
	\begin{tikzpicture}
		\draw (0,0.5)--(0,-0.5);
		\draw (14,0.5)--(14,-0.5);
		\node[left] at (0,0) {$2^{n}$};
		\node at (1, 0) {$+2(n+1)$};
		\node at (3, 0) {$-2(n-2),$};
		\node at (4.5, 0) {$+2(n-2)$};
		\node at (6.5, 0) {$-2(n-1),$};
		\node at (8, 0) {$+2(n-1)$};		
		\node at (10, 0) {$-2(n-2),$};
		\node at (11.5, 0) {$+2(n-2)$};	
		\node at (13.5, 0) {$-2n$};			
			
		\draw [decorate,decoration={brace,amplitude=3pt},rotate around={90:(2.5,-0.3)},=90] (2.5,-0.3) -- (2.5,0.5);	
		\draw [decorate,decoration={brace,amplitude=3pt},rotate around={90:(6,-0.3)},=90] (6,-0.3) -- (6,0.5);
		\draw [decorate,decoration={brace,amplitude=3pt},rotate around={90:(9.5,-0.3)},=90] (9.5,-0.3) -- (9.5,0.5);
		\draw [decorate,decoration={brace,amplitude=3pt},rotate around={90:(13,-0.3)},=90] (13,-0.3) -- (13,0.5);		
						
		\node[below] at (2.1, -0.4) {$2^{n-2}-2$};
		\node[below] at (5.6, -0.4) {$2^{n-2}-2$};
		\node[below] at (9.1, -0.4) {$2^{n-2}-2$};
		\node[below] at (12.6, -0.4) {$2^{n-2}-2$};
		
	\end{tikzpicture}
	\caption{Generation of the $2^{n}$-th row of Table \ref{tbl:patofpat}}
	\label{fig:pat2n}
\end{figure}
%
\begin{table}[!htbp]
		\caption{Pattern of Patterns}
		\centering
			{\scalebox{0.85}{\begin{tabular}{c | c | c | c | c | c | c | c | c | c | c | c | c | c | c | c | c | c}\toprule
			Appearance & \multicolumn{16}{c |}{Pattern} & indexes of Appearance \\ \midrule
			$1=2^{0}$ & \multicolumn{16}{| c |}{$+2$} & $1$ \\
			$2=2^{1}$ & \multicolumn{8}{| c |}{$+4$} & \multicolumn{8}{| c |}{$-2$} & $1,2$ \\
			$3$ & \multicolumn{8}{| c |}{$+4$} & \multicolumn{4}{| c |}{$+2$} & \multicolumn{4}{| c |}{$-4$} & $1,2,3$ \\
			$4=2^{2}$ & \multicolumn{4}{| c |}{$+6$} & \multicolumn{4}{| c |}{$-2$} & \multicolumn{4}{| c |}{$+2$} & \multicolumn{4}{| c |}{$-4$} & $1,4,2,3$ \\
			$5$ & \multicolumn{4}{| c |}{$+6$} & \multicolumn{4}{| c |}{$-2$} & \multicolumn{4}{| c |}{$+2$} & \multicolumn{2}{| c |}{$+2$} & \multicolumn{2}{| c |}{$-6$} & $1,4,2,3,5$ \\
			$6$ & \multicolumn{4}{| c |}{$+6$} & \multicolumn{4}{| c |}{$-2$} & \multicolumn{2}{| c |}{$+4$} & \multicolumn{2}{| c |}{$-2$} & \multicolumn{2}{| c |}{$+2$} & \multicolumn{2}{| c |}{$-6$} & $1,4,2,6,3,5$ \\
			$7$ & \multicolumn{4}{| c |}{$+6$} & \multicolumn{2}{| c |}{$+2$} & \multicolumn{2}{| c |}{$-4$} & \multicolumn{2}{| c |}{$+4$} & \multicolumn{2}{| c |}{$-2$} & \multicolumn{2}{| c |}{$+2$} & \multicolumn{2}{| c |}{$-6$} & $1,4,7,2,6,3,5$ \\
			$8=2^{3}$ & \multicolumn{2}{| c |}{$+8$} & \multicolumn{2}{| c |}{$-2$} & \multicolumn{2}{| c |}{$+2$} & \multicolumn{2}{| c |}{$-4$} & \multicolumn{2}{| c |}{$+4$} & \multicolumn{2}{| c |}{$-2$} & \multicolumn{2}{| c |}{$+2$} & \multicolumn{2}{| c |}{$-6$} & $1,8,4,7,2,6,3,5$ \\
			$9$ & \multicolumn{2}{| c |}{$+8$} & \multicolumn{2}{| c |}{$-2$} & \multicolumn{2}{| c |}{$+2$} & \multicolumn{2}{| c |}{$-4$} & \multicolumn{2}{| c |}{$+4$} & \multicolumn{2}{| c |}{$-2$} & \multicolumn{2}{| c |}{$+2$} & $+2$ & $-8$ & $1,8,4,7,2,6,3,5,9$ \\
			$10$ & \multicolumn{2}{| c |}{$+8$} & \multicolumn{2}{| c |}{$-2$} & \multicolumn{2}{| c |}{$+2$} & \multicolumn{2}{| c |}{$-4$} & \multicolumn{2}{| c |}{$+4$} & \multicolumn{2}{| c |}{$-2$} & $+4$ & $-2$ & $+2$ & $-8$ & $1,8,4,7,2,6,3,10,5,9$ \\
			$11$ & \multicolumn{2}{| c |}{$+8$} & \multicolumn{2}{| c |}{$-2$} & \multicolumn{2}{| c |}{$+2$} & \multicolumn{2}{| c |}{$-4$} & \multicolumn{2}{| c |}{$+4$} & $+2$ & $-4$ & $+4$ & $-2$ & $+2$ & $-8$ & $1,8,4,7,2,6,11,3,10,5,9$ \\
			$12$ & \multicolumn{2}{| c |}{$+8$} & \multicolumn{2}{| c |}{$-2$} & \multicolumn{2}{| c |}{$+2$} & \multicolumn{2}{| c |}{$-4$} & $+6$ & $-2$ & $+2$ & $-4$ & $+4$ & $-2$ & $+2$ & $-8$ & $1,8,4,7,2,12,6,11,3,10,5,9$ \\
			$13$ & \multicolumn{2}{| c |}{$+8$} & \multicolumn{2}{| c |}{$-2$} & \multicolumn{2}{| c |}{$+2$} & $+2$ & $-6$ & $+6$ & $-2$ & $+2$ & $-4$ & $+4$ & $-2$ & $+2$ & $-8$ & $1,8,4,7,13,2,12,6,11,3,10,5,9$ \\
			$14$ & \multicolumn{2}{| c |}{$+8$} & \multicolumn{2}{| c |}{$-2$} & $+4$ & $-2$ & $+2$ & $-6$ & $+6$ & $-2$ & $+2$ & $-4$ & $+4$ & $-2$ & $+2$ & $-8$ & $1,8,4,14,7,13,2,12,6,11,3,10,5,9$ \\
			$15$ & \multicolumn{2}{| c |}{$+8$} & $+2$ & $-4$ & $+4$ & $-2$ & $+2$ & $-6$ & $+6$ & $-2$ & $+2$ & $-4$ & $+4$ & $-2$ & $+2$ & $-8$ & $1,8,15,4,14,7,13,2,12,6,11,3,10,5,9$ \\
			$16=2^{4}$ & $+10$ & $-2$ & $+2$ & $-4$ & $+4$ & $-2$ & $+2$ & $-6$ & $+6$ & $-2$ & $+2$ & $-4$ & $+4$ & $-2$ & $+2$ & $-8$ & $1,16,8,15,4,14,7,13,2,12,6,11,3,10,5,9$\\ \bottomrule
		\end{tabular}}}
		\label{tbl:patofpat}
	\end{table}
\end{landscape}				
		
\subsection{Feigenbaum Universlity in General Classes}
\label{subsec:bifurcations}

Numerical results demonstrate that all the odd orbits which appear in the parameter window $(\lambda_0^\infty, \lambda_0^3]$ are going to go through infinitely many period doubling transformations when $\lambda$ decreases towards $\lambda_\infty$. This is demonstrated in the diagram \eqref{3rd appear2^n}-\eqref{3rd appear2^0} if we consider periods up to 3rd appearances. Let us fix any positive integer $J$ as highest appearance index, and deduce from the $J$th row of the Table \ref{tbl:patofpat} the distribution of all the odd orbits up to $J$th appearance in the parameter window $(\lambda_0^\infty, \lambda_0^3]$ (e.g. if $J=9$ then the portion of the odd orbits up to 9th appearnces between $3_1$ and $5_1$ are demonstrated in \eqref{eq:ninthrowex}). All these orbits are going to go through infinitely many bifurcations when $\lambda$ decreases towarsd $\lambda_\infty$, and for any positive integer $s$, the $s$th bifurcation appears in the parameter window $(\lambda_{s+1}^\infty, \lambda_s^\infty)$.  
It is fascinating that all these transition routes to chaos follow Feigenbaum universality. In particular, it is revealed that the Feigenbaum universality is relevant in very general classes of maps beyond the unimodal smooth endomorphisms.

Let integers $k\geq 1$ and $j\in [1, |\Lambda^{2k+1}_0|]$ be fixed. Recall that $\lambda_{0,j}^{2k+1}$ is the value of the parameter $\lambda$ when superstable $(2k+1)$-orbit appears $j$th time while increasing $\lambda$ from $\lambda_\infty$ to $\lambda^3_0$. Numerical results demonstrate that for all four maps we have
\begin{align}
\lambda^{2k+1}_{s,j} \downarrow \lambda_\infty, {\quad as} \ s \uparrow \infty; \ \label{perioddoubling1}\\
 \lim_{s\to\infty}\frac{\lambda^{2k+1}_{s-1,j}-\lambda^{2k+1}_{s-2,j}}{\lambda^{2k+1}_{s,i}-\lambda^{2k+1}_{s-1,j}}=\delta, \ \label{perioddoubling2}\\
\lambda^{2k+1}_{s,j}-\lambda_\infty \sim C \delta^{-s}, {\quad as} \ s\uparrow \infty, \ \label{perioddoubling3}
\end{align}
where $\delta=4.6692...$ in the case of logistic, sine and cubic maps (Tables \ref{tbl:logFeig}, \ref{tbl:sineFeig}, \ref{tbl:cubFeig}); $\delta=7.31...$ in the case of the quartic map (Table \ref{tbl:quaFeig});  $C>0$. Hence, we see that the convergence rate of the sequence of parameter values for superstable $(2^s(2k+1))_j$-orbits to critical value $\lambda_\infty$ as $s\to +\infty$ from above is the same as the convergence rate of the sequence of parameter values for the superstable $2^s$-orbits to the same value $\lambda_\infty$ from below. To clarify if Feigenbaum universality mechanism is indeed relevant we check asymptotical properties of the scaling factor for successful period doublings from $(2^s(2k+1))_j$- to $(2^{s+1}(2k+1))_j$-orbits. Let $d^{2k+1}_{s,j}$ be a directed distance from the maximum point of the map to the closest element of the superstable $(2^{s+1}(2k+1))_j$-orbit, i.e.
\begin{equation}\label{definitionofd}
d^{2k+1}_{s,j}=x_{max}-f^{2^s(2k+1)}_{\lambda^{2k+1}_{s+1,j}}(x_{max}). 
\end{equation}
Numerical results in Tables \ref{tbl:logFeig}-\ref{tbl:quaFeig} demonstrate that for all four models we have
\begin{equation}\label{scaleconvergence}
\lim_{s\to\infty}\frac{d^{2k+1}_{s-1,j}}{d^{2k+1}_{s,j}}=\alpha, 
\end{equation}
where $\alpha=-2.5029...$ in the case of logistic, sine and cubic maps (Tables \ref{tbl:logFeig}, \ref{tbl:sineFeig}, \ref{tbl:cubFeig}); $\alpha=-1.69...$ in the case of the quartic map (Table \ref{tbl:quaFeig}).  Hence, we see that the scaling factor of the successive bifurcations of the superstable $(2^s(2k+1))_j$-orbits when $\lambda$ converges to critival value $\lambda_\infty$ as $s\to +\infty$ from above is the same as the scaling factor of the successive bifurcations of the superstable $2^s$-orbits when $\lambda$ converges to critical value $\lambda_\infty$ from below. This indicates that the doubling transformation \eqref{doubling} with scaling factor $a=\alpha$ is a driving force for the transition to chaos through successful bifurcations of superstable $(2^s(2k+1))_j$-orbits for $s=0,1,2,...$. Therefore, Feigenbaum's universality theory should be valid beyond the class of $\mathscr{C}^1$-unimodal maps - the classes of maps whose structure is defined with $q$th iteration of unimodal maps, where $q=2k+1$ is any fixed odd number. Following Feigenbaum \cite{Feigenbaum1978, Feigenbaum1979} define the functions
\begin{equation}
	g_{m}^{2k+1}(x) = \lim\limits_{s\to\infty}\left ( -\alpha \right )^{s}f^{2^s(2k+1)}_{\lambda_{s+m,j}^{2k+1}}\left ( \frac{x}{\left ( -\alpha \right )^{s}} \right ), \ m=1,2,...
\label{eq:univg}
\end{equation}
\begin{equation}
	g^{2k+1}(x)=\lim\limits_{m\to \infty}g_{m}^k(x) = \lim\limits_{s\to\infty}\left ( -\alpha \right )^{s}f^{2^s(2k+1)}_{\lambda_{\infty}}\left ( \frac{x}{\left ( -\alpha \right )^{s}} \right ), \ k=1,2,...
	\label{eq:univg-g}
\end{equation}
Numerical results demonstarte that for any fixed non-negative integer $k$, family of functions in \eqref{eq:univg}, \eqref{eq:univg-g} are universal functions. The case $k=0$
in \eqref{eq:univg}, \eqref{eq:univg-g} is a particular case of classical Feigenbaum universality theory explaining the transition from $2^s$-orbits, $s=0,1,2,...$ to chaos through successful period doublings (see \eqref{eq:deltafeig} and following description of the rigorous universality theory). In this case $g^{(1)}=g$ is a fixed point of the doubling operator $\mathscr{F}$ as in \eqref{g}; each $g^{(1)}_m=g_m$ is the intersection of $\Sigma_m$ with the one-dimensional unstable manifold passing through $g$. Figure~ \ref{fig:prdDblMech21} demonstrates the convergence to the universal function $g_1$ in Figure~\ref{fig:univ}(a) for the logistic map after calculation of the few terms under the limit sign in \eqref{eq:univg}. 

Figures~\ref{fig:prdDblMech31}, \ref{fig:prdDblMech71}, \ref{fig:prdDblMech72}, \ref{fig:prdDblMech91}, \ref{fig:prdDblMech92} demonstrate the convergence in \eqref{eq:univg} to universal functions $g_1^{2k+1}$ via successive bifurcations of superstable $(2^s3)_1$-, $(2^s7)_1$-, $(2^s7)_2$-, $(2^s9)_1$-,  $(2^s9)_2$-orbits respectively under the transition \eqref{perioddoubling1}-\eqref{perioddoubling3} for the logistic map. Note that the side length of the green square in each of the Figures~\ref{fig:prdDblMech31}-\ref{fig:prdDblMech92} is equal to respective value of $d_{s,j}^{2k+1}$, and Figures~\ref{fig:univ}(b),(c),~\ref{fig:univ2}(a)(b) demonstrate the first four terms of the limit expression in \eqref{eq:univg}.

In fact, the universal functions $g^{2k+1}, k=1,2,...$ in \eqref{eq:univg-g} are fixed points of the doubling operator $\mathscr{F}$, and solve the functional equation in \eqref{g}. That is the reason that the convergence rate of parameter sequences in  \eqref{perioddoubling1}-\eqref{perioddoubling3} is the same universal constant $\delta$. Moreover, It is easy to prove that if function $g$ solves functional equation \eqref{g}, then any iteration of $g$ is also a solution of the same equation. The normalization condition $g(0)=1$ can be arranged by replacing $g$ with $g_\mu=\mu g(x/\mu)$, and by choosing the constant $\mu$ appropriately. Indeed, for arbitrary $\mu\neq 0$, $g_\mu$ is a solution of the functional equation \eqref{eq:univg} if $g$ is so.  
Hence, universal functions $g^{2k+1}$ must be exactly $2k+1$st iterations of the universal function $g$ (which is the justification of our notation), which is the fixed point of the doubling operator in the class of  $\mathscr{C}^1$-unimodal maps. For any fixed $k=1,2,...$, $g^{2k+1}$ represents a fixed point of the doubling operator (and hence solving the functional equation in \eqref{g}) in the more complicated class of maps which is the $(2k+1)$st iteration of the class of  $\mathscr{C}^1$-unimodal maps.

\begin{figure}
	\centering
	\includegraphics[width=0.34\textwidth]{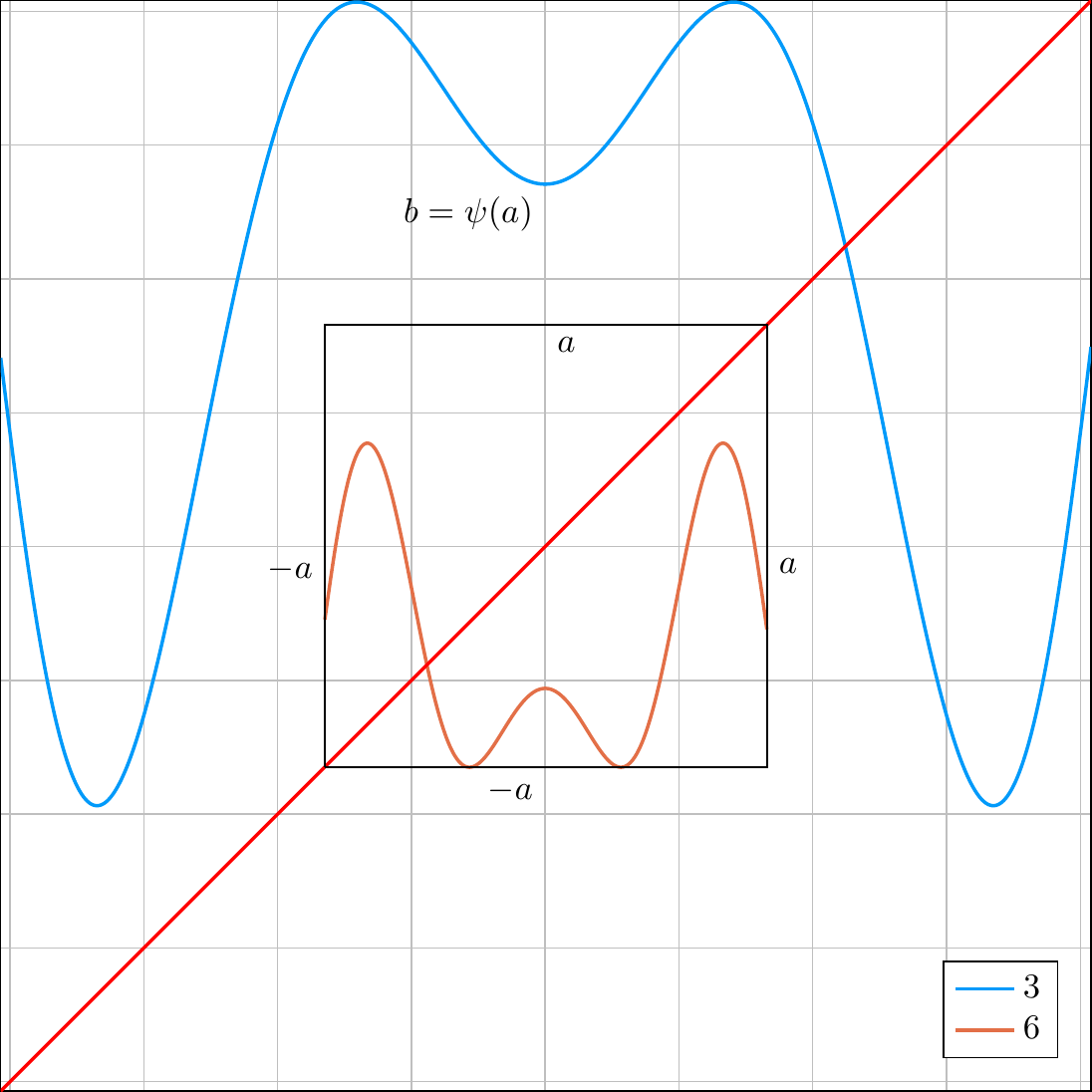}
	\caption{Period doubling mechanism for $3_1$ showing the scaling}
\label{graphofpsi3}
\end{figure}

Hence, the numerical analysis suggests that the known rigorous universality theory (\cite{Collet1980}) must be true in a much larger class of maps than $\mathscr{C}^1$-unimodal maps, and this generalization is a driving force of infinitely many Feigenbaum scenarios of transition to chaos through successive bifurcations of all possible odd orbits as it is outlined in \eqref{perioddoubling1}-\eqref{eq:univg-g}. We end our presentation with the description of the anticipated rigorous universality theory in the particular case of $k=1$, or in the class of maps which is the 3rd iteration of the $\mathscr{C}^1$-unimodal maps. Let
\[ \varmathbb{P}^3=\{\phi: \phi=\psi^3, \psi \in \varmathbb{P} \}. \]
Assume that $\psi \in  \varmathbb{P}$ satisfies \eqref{condition for rescaling}. Since $\psi$ is continuous, there exists $e\in(a,1)$ such that $\psi(e)=0$, and
$\psi^2$ is increasing and maps $[0,e]$ onto $[-a,1]$; $\psi^2$ is decreasing and maps $[e,1]$ onto $[b,1]$. By continuity there exists $d\in (0,a)$ such that $\psi^2(d)=0$. Now consider symmetric function $\phi=\psi^3$. $\phi$ is increasing and maps $[0,d]$ onto $[b,1]$; $\phi$ is decreasing and maps $[d,e]$ onto $[-a,1]$; $\phi$ is increasing and maps $[e,1]$ onto $[-a,\psi(b)]$; This guarantees that the second iteration $\phi^2=\psi^6$ maps $[-a,a]$ to itself. Indeed, first of all from \eqref{condition for rescaling} it follows that $\phi(a)>b$, and hence, $\phi$ maps $[-a,a]$ to $[b,1]$. Also, since $\psi^2$ maps $[-a,a]$ to itself, we have $\phi(b)=\psi^4(a)\leq a$. Accordingly, $\phi$ maps $[b,1]$ into $[-a,a]$, and hence $\phi^2$ maps $[-a,a]]$ to itself. Therefore, the doubling transformation $\mathscr{F}$ maps $[-1,1]$ into itself. Figure~\ref{graphofpsi3} demonstrates the structure of $\phi=\psi^3$ and $\phi^2=\psi^6$ under the condition \eqref{condition for rescaling}. 

The following properties of $\mathscr{F}$ are key features of the universality theory in the class $\varmathbb{P}^3$:
\begin{itemize}
\item $\mathscr{F}$ has a fixed point $g^3$ with $a=-\alpha^{-1}$. Namely, $g^3$ solves the functional equation in \eqref{g} in the class $\varmathbb{P}^3$.
In fact, $g^3$ is precisely 3rd iteration of the fixed point of the doubling operator $\mathscr{F}$ in the class of $\mathscr{C}^1$-unimodal maps, defined in \eqref{g}.
\item The Frechet derivative of $\mathscr{F}$ at the fixed point $g^3$ has a simple eigenvalue equal to $\delta$; the remainder of the spectrum is contained in the open unit disk. Therefore, $\mathscr{F}$ has a one-dimensional unstable manifold $W_u$ and a codimension one stable manifold $W_s$ at $g^3$.
\item $W_u$ intersects transversally the codimension-one surface $\Sigma^3_1$ of maps with superstable 2-orbits:
\[ \Sigma^3_1=\{\phi \in \varmathbb{P}^3: \phi^2(0)=0 \} \]
\item Consider a set $\Sigma^3_m$ of maps with superstable $2^m$-orbits (inverse images of $\Sigma^3_1$), i.e.
\[ \Sigma^3_m=\mathscr{F}^{-(m-1)}\Sigma^3_1=\{\phi: \phi=\mathscr{F}^{m-1}\phi_0, \phi_0 \in \Sigma^3_1\}, \ m=2,3,... \]
Then the distance between $\Sigma^3_m$ and $W_s$ decreases like $\delta^{-m}$ for large $m$.
\item Consider arbitrary one-parameter family $\mu \to \phi_\mu$ of maps and treat it as a curve in $\varmathbb{P}^3$. Assume that this curve crosses stable manifold $W_s$ at $\mu_{\infty}$ with non-zero transverse velocity. This implies that for all large $m$, there will be a unique $\mu_m$ near $\mu_{\infty}$ such that $\phi_{\mu_m} \in \Sigma^3_m$ is a map with superstable $2^m$-orbit. Then 
\[ \lim_{j\to \infty} \mathscr{F}^{j-m}\phi_{\mu_j}=g^3_m, m=1,2,3,... \  \  \lim_{j\to \infty} \mathscr{F}^{j}\psi_{\mu_{\infty}}=g^3  \]
where $g^3_m$ is an intersection of $\Sigma_m$ with $W_u$; $g^3$ is a fixed point of $\mathscr{F}$ which solves functional equation in \eqref{g} in the class $\varmathbb{P}^3$.  All the functions $g^3_m$ and $g^3$ are universal functions.
\end{itemize}
For example, numerical calculation of the universal function $g^3_1$ is demonstrated in Figure~\ref{fig:univ2}(b). 

Similar description of the rigorous universality theory can be outlined in various classes
\[ \varmathbb{P}^{2k+1}=\{\phi: \phi=\psi^{2k+1}, \psi \in \varmathbb{P} \}, \ k=2,3,4,... \]

\captionsetup[subfigure]{skip=0pt}
\begin{figure}[htb]
	\centering		
	\subcaptionbox{$f$ with superstable $2$}{\includegraphics[width=\figwidth\textwidth, keepaspectratio=true]{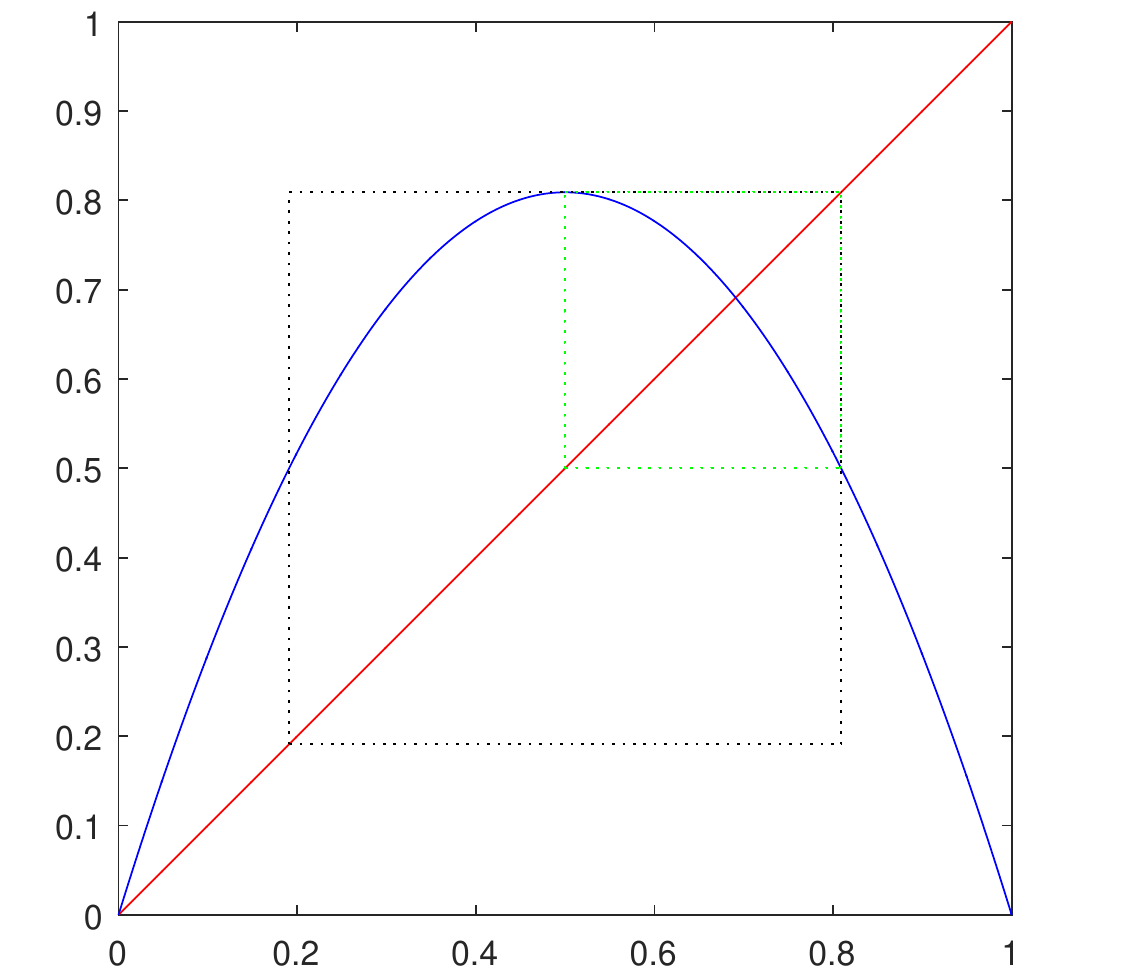}} 					
	\subcaptionbox{$f^{2}$ with superstable $2$}{\includegraphics[width=\figwidth\textwidth, keepaspectratio=true]{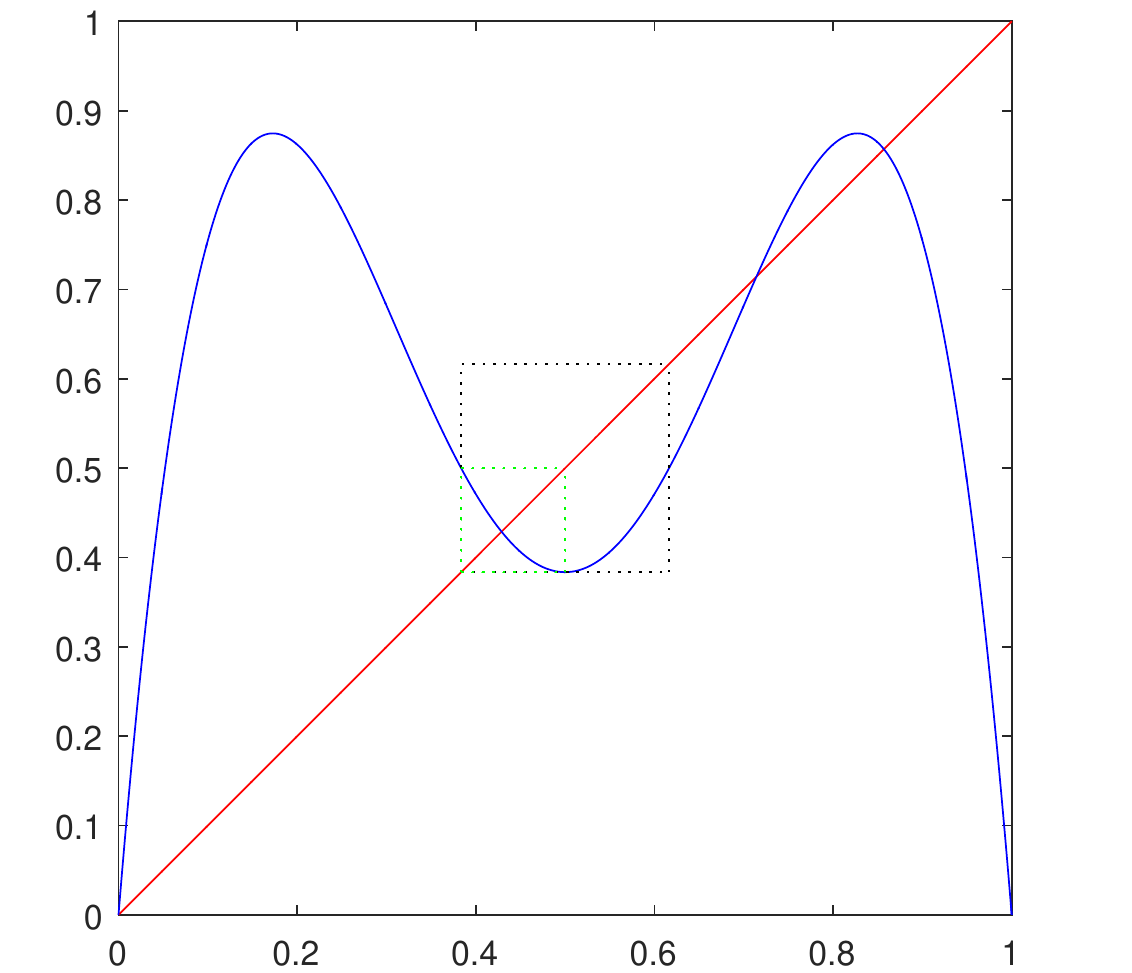}}	
	\subcaptionbox{$f^{4}$ with superstable $2$}{\includegraphics[width=\figwidth\textwidth, keepaspectratio=true]{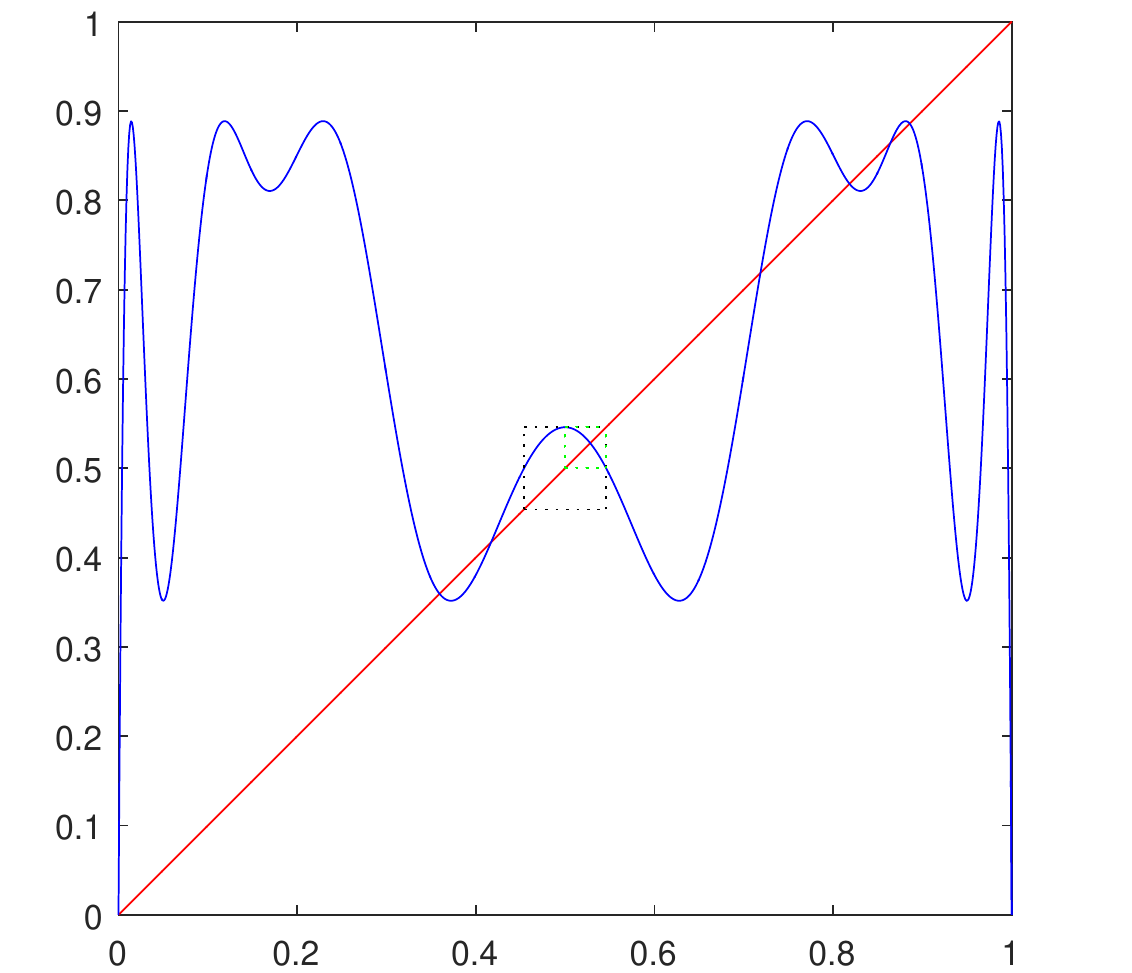}}
	\caption{The Period Doubling Mechanism}
	\label{fig:prdDblMech21}
\end{figure}

\begin{figure}[htb]
	\centering		
	\subcaptionbox{$f^{3}$ with superstable $2$}{\includegraphics[width=\figwidth\textwidth, keepaspectratio=true]{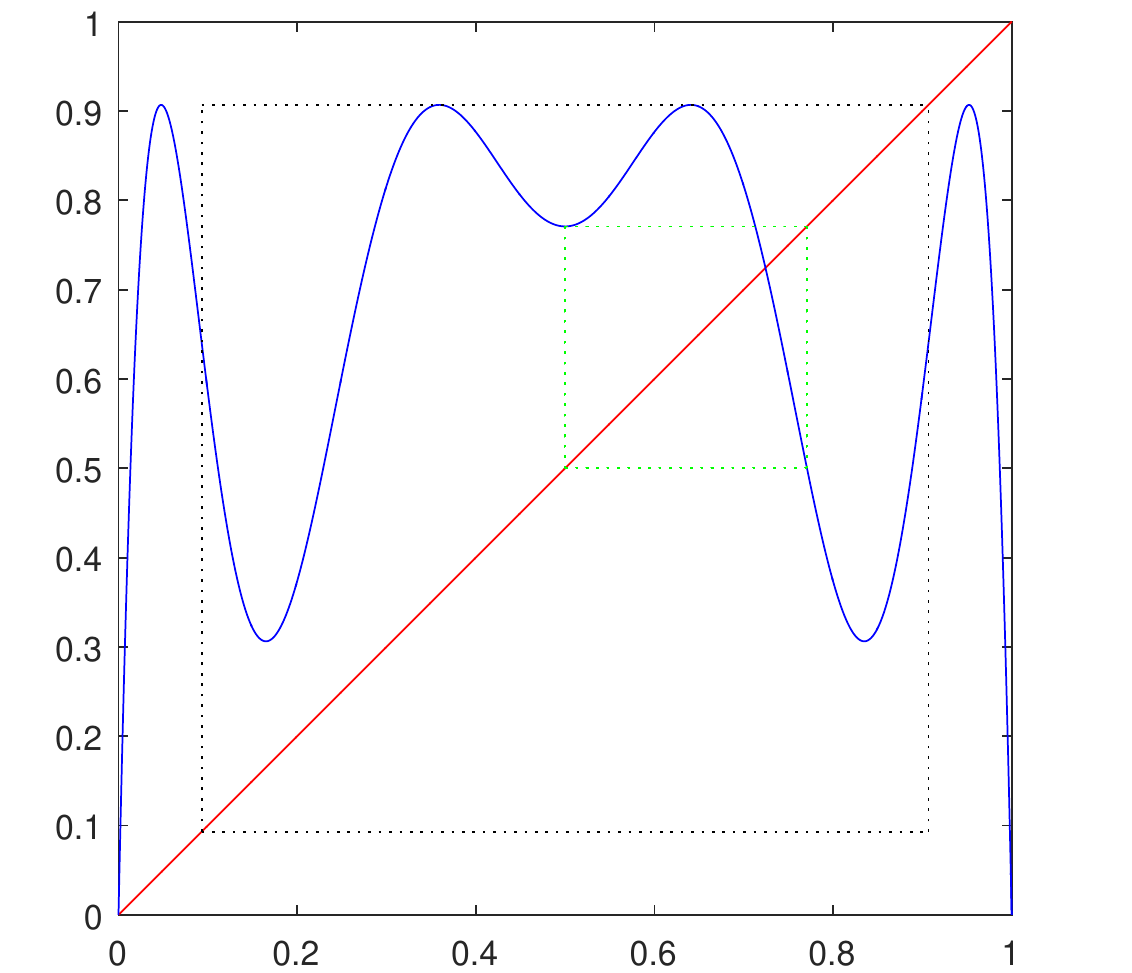}} 					
	\subcaptionbox{$f^{6}$ with superstable $2$}{\includegraphics[width=\figwidth\textwidth, keepaspectratio=true]{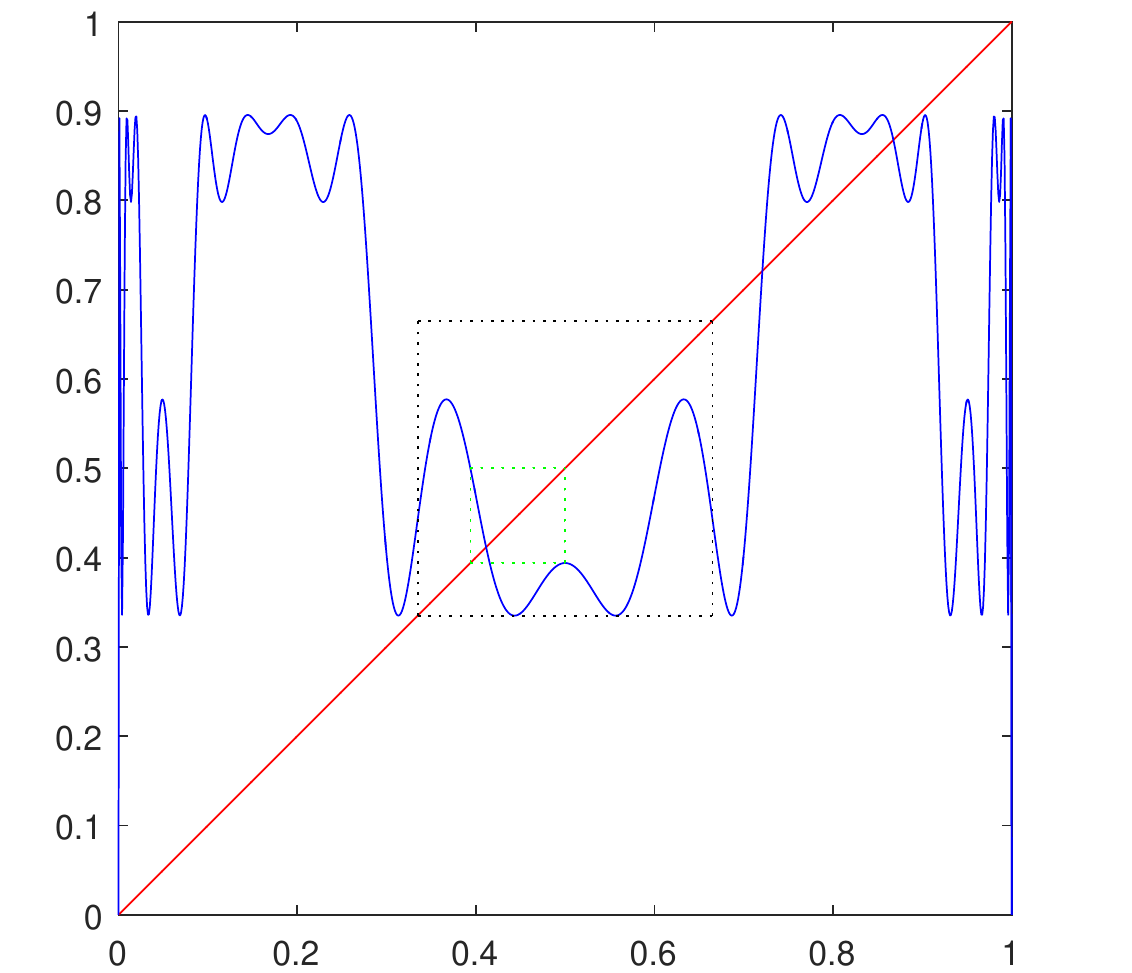}}	
	\subcaptionbox{$f^{12}$ with superstable $2$}{\includegraphics[width=\figwidth\textwidth, keepaspectratio=true]{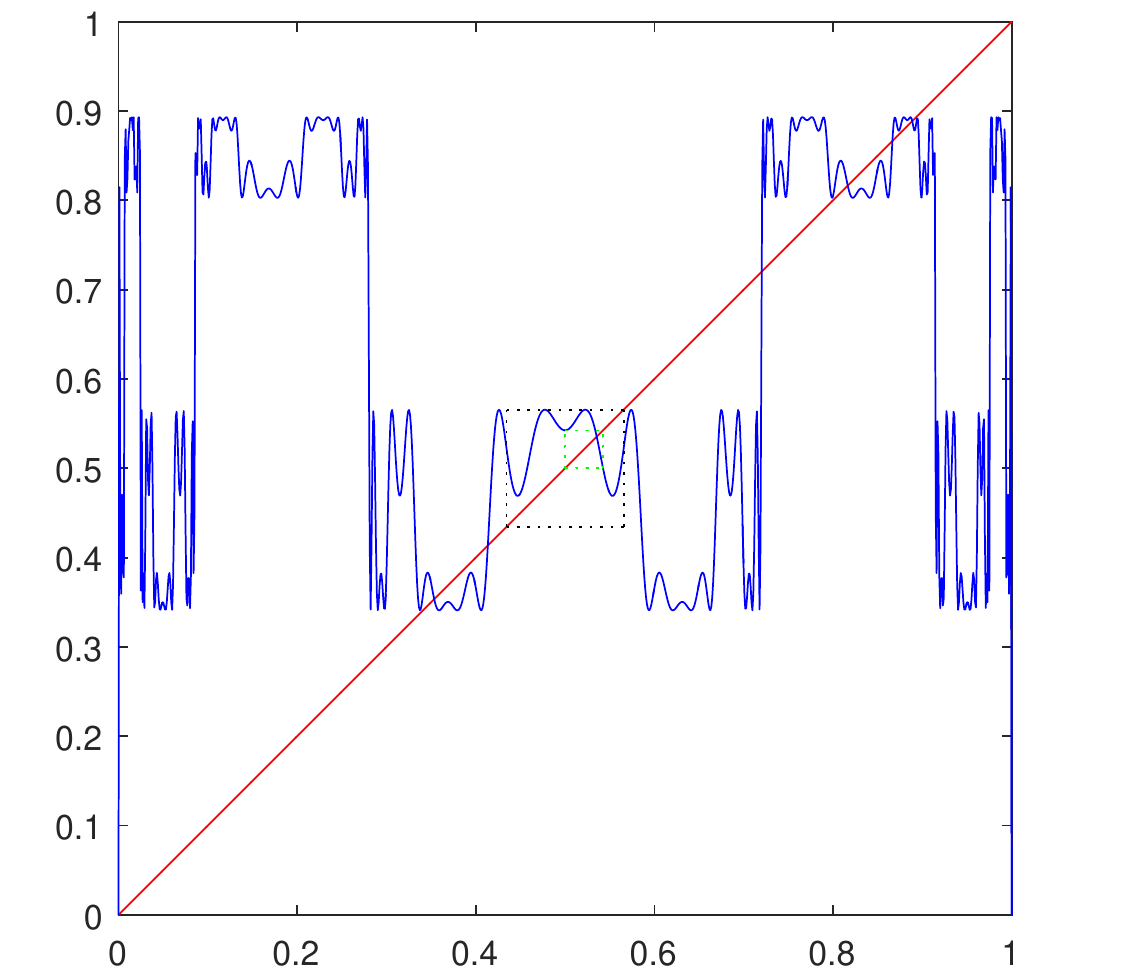}}
	\caption{The Period Doubling Mechanism for $3_1$}
	\label{fig:prdDblMech31}
\end{figure}

\begin{figure}[htb]
	\centering		
	\subcaptionbox{$f^{7}$ with superstable $2$}{\includegraphics[width=\figwidth\textwidth, keepaspectratio=true]{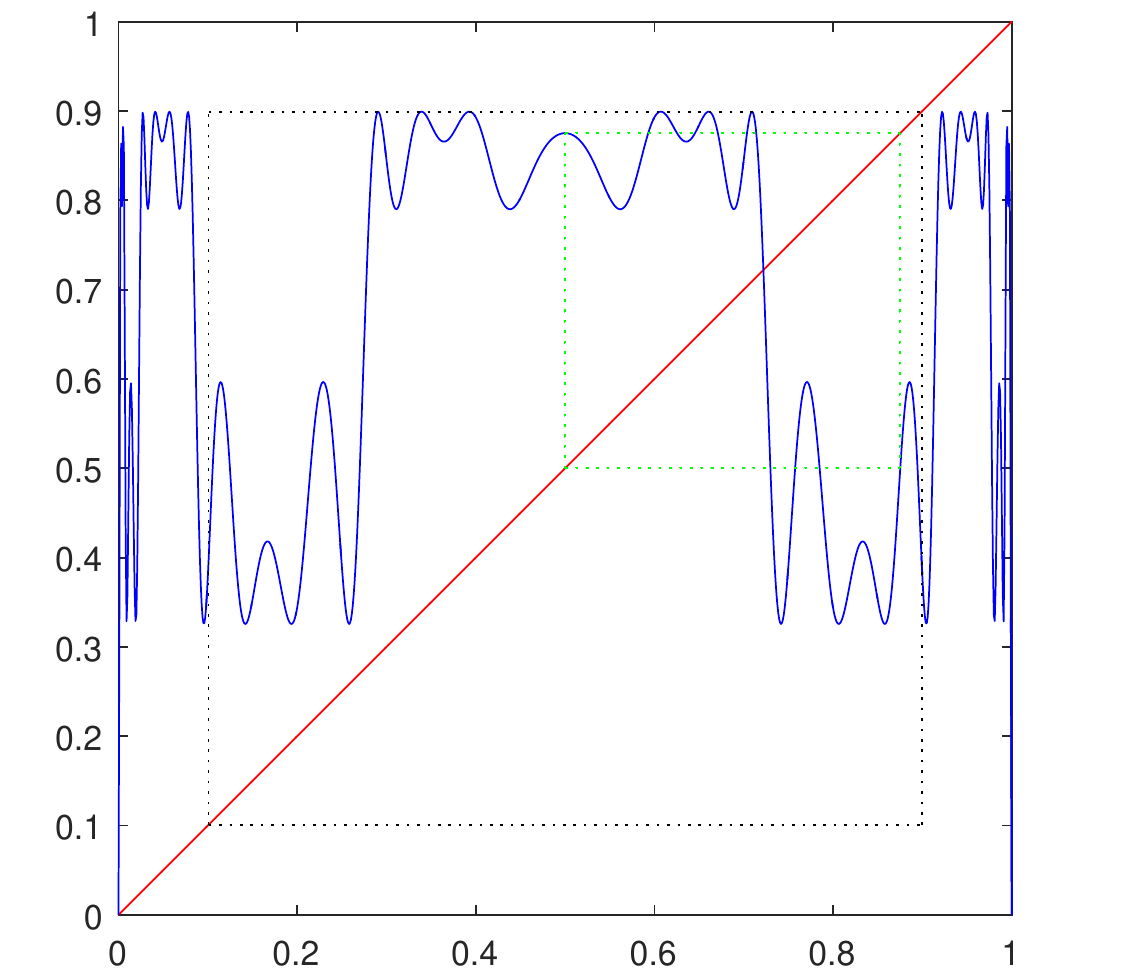}} 					
	\subcaptionbox{$f^{14}$ with superstable $2$}{\includegraphics[width=\figwidth\textwidth, keepaspectratio=true]{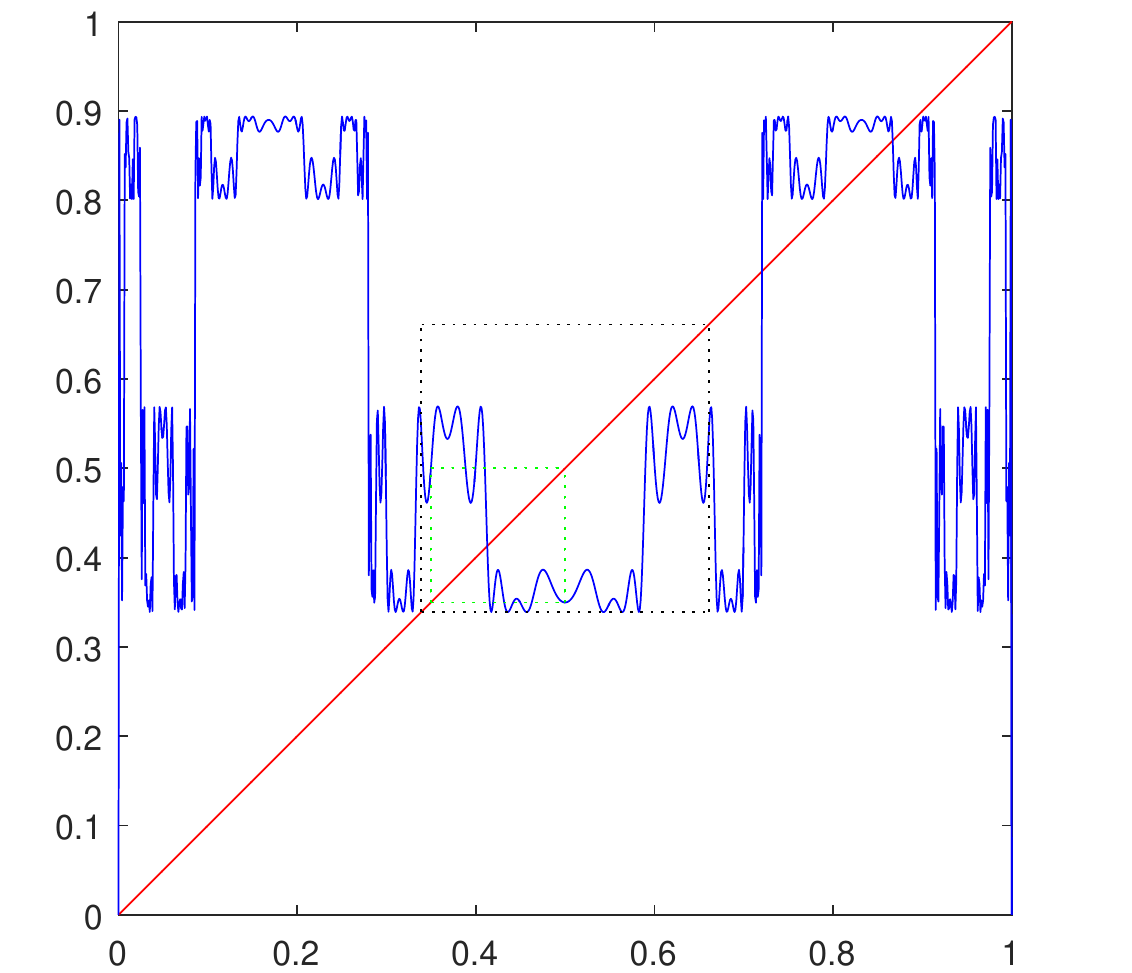}}	
	\subcaptionbox{$f^{28}$ with superstable $2$}{\includegraphics[width=\figwidth\textwidth, keepaspectratio=true]{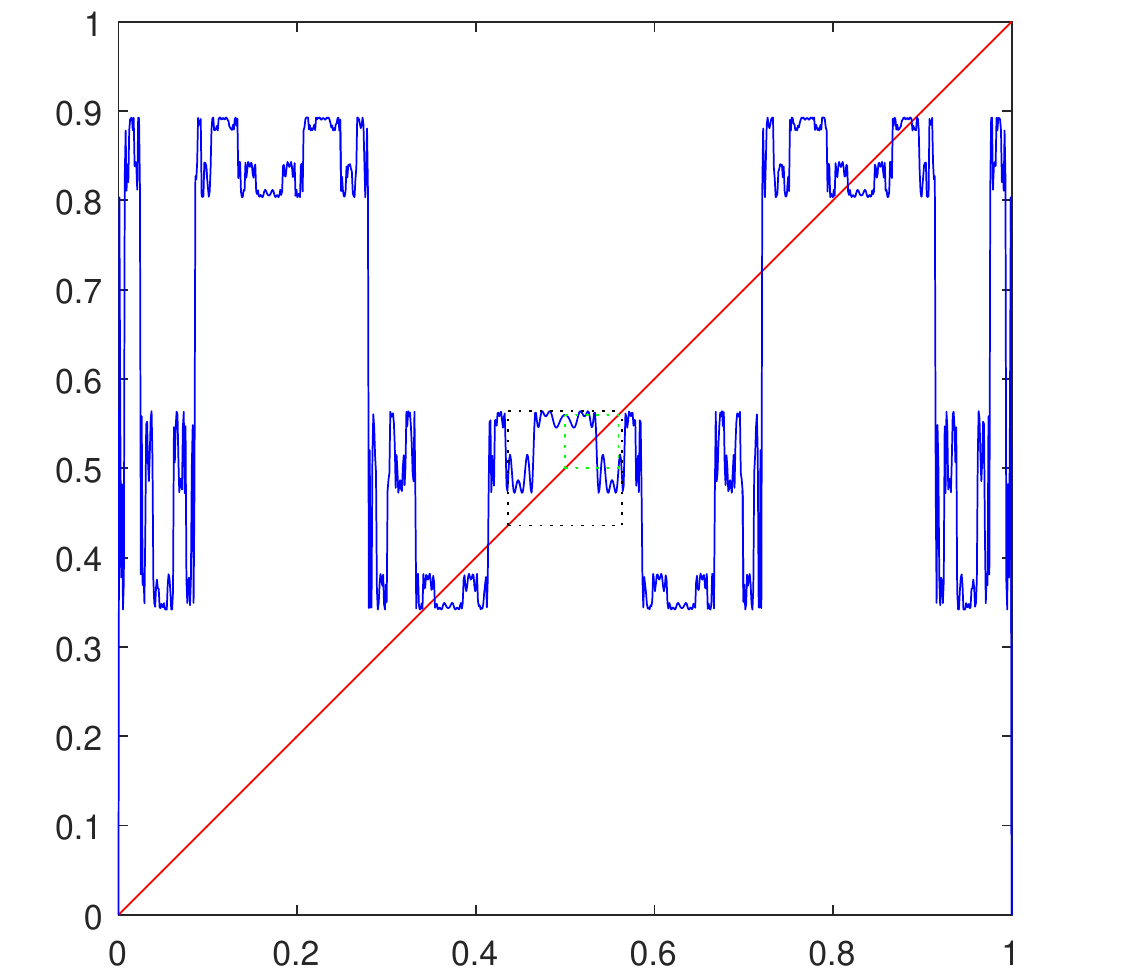}}
	\caption{The Period Doubling Mechanism for $7_1$}
	\label{fig:prdDblMech71}
\end{figure}

\captionsetup[subfigure]{skip=0pt}
\begin{figure}[htb]
	\centering		
	\subcaptionbox{$f^{7}$ with superstable $2$}{\includegraphics[width=\figwidth\textwidth, keepaspectratio=true]{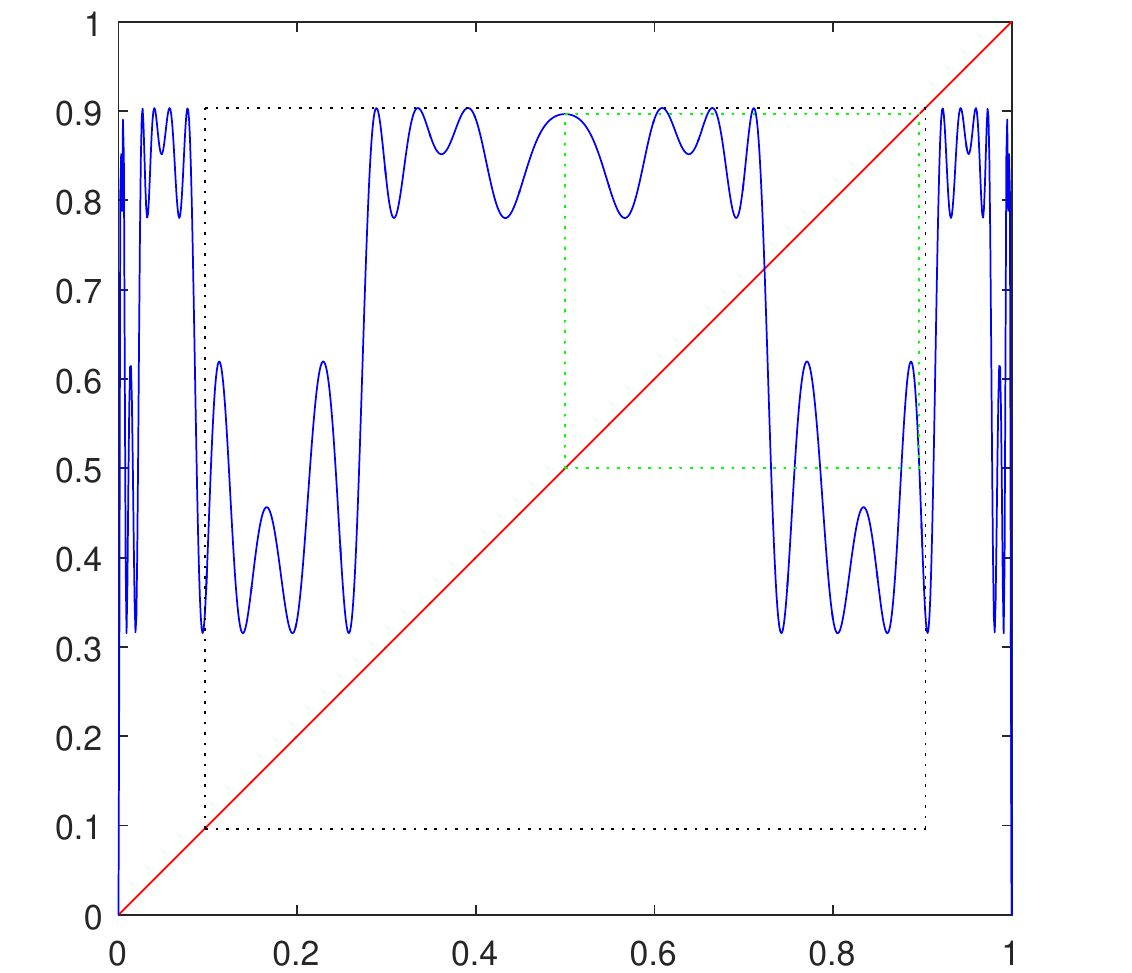}} 					
	\subcaptionbox{$f^{14}$ with superstable $2$}{\includegraphics[width=\figwidth\textwidth, keepaspectratio=true]{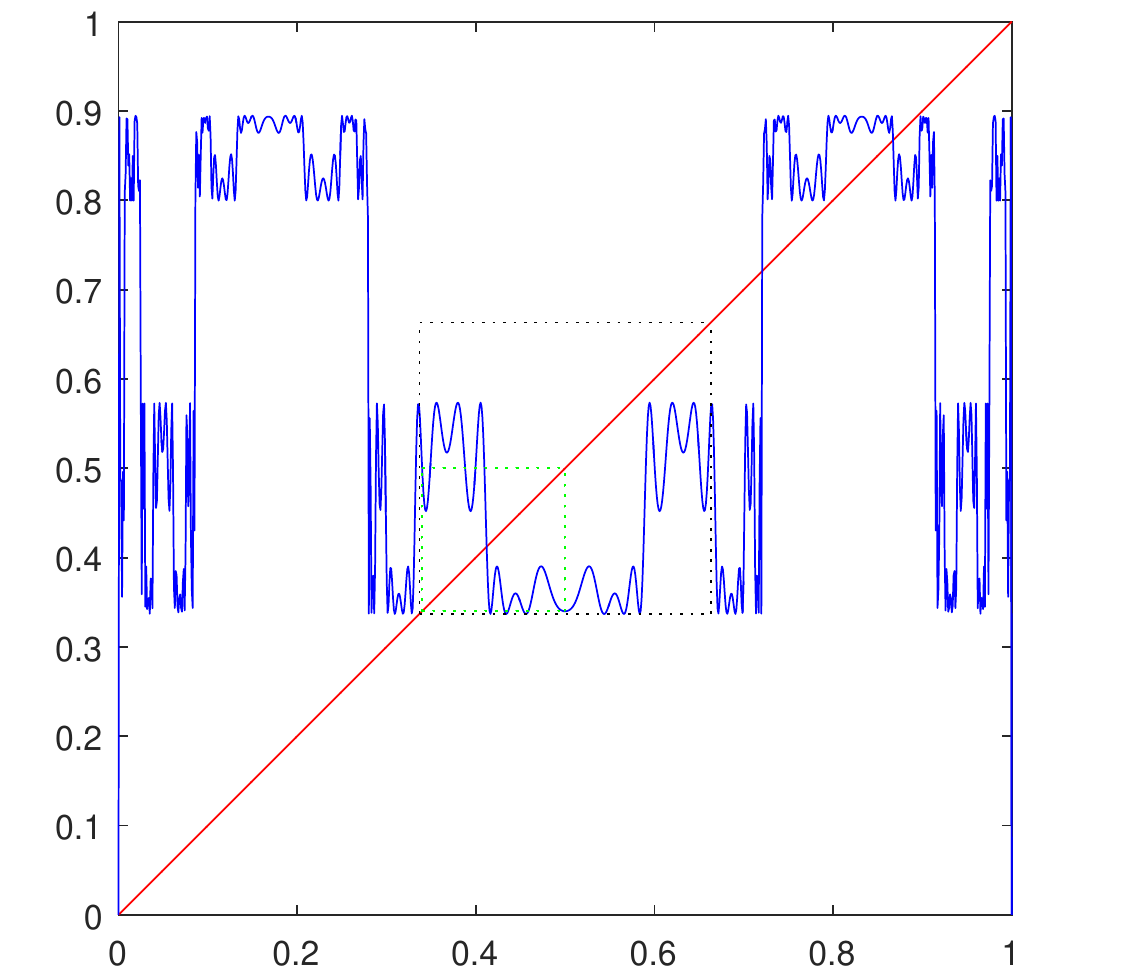}}	
	\subcaptionbox{$f^{28}$ with superstable $2$}{\includegraphics[width=\figwidth\textwidth, keepaspectratio=true]{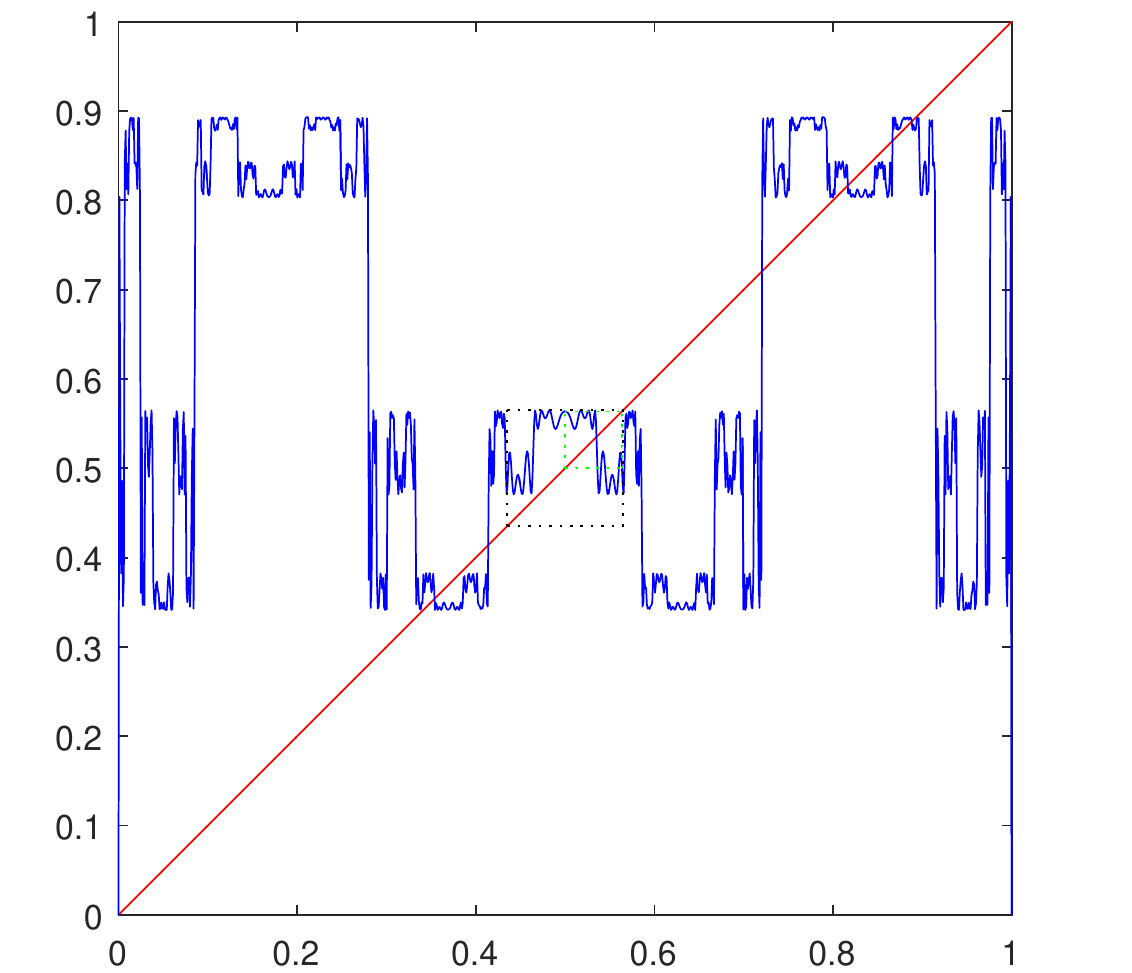}}
	\caption{The Period Doubling Mechanism for $7_2$}
	\label{fig:prdDblMech72}
\end{figure}

\pagebreak

\begin{figure}[htb]
	\centering		
	\subcaptionbox{$f^{9}$ with superstable $2$}{\includegraphics[width=\figwidth\textwidth, keepaspectratio=true]{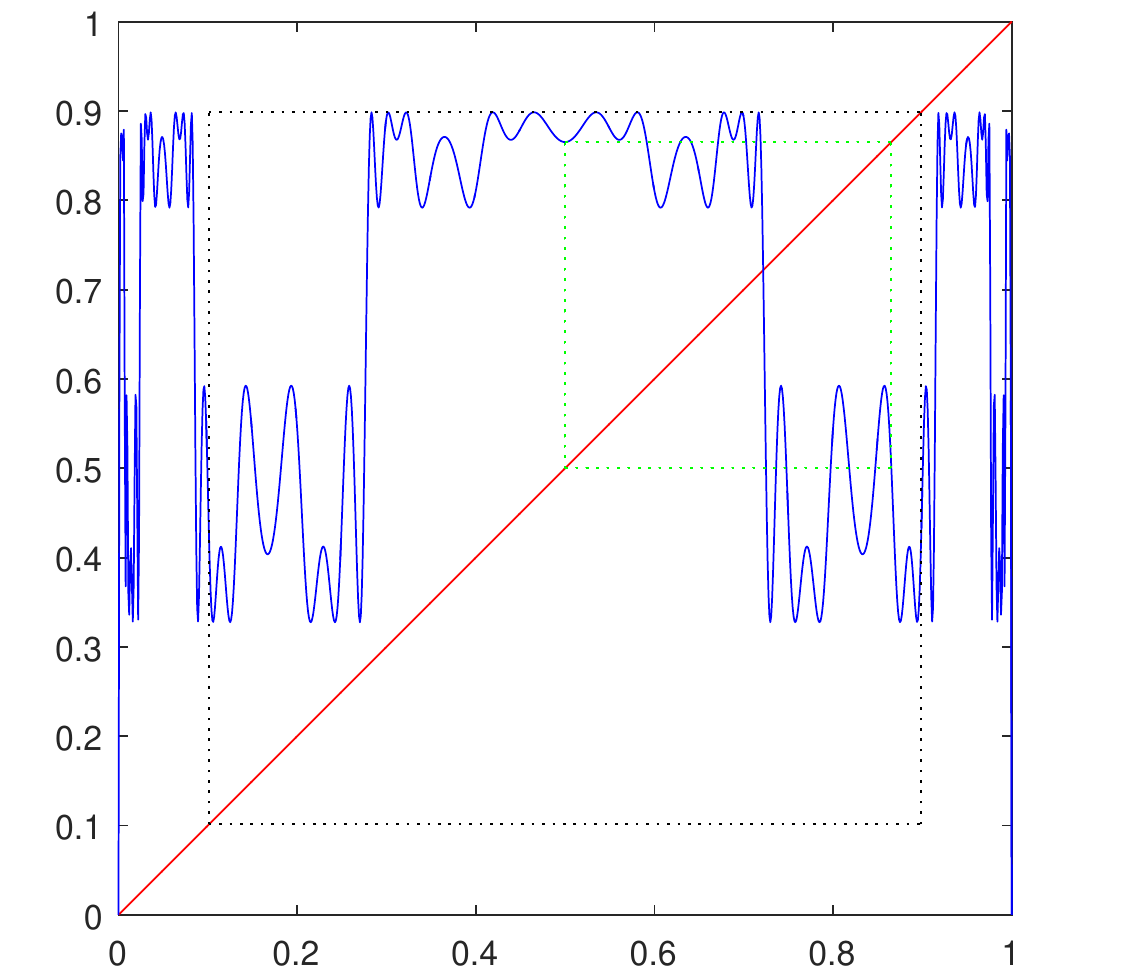}} 					
	\subcaptionbox{$f^{18}$ with superstable $2$}{\includegraphics[width=\figwidth\textwidth, keepaspectratio=true]{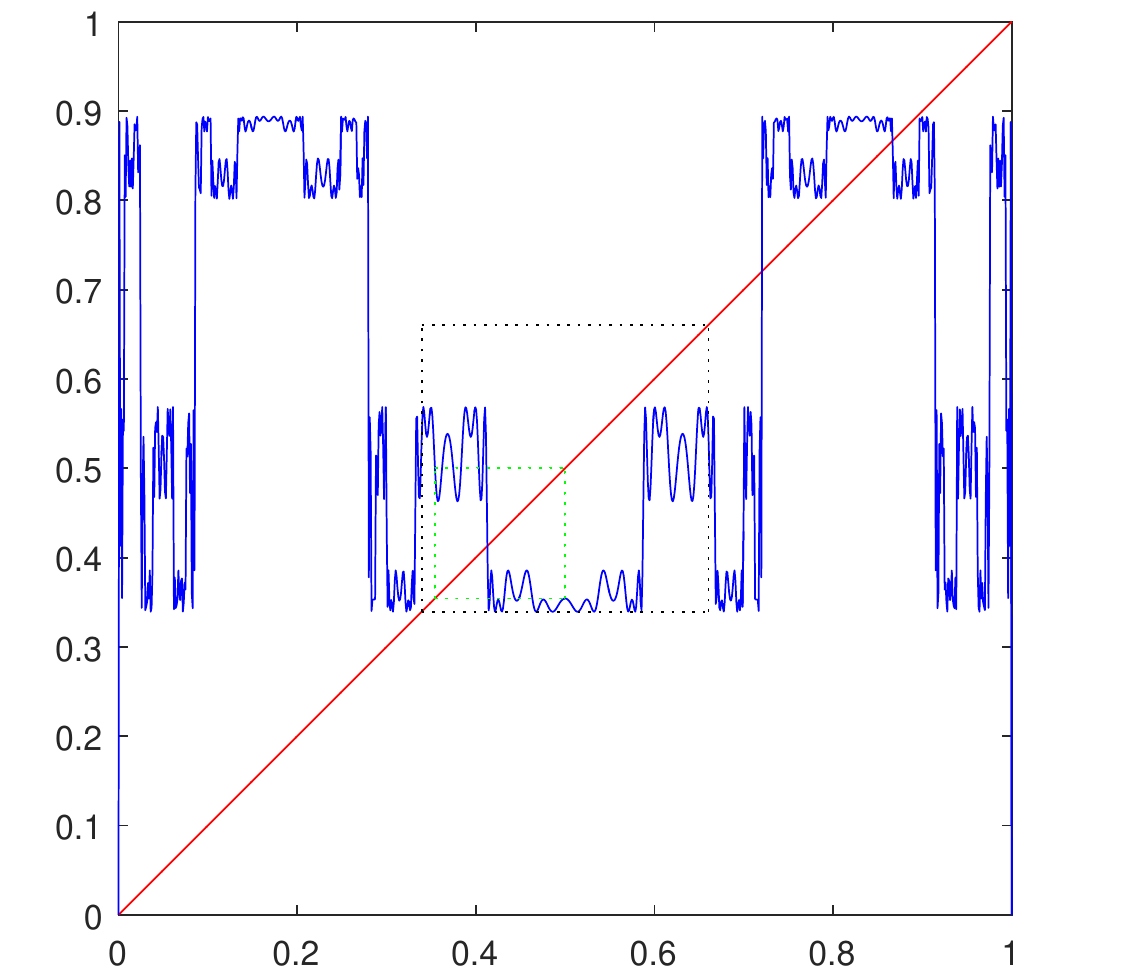}}	
	\subcaptionbox{$f^{36}$ with superstable $2$}{\includegraphics[width=\figwidth\textwidth, keepaspectratio=true]{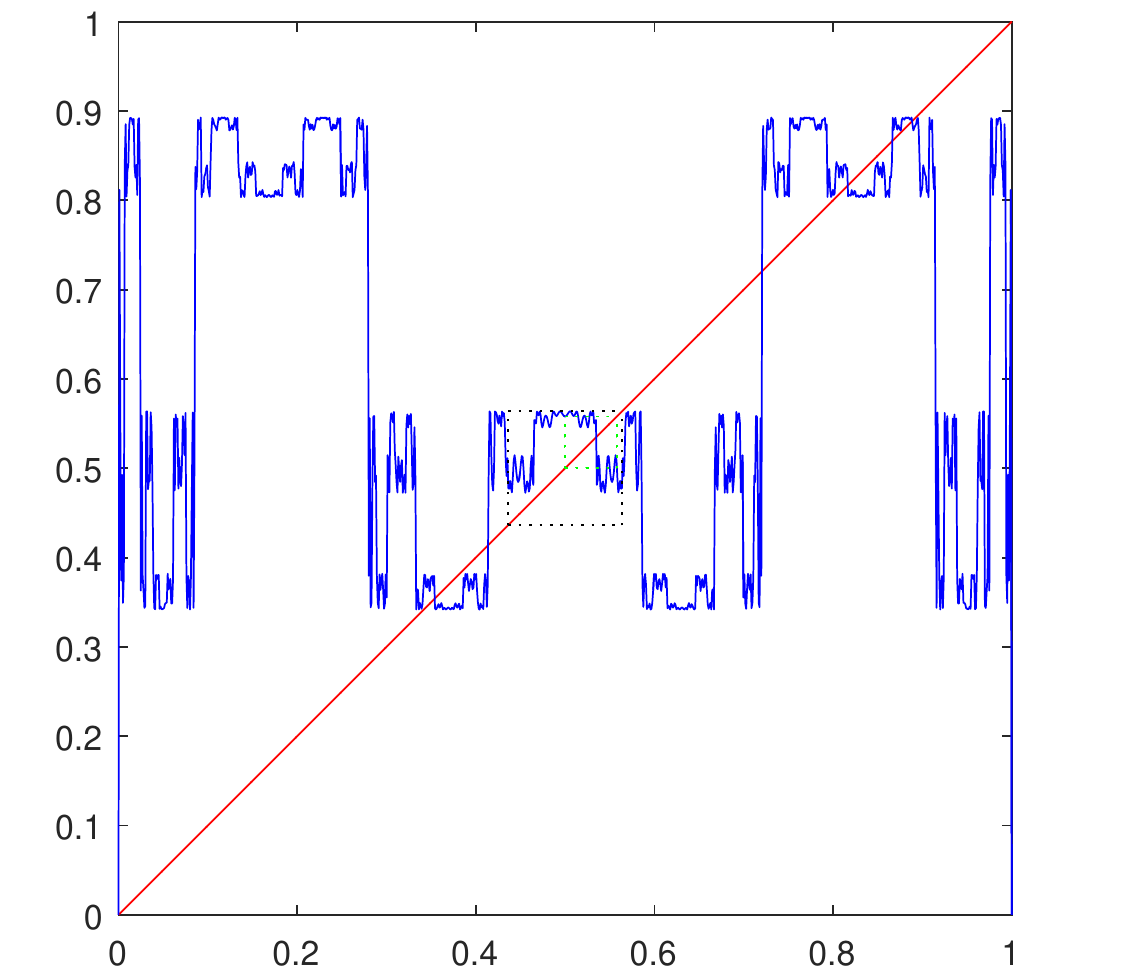}}
	\caption{The Period Doubling Mechanism for $9_1$}
	\label{fig:prdDblMech91}
\end{figure}

\begin{figure}[htb]
	\centering		
	\subcaptionbox{$f^{9}$ with superstable $2$}{\includegraphics[width=\figwidth\textwidth, keepaspectratio=true]{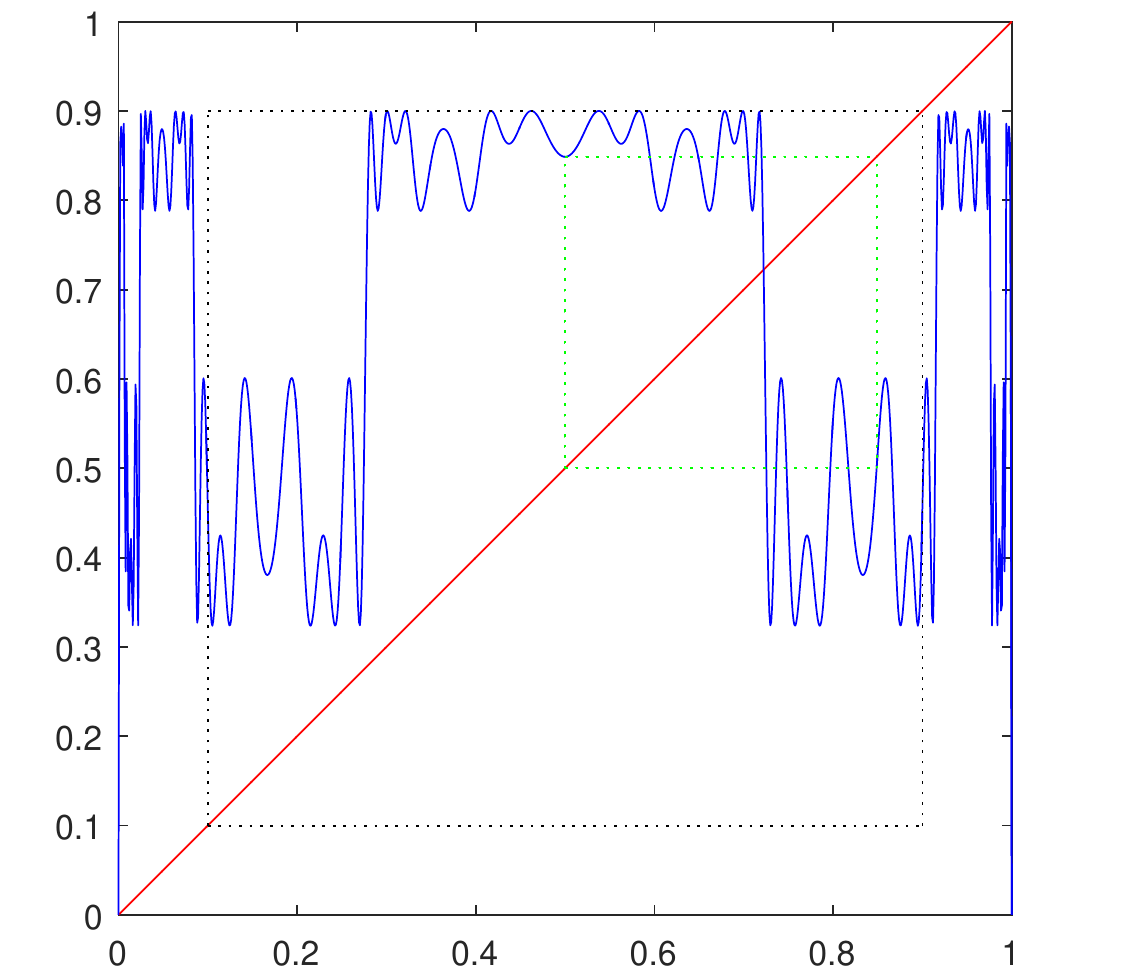}} 					
	\subcaptionbox{$f^{18}$ with superstable $2$}{\includegraphics[width=\figwidth\textwidth, keepaspectratio=true]{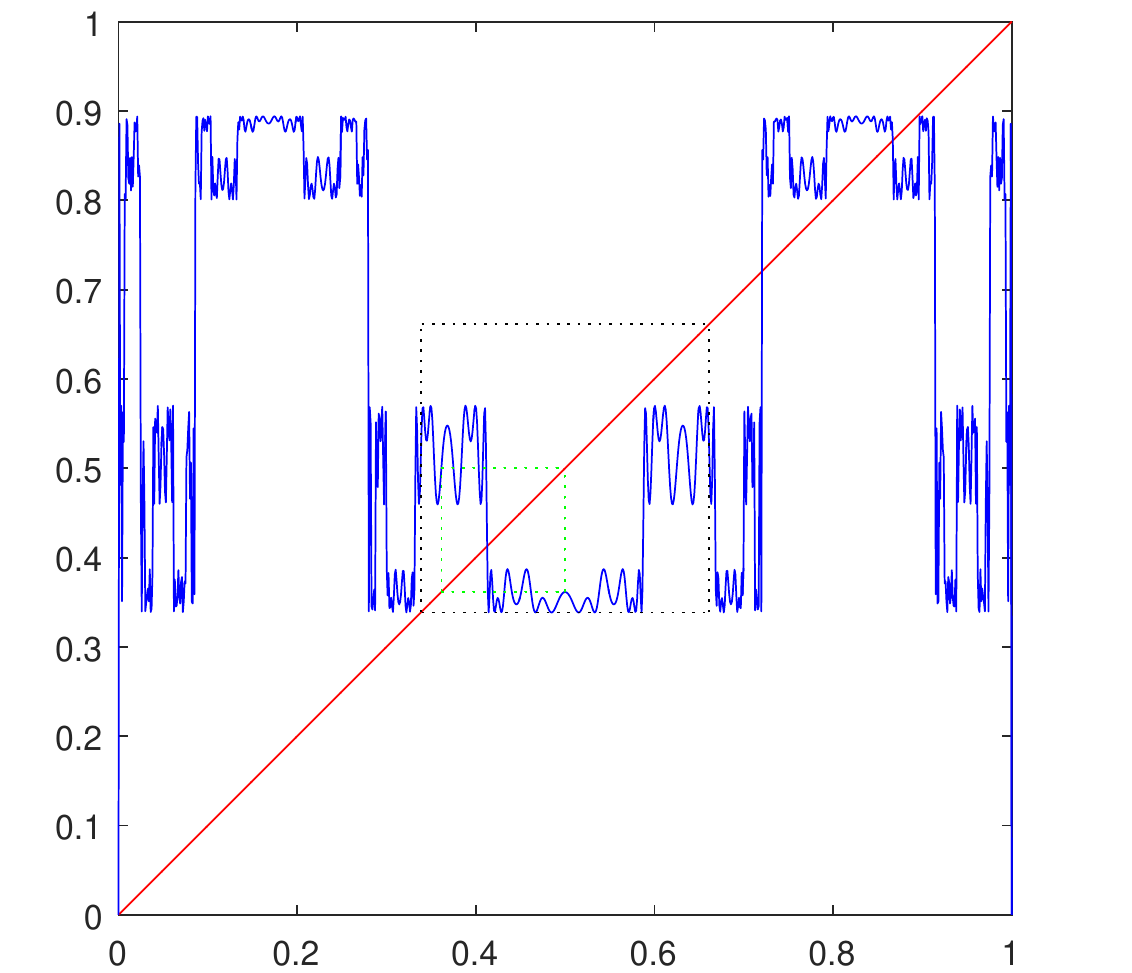}}	
	\subcaptionbox{$f^{36}$ with superstable $2$}{\includegraphics[width=\figwidth\textwidth, keepaspectratio=true]{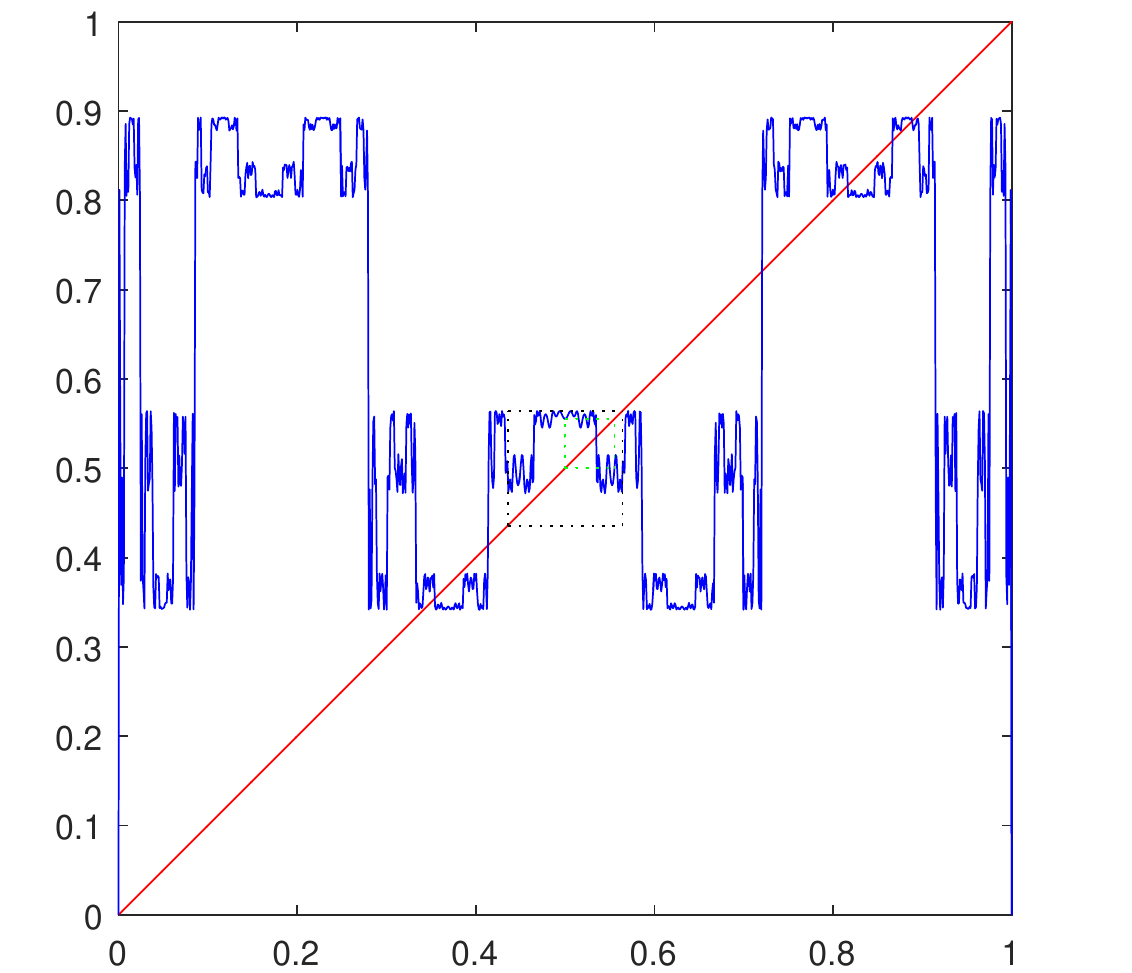}}
	\caption{The Period Doubling Mechanism for $9_2$}
	\label{fig:prdDblMech92}
\end{figure}

\begin{figure}[htb]
	\centering
	\begin{minipage}[b]{\textwidth}
		\begin{minipage}[b]{.45\textwidth}
				\includegraphics[scale=.65]{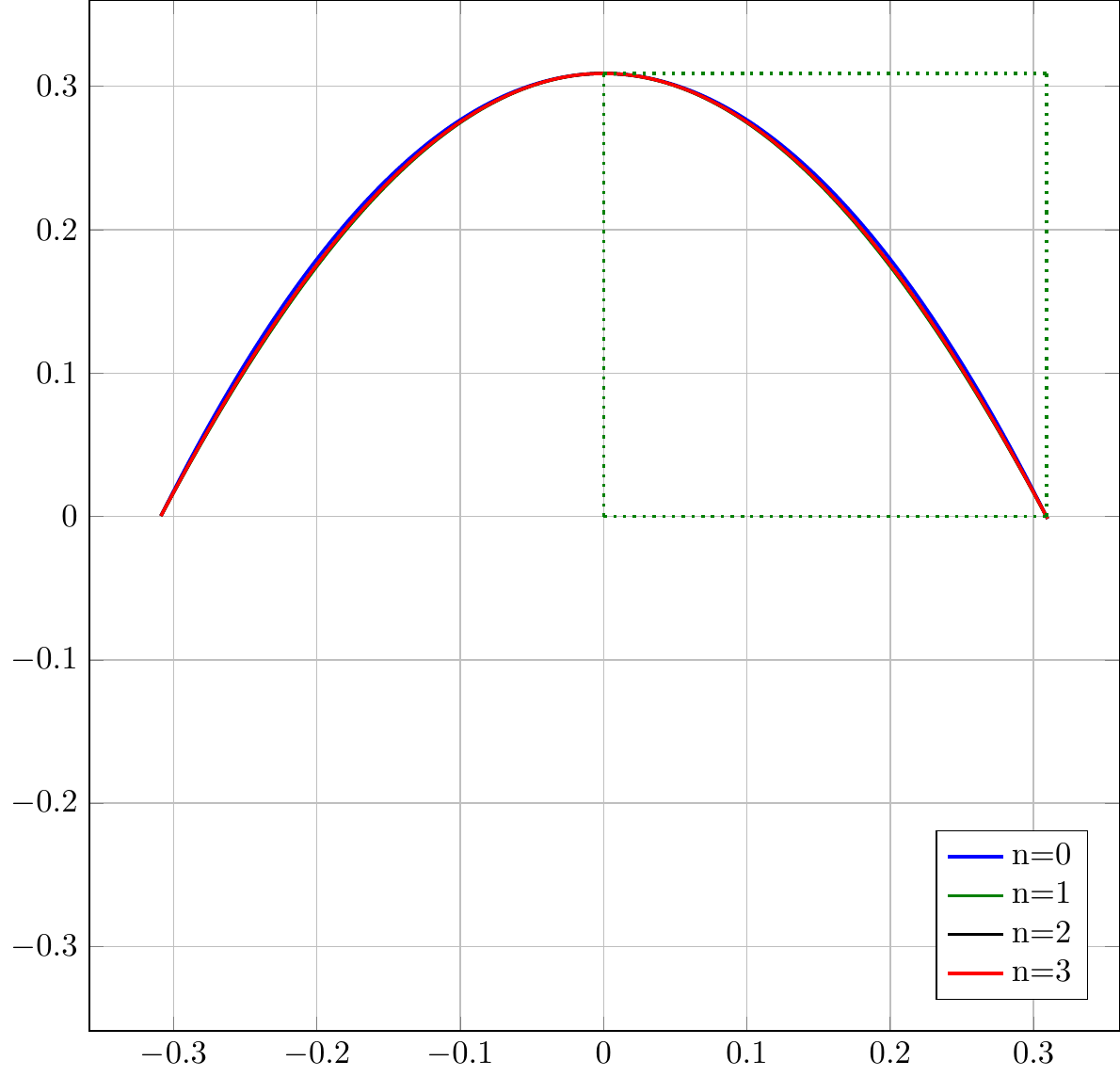}
		\subcaption{Universal function $g_{1}$ for Period $2^n$,\\ $g_{1} = \lim\limits_{n\to\infty}\left ( -\alpha \right )^{n}f^{2^{n}}_{\lambda_{n+1}}\left 
		( \frac{x}{\left ( -\alpha \right )^{n}} \right )$}
		\end{minipage}
		\hfill
		\begin{minipage}[b]{.45\textwidth}
				\includegraphics[scale=.65]{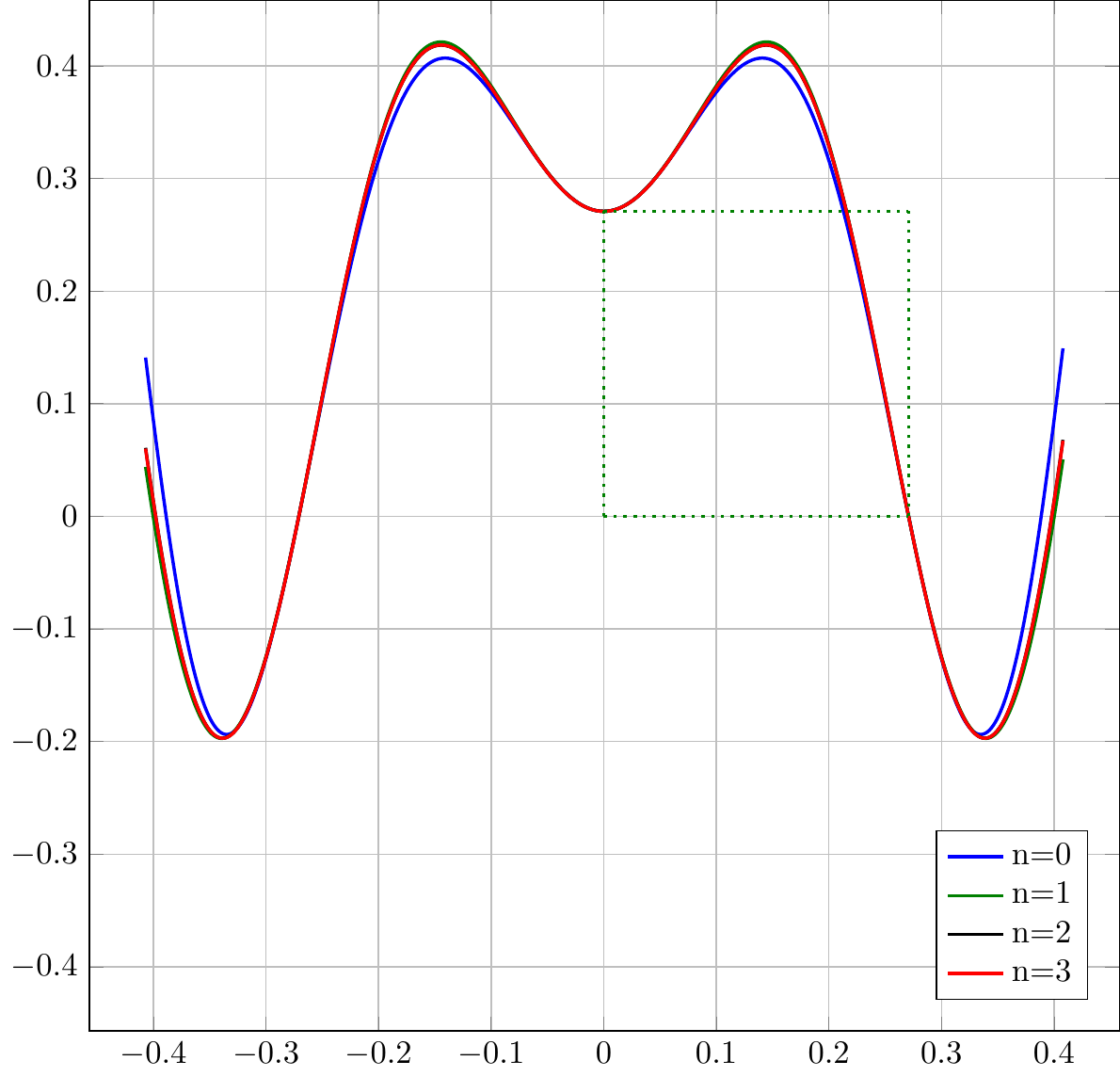}
		\subcaption{Universal function $g_{1}$ for Period $3\cdot 2^n$,\\ $g_{1} = \lim\limits_{n\to\infty}\left ( -\alpha \right )^{n}f^{3_{1}\cdot 2^{n}}_{\lambda_{n+1}}\left 
		( \frac{x}{\left ( -\alpha \right )^{n}} \right )$}
		\end{minipage}
	\end{minipage}	

	\vskip 2em

	\begin{minipage}[b]{\textwidth}
		\begin{minipage}[b]{.45\textwidth}
				\includegraphics[scale=.65]{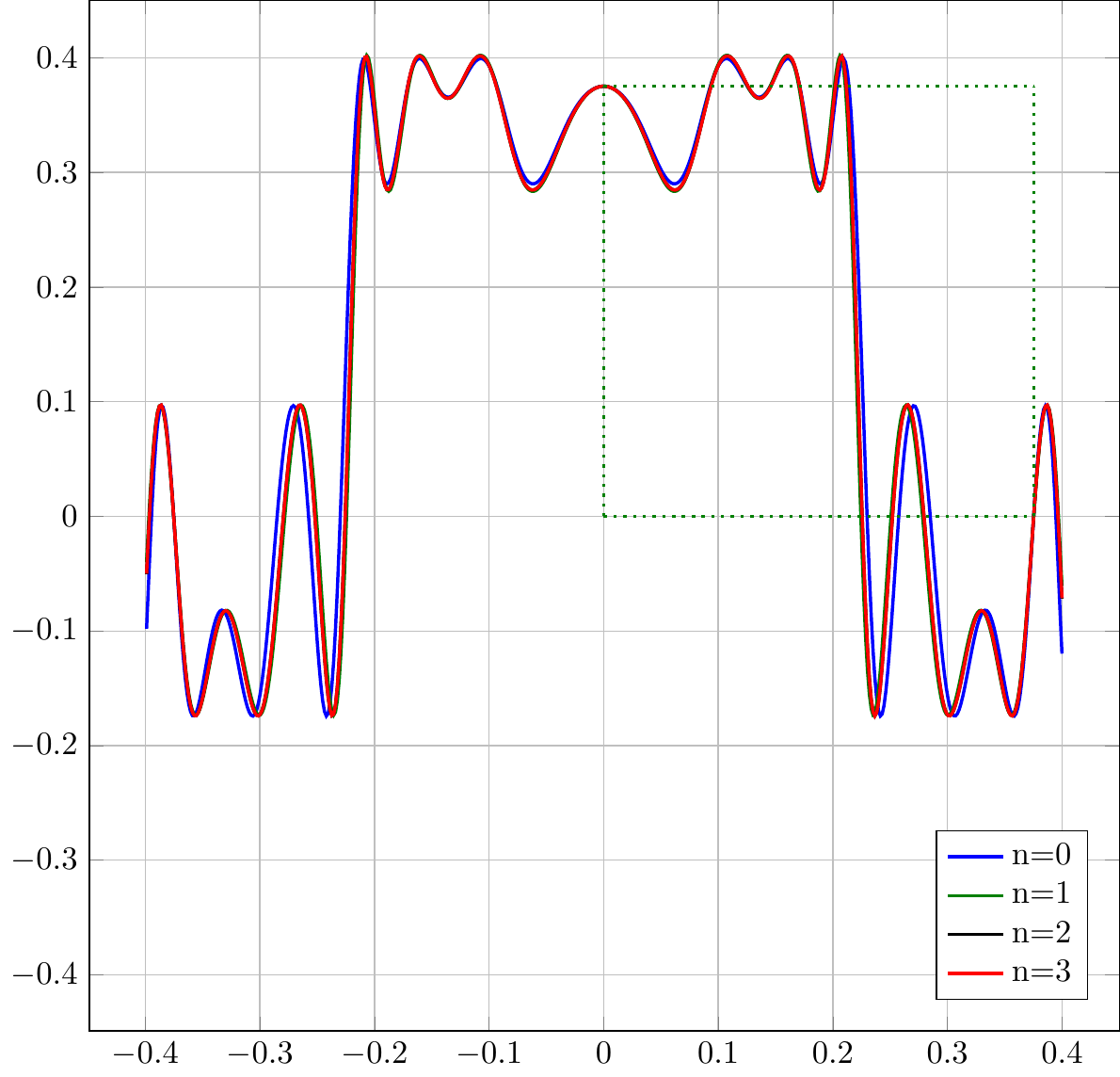}
		\subcaption{Universal function $g_{1}$ for Period $7_1\cdot 2^n$,\\ $g_{1} = \lim\limits_{n\to\infty}\left ( -\alpha \right )^{n}f^{7_{1}\cdot 2^{n}}_{\lambda_{n+1}}\left 
		( \frac{x}{\left ( -\alpha \right )^{n}} \right )$}
		\end{minipage}
		\hfill
		\begin{minipage}[b]{.45\textwidth}
				\includegraphics[scale=.65]{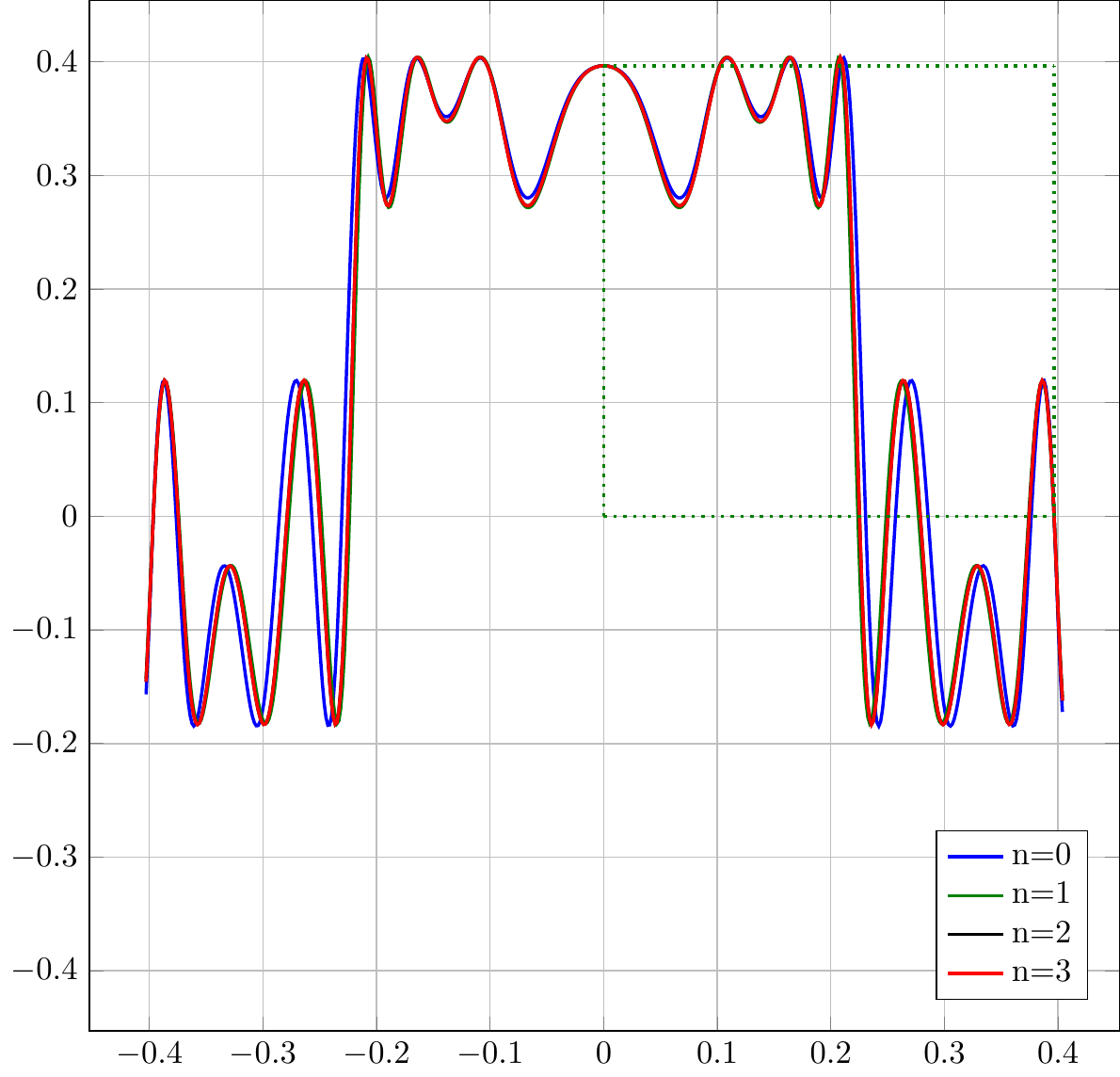}
			\subcaption{Universal function $g_{1}$ for Period $7_2\cdot 2^n$,\\ $g_{1} = \lim\limits_{n\to\infty}\left ( -\alpha \right )^{n}f^{7_{2}\cdot 2^{n}}_{\lambda_{n+1}}\left 
			( \frac{x}{\left ( -\alpha \right )^{n}} \right )$}
		\end{minipage}	
	\end{minipage}
	\caption{Universal Function $g_1$ for first appearance odds, Logistic Map.}
	\label{fig:univ}	
\end{figure}
	
\begin{figure}[htb]
	\centering
	\begin{minipage}[b]{\textwidth}
		\begin{minipage}[b]{.45\textwidth}
				\includegraphics[scale=.65]{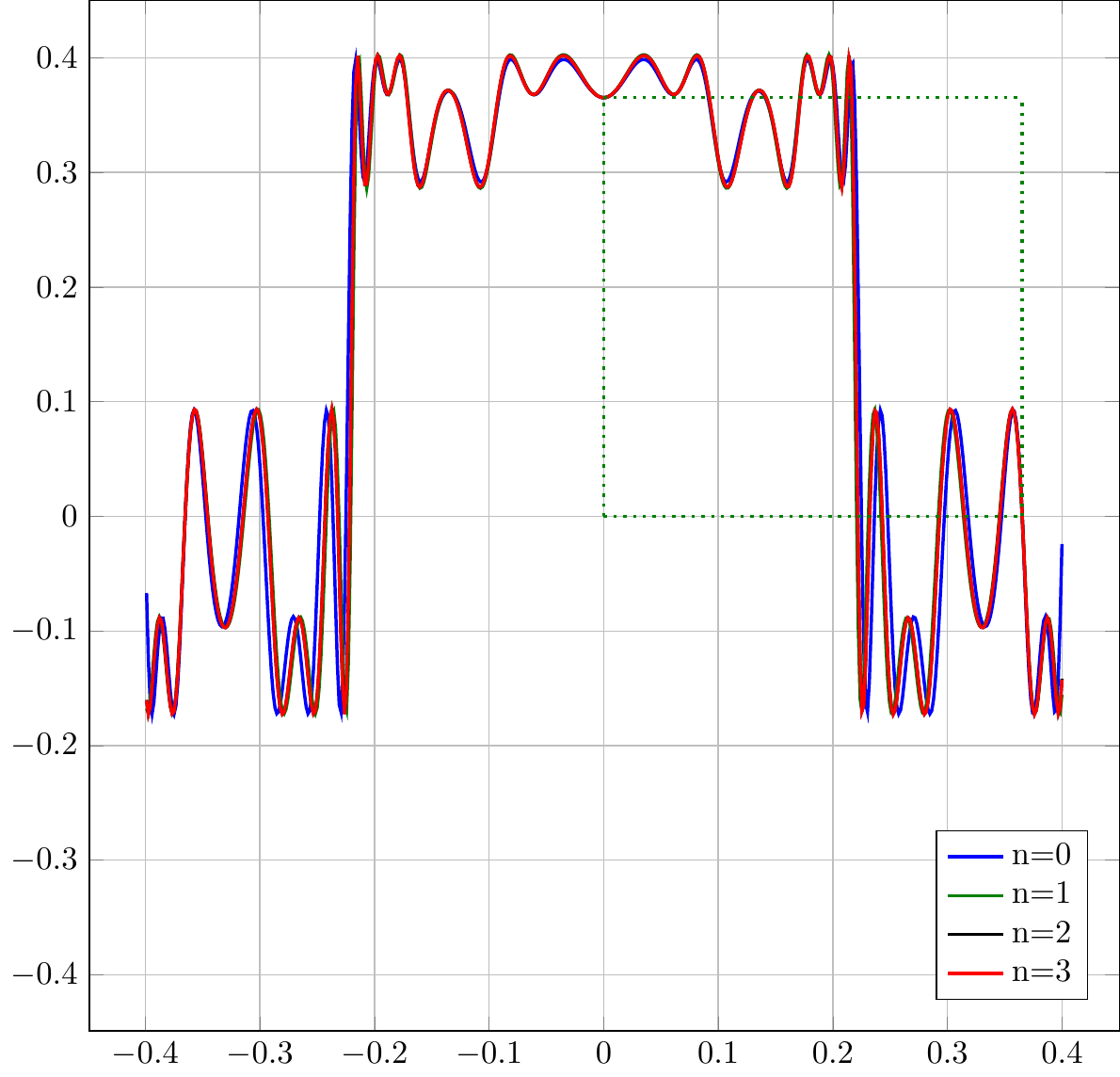}
		\subcaption{Universal function $g_{1}$ for Period $9_{1}\cdot 2^n$,\\ $g_{1} = \lim\limits_{n\to\infty}\left ( -\alpha \right )^{n}f^{9_{1}\cdot 2^{n}}_{\lambda_{n+1}}\left 
		( \frac{x}{\left ( -\alpha \right )^{n}} \right )$}
		\end{minipage}
		\hfill
		\begin{minipage}[b]{.45\textwidth}
				\includegraphics[scale=.65]{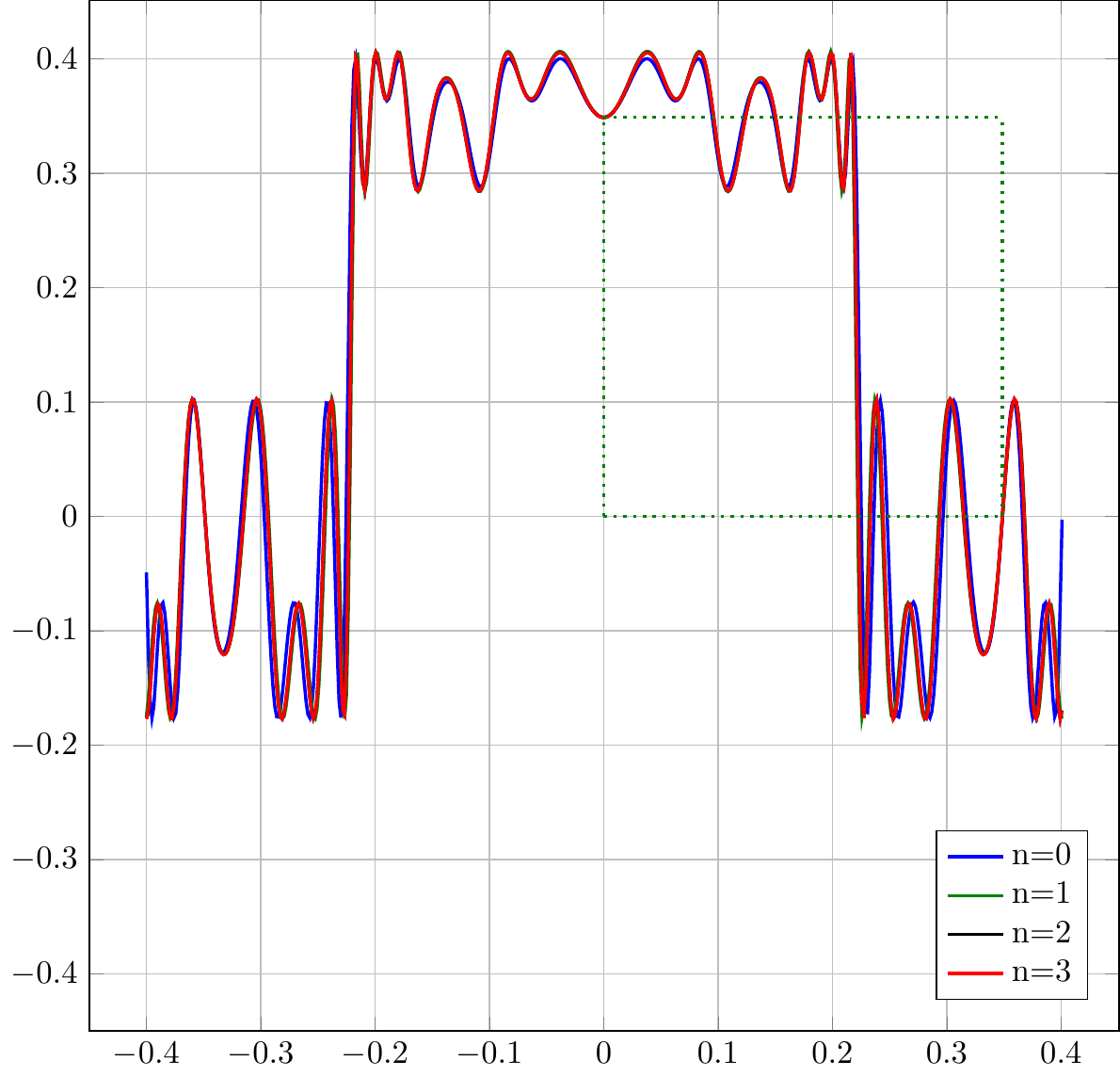}
		\subcaption{Universal function $g_{1}$ for Period $9_{2}\cdot 2^n$,\\ $g_{1} = \lim\limits_{n\to\infty}\left ( -\alpha \right )^{n}f^{9_{2}\cdot 2^{n}}_{\lambda_{n+1}}\left 
		( \frac{x}{\left ( -\alpha \right )^{n}} \right )$}
		\end{minipage}
	\end{minipage}	

	\vskip 2em

	\begin{minipage}[b]{\textwidth}
		\begin{minipage}[b]{.45\textwidth}
				\includegraphics[scale=.65]{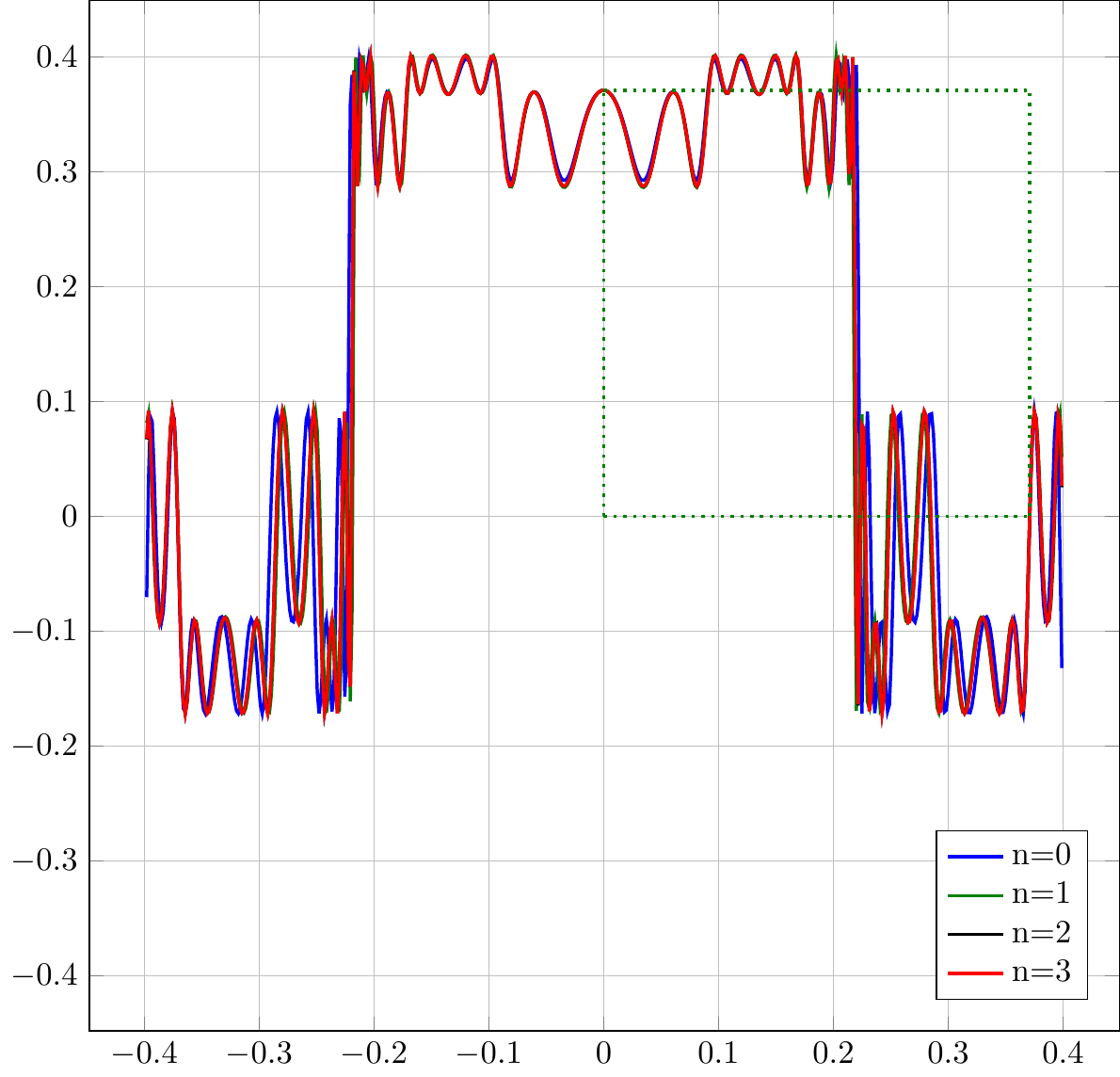}
		\subcaption{Universal function $g_{1}$ for Period $11_1\cdot 2^n$,\\ $g_{1} = \lim\limits_{n\to\infty}\left ( -\alpha \right )^{n}f^{11_{1}\cdot 2^{n}}_{\lambda_{n+1}}\left 
		( \frac{x}{\left ( -\alpha \right )^{n}} \right )$}
		\end{minipage}
		\hfill
		\begin{minipage}[b]{.45\textwidth}
				\includegraphics[scale=.65]{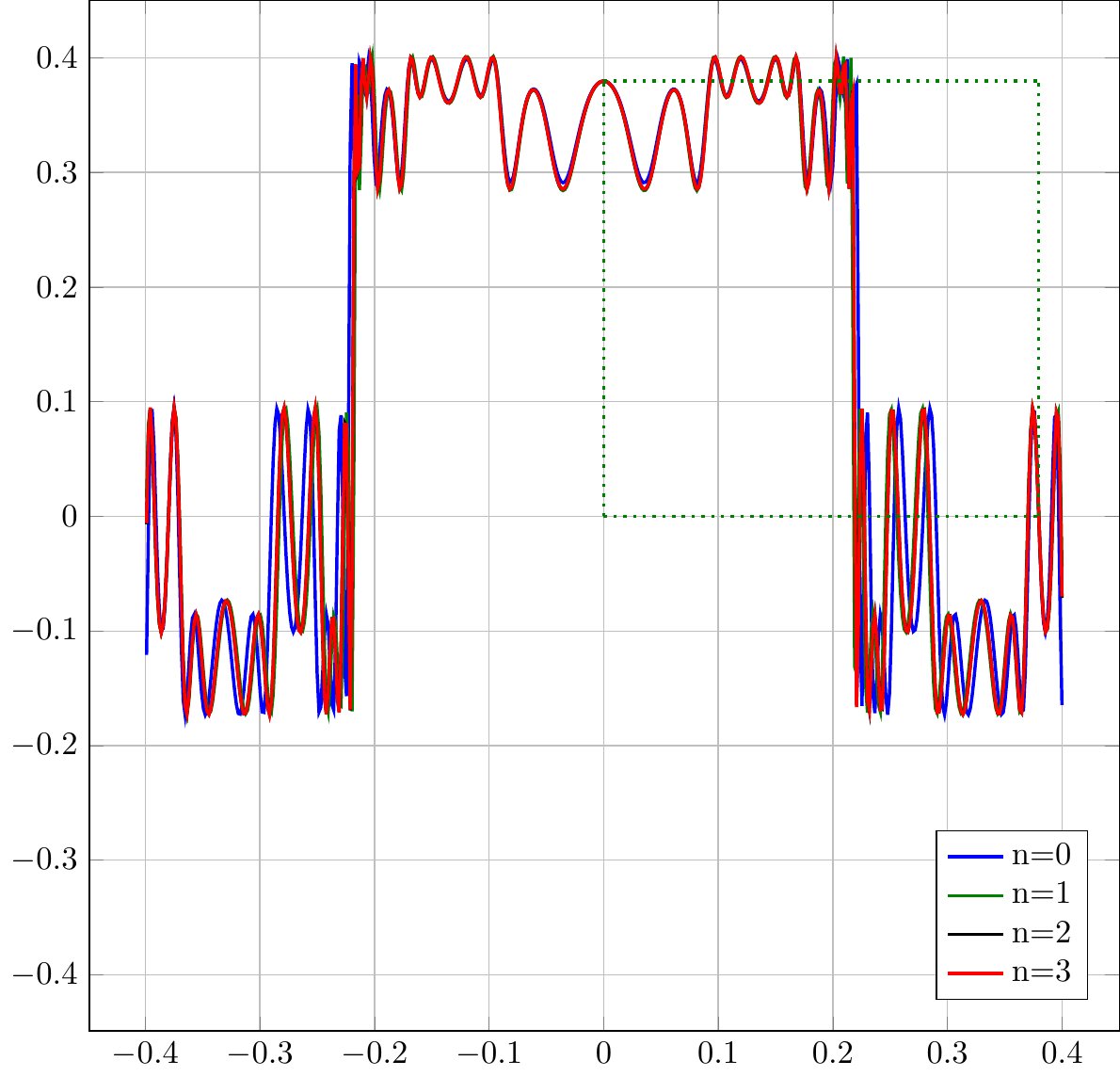}
			\subcaption{Universal function $g_{1}$ for Period $11_2\cdot 2^n$,\\ $g_{1} = \lim\limits_{n\to\infty}\left ( -\alpha \right )^{n}f^{11_{2}\cdot 2^{n}}_{\lambda_{n+1}}\left 
			( \frac{x}{\left ( -\alpha \right )^{n}} \right )$}
		\end{minipage}	
	\end{minipage}
	\caption{Universal Function $g_1$ for first appearance odds, Logistic Map.}
	\label{fig:univ2}	
\end{figure}

\begin{figure}[htb]
	\centering
	\begin{minipage}[b]{\textwidth}
		\begin{minipage}[b]{.45\textwidth}
				\includegraphics[scale=.65]{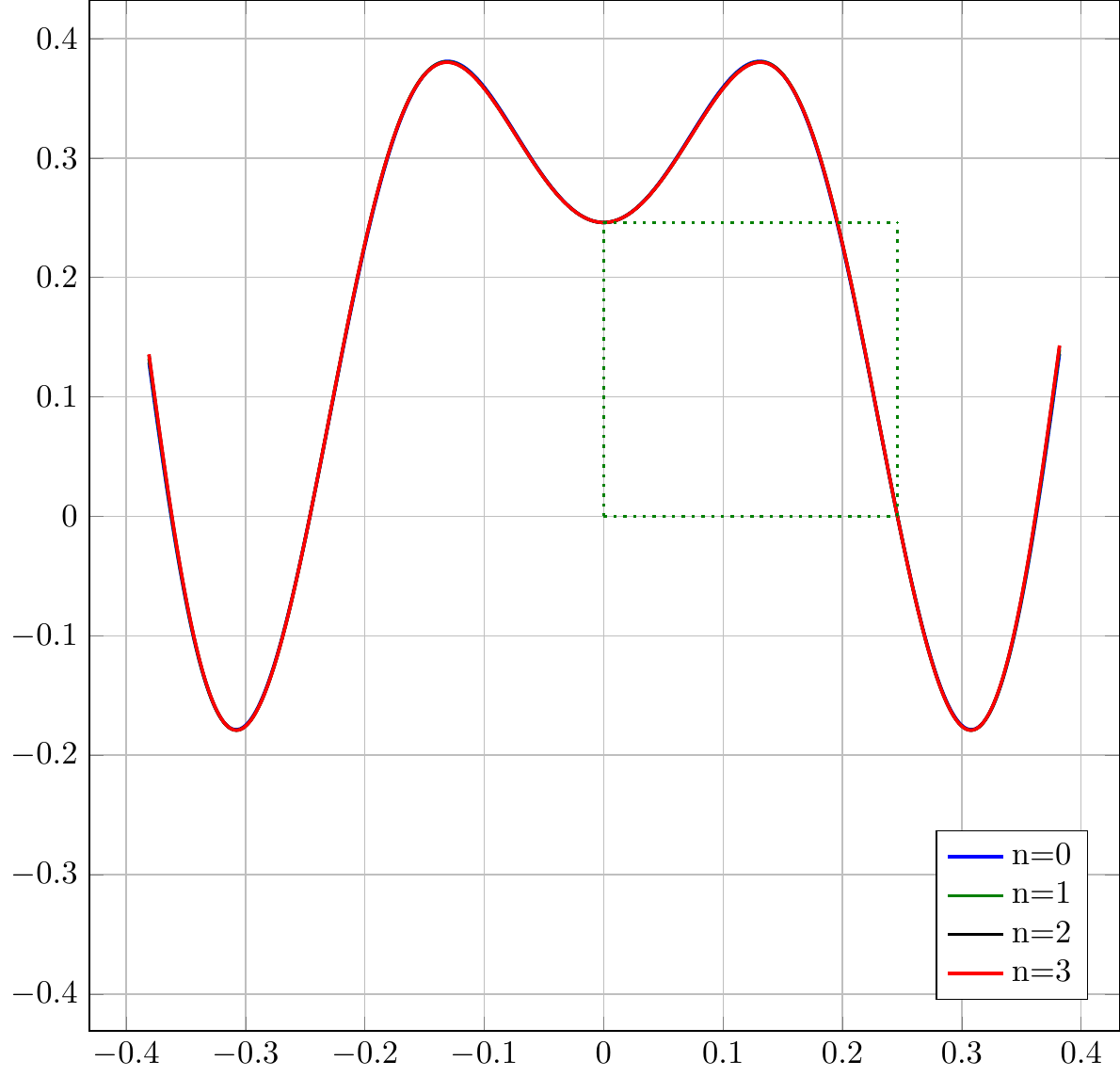}
		\subcaption{Universal function $g_{1}$ for Period $3_{1}\cdot 2^n$,\\ $g_{1} = \lim\limits_{n\to\infty}\left ( -\alpha \right )^{n}f^{3_{1}\cdot 2^{n}}_{\lambda_{n+1}}\left 
		( \frac{x}{\left ( -\alpha \right )^{n}} \right )$, Sine map.}
		\end{minipage}
		\hfill
		\begin{minipage}[b]{.45\textwidth}
				\includegraphics[scale=.65]{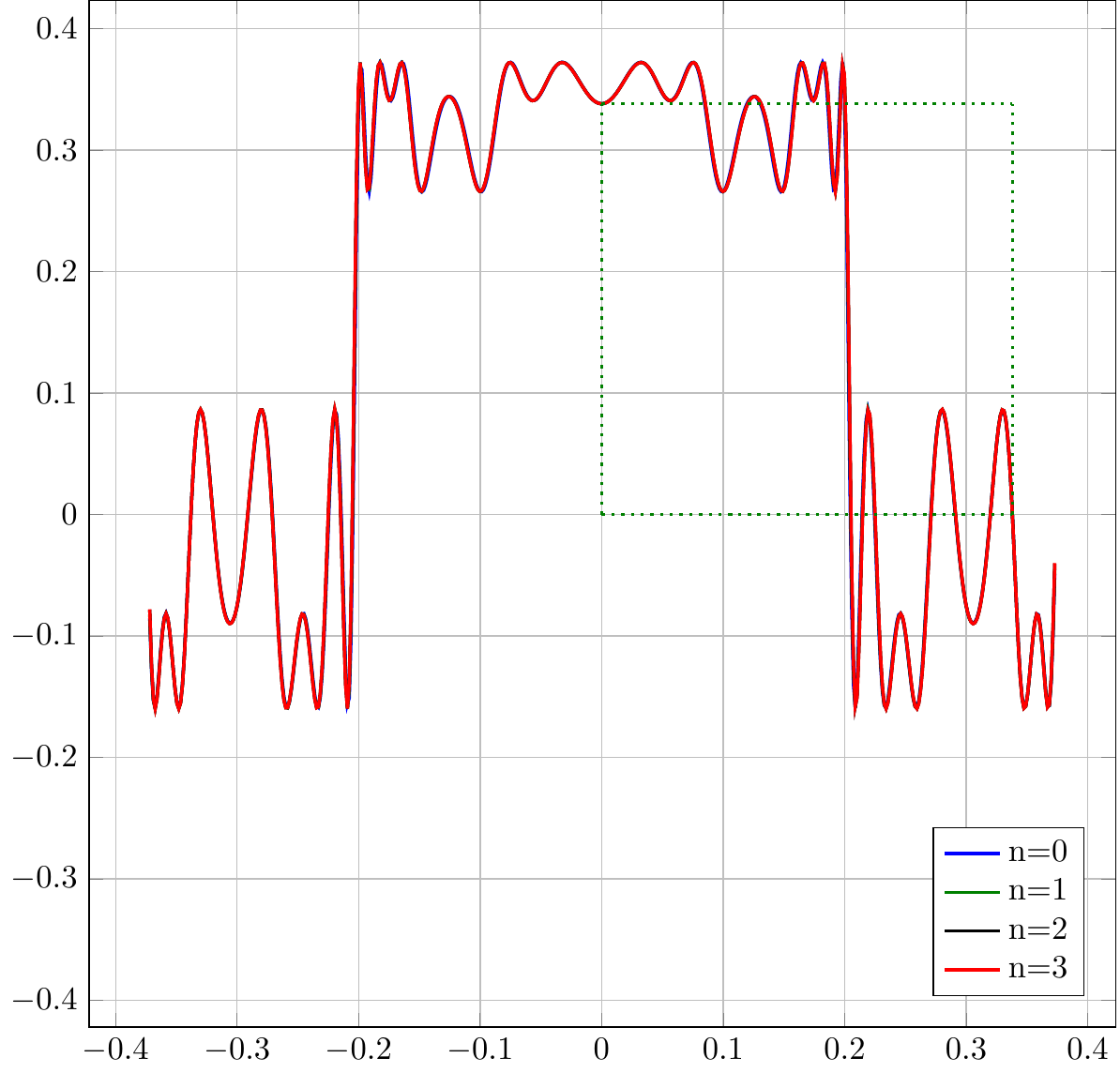}
		\subcaption{Universal function $g_{1}$ for Period $9_{1}\cdot 2^n$,\\ $g_{1} = \lim\limits_{n\to\infty}\left ( -\alpha \right )^{n}f^{9_{1}\cdot 2^{n}}_{\lambda_{n+1}}\left 
		( \frac{x}{\left ( -\alpha \right )^{n}} \right )$, Sine map.}
		\end{minipage}
	\end{minipage}	

	\vskip 2em

	\begin{minipage}[b]{\textwidth}
		\begin{minipage}[b]{.45\textwidth}
				\includegraphics[scale=.65]{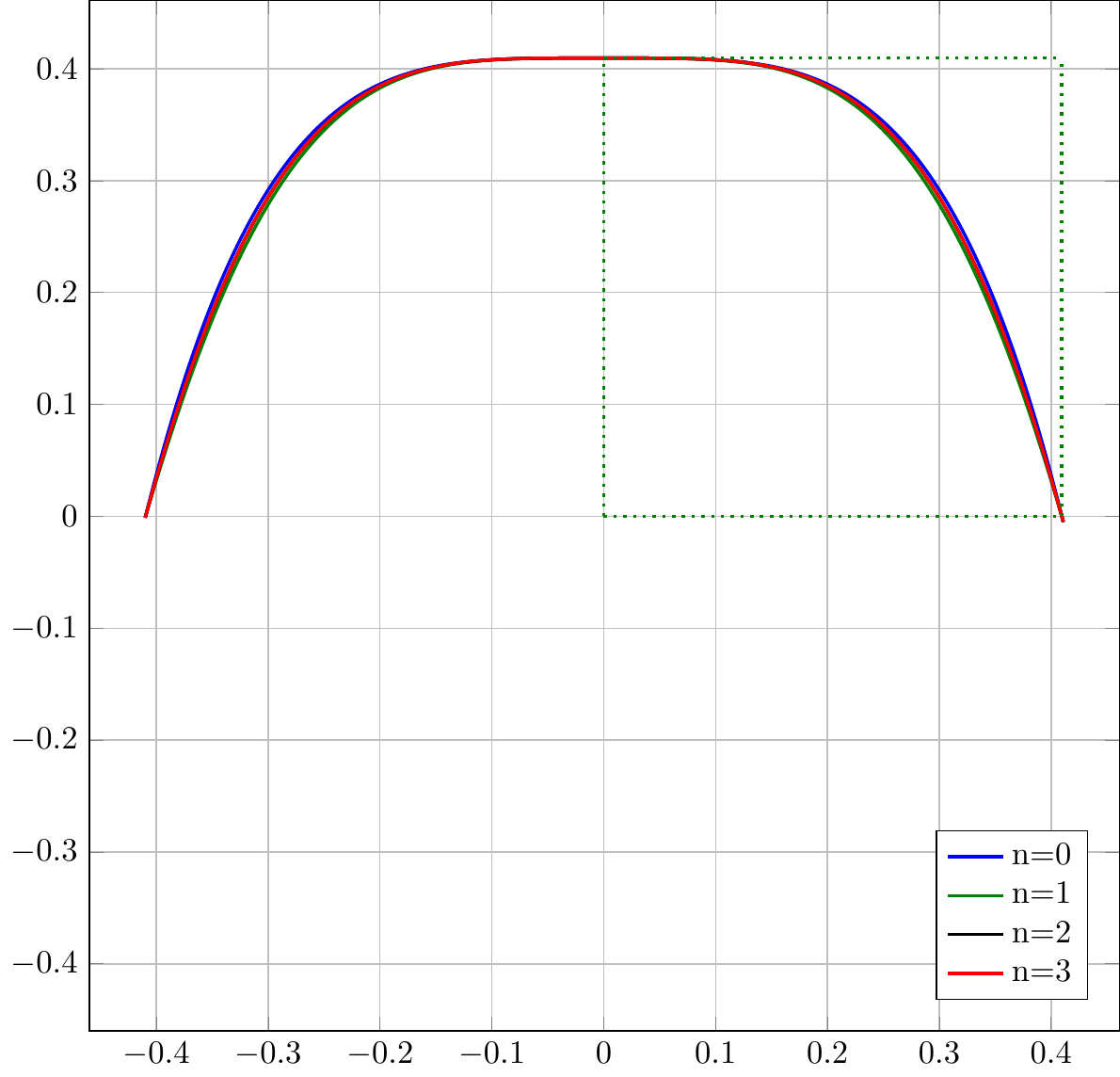}
		\subcaption{Universal function $g_{1}$ for Period $2^n$,\\ $g_{1} = \lim\limits_{n\to\infty}\left ( -\alpha \right )^{n}f^{2^{n}}_{\lambda_{n+1}}\left 
		( \frac{x}{\left ( -\alpha \right )^{n}} \right )$, Quartic map.}
		\end{minipage}
		\hfill
		\begin{minipage}[b]{.45\textwidth}
				\includegraphics[scale=.65]{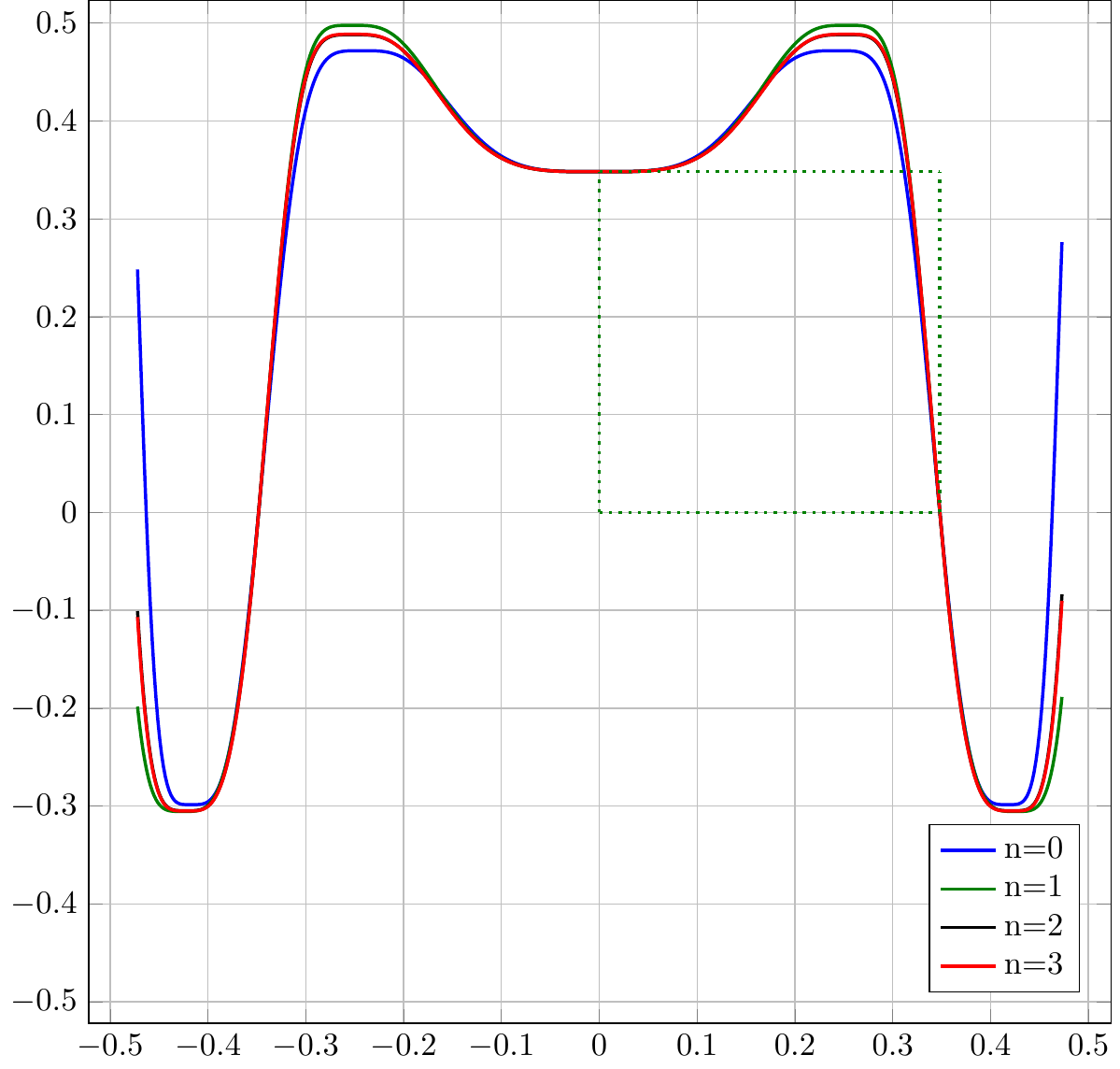}
			\subcaption{Universal function $g_{1}$ for Period $3_1\cdot 2^n$,\\ $g_{1} = \lim\limits_{n\to\infty}\left ( -\alpha \right )^{n}f^{3_{1}\cdot 2^{n}}_{\lambda_{n+1}}\left 
			( \frac{x}{\left ( -\alpha \right )^{n}} \right )$, Quartic map.}
		\end{minipage}	
	\end{minipage}
	\caption{Universal Function $g_1$ for first appearance odds, Sine and Quartic Maps.}
	\label{fig:univ41}	
\end{figure}
\FloatBarrier
	
\section{Conclusions}
\label{sec:conclusion}

The following are the main conclusions of this paper:
\begin{itemize}
\item We indroduced the notion of a second minimal orbit with respect to the Sharkovski ordering, for continuous endomorphisms on the real line. It is proved that there are 9-types of second minimal orbits up to their inverses. It is conjectured that there are $4k-3$-types of second minimal $(2k+1)$-orbits, with accuracy up to their inverses. The proof of this conjecture is addressed in a forthcoming paper.
\item We demonstrate the numerical results which reveal a fascinating universal pattern of the distribution of periodic orbits within the chaotic regime of the bifurcation diagram of the one-parameter family of unimodal maps, when the parameter changes in the range between the Feigenbaum transition point to chaos and the value when the superstable 3-orbit appears for the first time. Numerical results demonstarte that this parameter range is divided into infinitely many Sharkovski $s$-blocks where all the $2^s(2k+1)$-orbits are distributed and the pattern is independent of $s$.
\item The first appearance of any orbit in the indicated parameter range is always a minimal orbit \cite{abdulla2013}. Numerical results of this paper demonstrate that the second appearances of all odd orbits are always second minimal orbits with a Type 1 digraph. The reason for the relevance of exactly Type 1 second minimal (2k+1)-orbits are hidden in the fact that the topological structure of the single maximum unimodal map is equivalent to the topological structure of the piecewise monotonic map associated with Type 1 second minimal (2k+1)-orbits.
\item Numerical results demonstrate that the convergence of the successive parameter values for superstable $2^s(2k+1)$-orbits within each $s$-block is exponential with a rate independent of the appearance index. In particular, for any fixed two appearance indices, the ratio of distances of parameter values for respective appearances of superstable $2^s(2k+1)$-orbits is asymptotically constant for large $k$. Otherwise speaking, there is an asymptotically constant shift in appearances. 
\item Numerical results demonstrate that any superstable odd orbits in the indicated parameter range are going through successful period doublings, according to the Feigenbaum scenario when the parameter decreases to the critical transition point. In particular, this reveals that the Feigenbaum universality is true in very general classes of maps, such as the class of maps which are the $(2k+1)$st iteration of the class of $\mathscr{C}^1$-unimodal maps. This generalization is a driving force of infinitely many Feigenbaum scenarios of transition to chaos through successive bifurcations of all possible odd orbits in the indicated range when the parameter decreases towards the first transition value to chaos.
\item This paper outlines the elements of the rigorous Feigenbaum universality theory in the general class of maps, which are the $(2k+1)$st iteration of the class of 
$\mathscr{C}^1$-unimodal maps.    
\end{itemize} 
\section*{Acknowledgement}
This research was funded by National Science Foundation: grant \#1359074--REU Site: Partial Differential Equations and Dynamical Systems at Florida Institute of Technology (Principal Investigator Ugur G. Abdulla). 
Students Andy Ruden, Batul Kanawati (REU-2014), Alyssa Turnquist (REU-2015) and Emily Ribando-Gros (REU-2016) participated and contributed to numerical calculation of the superstable periodic orbits of the four model examples reported in the paper.


	\clearpage
\appendix
\section{\small{Topological Structure and Digraphs of Second Minimal $7$ Periodic Orbits}}\label{app:sevenDigraph}
\FloatBarrier
\captionsetup[subfigure]{skip=-3pt}
\begin{figure}[!htbp]		
	\centering
	\subcaptionbox{Topological Structure $1$: max}{\resizebox{0.22\textwidth}{!}{
}} 					
	\subcaptionbox{Digraph}{\resizebox{0.23\textwidth}{!}{%
}} 		
	\subcaptionbox{Topological Structure $2$: max-min-max}{\resizebox{0.22\textwidth}{!}{%
}} 					
	\subcaptionbox{Digraph}{\resizebox{0.23\textwidth}{!}{%
}} 				
\end{figure}		
\begin{figure}[!htbp]
	\centering
	\subcaptionbox{Topological Structure $3$: min-max}{\resizebox{0.22\textwidth}{!}{%
}} 					
	\subcaptionbox{Digraph}{\resizebox{0.23\textwidth}{!}{%
}}
	\subcaptionbox{Topological Structure $4$: min-max-min}{\resizebox{0.22\textwidth}{!}{%
}} 				
\end{figure}			
\begin{figure}[!htbp]
	\centering
	\subcaptionbox{Topological Structure $5$: max-min-max}{\resizebox{0.22\textwidth}{!}{%
}} 		
	\subcaptionbox{Topological Structure $6$: max-min}{\resizebox{0.22\textwidth}{!}{%
}} 		
\end{figure}			
\begin{figure}[!htbp]
	\centering
	\subcaptionbox{Topological Structure $7$: max-min}{\resizebox{0.22\textwidth}{!}{%
}} 		
	\subcaptionbox{Topological Structure $8$: max-min-max}{\resizebox{0.22\textwidth}{!}{%
}} 		
\end{figure}			
\begin{figure}[!htbp]
	\centering
	\subcaptionbox{Topological Structure $9$: max-min-max-min-max}{\resizebox{0.22\textwidth}{!}{%
}} 		
\end{figure}	
	
\FloatBarrier
\section{Period Doubling Universality}\label{app:feiguni}
\FloatBarrier
	\begin{longtable}[!htbp]{c | c | l | l| l}
	\caption{Logistic Map, Calculation of $\delta$, $\alpha$ and $\lambda_{\infty}$ for Period Doubling Starting at $(2k+1)_j$-orbit.}\label{tbl:logFeig}\\
	\toprule
		$s$ & $(2^{s}(2k+1))_{j}$ & $\delta=\frac{\lambda^{2k+1}_{s-1,j}-\lambda^{2k+1}_{s-2,j}}{\lambda^{2k+1}_{s,i}-\lambda^{2k+1}_{s-1,j}}$ & $\frac{d^{2k+1}_{s-1,j}}{d^{2k+1}_{s,j}}$ & $\lambda_{\infty} = \frac{\lambda_{s,j}^{2k+1}-\lambda_{s-1,j}^{2k+1}}{\delta} + \lambda_{s,j}^{2k+1}$ \\ \midrule
		$0$  &  $3_1$  &  $$  &  $$  & $$\\\midrule
		$1$  &  $6_1$  &  $$  &  $-2.454268432041252$  & $0.89594979661707406416946181313487$\\\midrule
		$2$  &  $12_1$ &  $4.507542941$  &  $-2.488688613626316$  & $0.89313050797449905484309938153867$\\\midrule
		$3$  &  $24_1$ &  $4.695932444$  &  $-2.499742045692276$  & $0.89262750317573161723282801462988$\\\midrule
		$4$  &  $48_1$ &  $4.667366742$ & $-2.502259346885118$ & $0.89251658164785687538727952773507$\\\midrule

		$0$  &  $5_1$   &  $$  &  $$  & $$\\
\bottomrule	
\end{longtable}


	%
	{\begin{longtable}{c | c | l | l | l}
	\caption{Sine Map, Calculation of $\delta$, $\alpha$ and $\lambda_{\infty}$ for Period Doubling Starting at $(2k+1)_j$-orbit.}\label{tbl:sineFeig}\\
	\toprule
		$s$ & $(2^{s}(2k+1))_{j}$ & $\delta=\frac{\lambda^{2k+1}_{s-1,j}-\lambda^{2k+1}_{s-2,j}}{\lambda^{2k+1}_{s,i}-\lambda^{2k+1}_{s-1,j}}$ & $\frac{d^{2k+1}_{s-1,j}}{d^{2k+1}_{s,j}}$ & $\lambda_{\infty} = \frac{\lambda_{s,j}^{2k+1}-\lambda_{s-1,j}^{2k+1}}{\delta} + \lambda_{s,j}^{2k+1}$ \\ \midrule

$0$ & $3_1$ & $$ & $$ & $$\\\midrule
$1$ & $6_1$ & $$ & $-2.458609821276883$ & $0.86873961050980452273981931736855$\\\midrule
$2$ & $12_1$ & $4.730427896$ & $-2.492752378392625$ & $0.86627858717481051811520136603587$\\\midrule
$3$ & $24_1$ & $4.689585511$ & $-2.500535459353473$ & $0.86573097100235693427121291937719$\\\midrule
$4$ & $48_1$ & $4.673272858$ & $-2.502381199234309$ & $0.86561183926284564115177937682167$\\\midrule

$0$ & $5_1$ & $$ & $$ & $$\\
\bottomrule

	\end{longtable}}
	%

	%
	{\begin{longtable}{c | c | l | l | l}
	\caption{Cubic Map, Calculation of $\delta$, $\alpha$ and $\lambda_{\infty}$ for Period Doubling Starting at $(2k+1)_j$-orbit.}\label{tbl:cubFeig}\\
	\toprule
		$s$ & $(2^{s}(2k+1))_{j}$ & $\delta=\frac{\lambda^{2k+1}_{s-1,j}-\lambda^{2k+1}_{s-2,j}}{\lambda^{2k+1}_{s,i}-\lambda^{2k+1}_{s-1,j}}$ & $\frac{d^{2k+1}_{s-1,j}}{d^{2k+1}_{s,j}}$ & $\lambda_{\infty} = \frac{\lambda_{s,j}^{2k+1}-\lambda_{s-1,j}^{2k+1}}{\delta} + \lambda_{s,j}^{2k+1}$ \\ \midrule
	
$0$ & $7_2$ & $$ & $$ & $$\\\midrule
$1$ & $14_2$ & $$ & $-2.372145654798267$ & $0.88804789214496225396997450138614$\\\midrule
$2$ & $28_2$ & $4.723797089$ & $-2.531319536769563$ & $0.88657074923765426350764650677041$\\\midrule
$3$ & $56_2$ & $4.679053349$ & $-2.485852964973201$ & $0.88624011561640477208179462882797$\\\midrule
$4$ & $112_2$ & $4.671403318$ & $-2.508524126657965$ & $0.8861687905436954155321253795437$\\\midrule

$0$ & $9_2$ & $$ & $$ & $$\\
\bottomrule
	\end{longtable}}
	%

	%
	{\begin{longtable}{c | c | l | l | l}
	\caption{Quartic Map, Calculation of $\delta$, $\alpha$ and $\lambda_{\infty}$ for Period Doubling Starting at $(2k+1)_j$-orbit.}\label{tbl:quaFeig}\\
	\toprule
		$s$ & $(2^{s}(2k+1))_{j}$ & $\delta=\frac{\lambda^{2k+1}_{s-1,j}-\lambda^{2k+1}_{s-2,j}}{\lambda^{2k+1}_{s,i}-\lambda^{2k+1}_{s-1,j}}$ & $\frac{d^{2k+1}_{s-1,j}}{d^{2k+1}_{s,j}}$ & $\lambda_{\infty} = \frac{\lambda_{s,j}^{2k+1}-\lambda_{s-1,j}^{2k+1}}{\delta} + \lambda_{s,j}^{2k+1}$ \\ \midrule
	
$0$ & $7_2$ & $$ & $$ & $$\\\midrule
$1$ & $14_2$ & $$ & $-1.685016151822741$ & $0.96795394159304505710938590181654$\\\midrule
$2$ & $28_2$ & $7.257463741$ & $-1.689616004337216$ & $0.96852661836071093093947880007475$\\\midrule
$3$ & $56_2$ & $7.286379135$ & $-1.690226809471704$ & $0.96860664476713929365715844052214$\\\midrule
$4$ & $112_2$ & $7.284942740$ & $-1.690235270778622$ & $0.96861762018986771419027828078051$\\\midrule

$0$ & $9_2$ & $$ & $$ & $$\\
\bottomrule
	\end{longtable}}

\section{Parameter Tables}
\label{app:params}

The following few pages contain parameter values that we used to construct the tables and figures in this document. Below is a key outlining the table headers.
\begin{itemize}
	\item \textbf{Parameter}: the numeric value of the parameter for the associated map
	\item \textbf{P}: the super stable periodic orbit corresponding to the above parameter value
	\item \textbf{A}: the appearance number of the periodic orbit corresponding to the parameter value
\end{itemize}

\begin{center}
\twocolumn
\topcaption{Parameter Values for Logistic Map $f_{\lambda}(x)=4\lambda x(1-x)$}
\tablefirsthead{%
\hline
Parameter & P & A \\
}
\tablehead{Parameter & P & A \\}
\tabletail{%
}
\tablelasttail{\hline}
\begin{supertabular}[H]{p{0.30\textwidth}p{0.05\textwidth}p{0.05\textwidth}}
0.5 & 1.0 & 1.0 \\ 
0.8090169943749475 & 2.0 & 1.0 \\ 
0.8746404248319252 & 4.0 & 1.0 \\ 
0.8886602156922059 & 8.0 & 1.0 \\ 
0.8916668449640671 & 16.0 & 1.0 \\ 
0.8925435410844437 & 176.0 & 1.0 \\ 
0.8925462153481506 & 144.0 & 1.0 \\ 
0.8925493727315035 & 176.0 & 2.0 \\
\end{supertabular}
\end{center}
\end{document}